\documentclass[a4paper, 10pt]{amsart}
\usepackage[mathscr]{eucal}
\usepackage{amssymb, amsfonts, amsmath, amsthm, graphicx, cases, xcolor, tikz} 

\newtheorem{proposition}{Proposition}[section]
\newtheorem{theorem}{Theorem}[section]
\newtheorem{lemma}[proposition]{Lemma}

\theoremstyle{definition}
\newtheorem{definition}{Definition}[section]

\newtheorem*{thank}{Acknowledgments}
\newtheorem*{availability}{Data Availability}
\newtheorem*{interest}{Conflict of interest}
\newtheorem*{ai}{Use of AI}
\theoremstyle{remark}
\newtheorem{remark}{Remark}[section]

\numberwithin{equation}{section}

\title{Constrained minimization problem 
and unbounded dynamics 
for the four-dimensional 
Klein-Gordon-Schr\"{o}dinger System}
\author{Masaru Hamano, Hiroaki Kikuchi 
and Kenta Tomioka}
\date{\today}
\newcommand{\R}{\mathbb{R}}
\newcommand{\N}{\mathbb{N}}

\newcommand{\e}{\varepsilon}
\newcommand{\p}{\partial}

\newcommand{\K}{\mathcal{K}}
\newcommand{\M}{\mathcal{M}}

\newcommand{\E}{\mathcal{E}}



\begin{document}
\maketitle

\tableofcontents

\begin{abstract}
In this paper, we study the Klein–Gordon–Schr\"{o}dinger system in four spatial dimensions. We consider a constrained minimization problem associated with the corresponding stationary system and establish the existence of a positive radial minimizer. We further investigate the asymptotic behavior and qualitative properties of minimizers. As an application, we analyze the dynamics of solutions 
with initial data in a potential well and prove that such solutions either blow up in finite time or exhibit unbounded growth in the energy norm as time tends to infinity. Since the minimizers lie on the boundary of the potential well, we also establish the strong instability of the associated standing waves.
\end{abstract}

\section{Introduction} 
In this paper, we consider the following 
Klein-Gordon-Schr\"{o}dinger system: 
	\begin{equation}\label{kgs}
	\begin{cases}
	i \p_{t} u + \Delta u = -2uv & \qquad \mbox{in $[0,\infty) \times \R^{d}$}, \\
	\p_{t}^{2} v - \Delta v + v = |u|^{2} 
	& \qquad \mbox{in $[0,\infty) \times \R^{d}$},  
	\end{cases}
	\end{equation}
where $1 \leq d \leq 5$, $u: [0,\infty) \times \mathbb{R}^{d} \to \mathbb{C}$, and $v: [0,\infty) \times \mathbb{R}^{d} \to \mathbb{R}$.
The system \eqref{kgs} arises as a model describing the interaction between a complex nucleon field $u$ and a real meson field 
$v$
(see e.g. \cite{MR515899, Yukawa}). 

The main purpose of this paper is to investigate the global dynamics of solutions to \eqref{kgs} in the four-dimensional case ($d=4$). 
To this end, we analyze a constrained minimization problem associated with \eqref{kgs} and establish several 
properties of its minimizers. 
As an application of these variational properties, we introduce a potential well and prove that solutions to \eqref{kgs} 
with initial data in this potential well either blow up in finite time or exhibit unbounded growth in the energy norm as time tends to infinity. 
Moreover, since the minimizers lie on the boundary of the potential well, we establish the strong instability of the associated standing waves.

We first consider the local well-posedness of 
\eqref{kgs}, which has been studied 
extensively
(see \cite{MR2289361, MR778979, MR519635, MR2403699, 
MR515899, MR2225086, MR900982, 
MR2035502, MR2906550} and references therein). 
We define the energy space $X$ as 
follows: 
    \[
    X := H^{1}(\R^{d}) \times 
    H^{1}(\R^{d}) \times L^{2}(\R^{d}). 
    \]
We have the following local well-posedness result in the energy space $X$.  
\begin{theorem}\label{thm1:local-wp}
Assume that $d = 4$. 
For any $(u_{0}, v_{0}, v_{1}) \in X$, 
there exist a maximal existence time $T_{\max} 
= T_{\max}(\|(u_0, v_0, v_{1})\|_{X})
> 0$ 
and a unique local solution 
$(u, v, \partial_{t} v) \in C([0, T_{\max}), X)$ 
with $(u, v, \partial_{t} v)|_{t = 0} = 
(u_{0}, v_{0}, v_{1})$. Moreover, 
the following conservation laws hold: 
\begin{equation}\label{eq:conserv}
\mathcal{E}(u(t), v(t), \partial_{t} v(t)) = 
\mathcal{E}(u_{0}, v_{0}, v_{1}), \quad 
\|u(t)\|_{L^{2}} = \|u_{0}\|_{L^{2}} 
\quad \mbox{for all $t \in [0, T_{\max})$},  
\end{equation}
where the energy $\mathcal{E}$ is defined by 
    \begin{equation*}
	\mathcal{E}(u, v, w) := \|\nabla u\|_{L^{2}}^{2} 
	 + \|\nabla v\|_{L^{2}}^{2} 
	+ \|v\|_{L^{2}}^{2} 
    + \|w\|_{L^{2}}^{2}
    - 2 \int_{\R^{d}} |u|^{2} v \, dx.  
	\end{equation*}
\end{theorem}
Theorem \ref{thm1:local-wp} follows from standard arguments based on Strichartz estimates. For completeness, we provide its proof in Appendix \ref{sec:cauchy} below.

We next consider standing wave solutions 
to \eqref{kgs}. 
By a \textit{standing wave}, we mean a solution 
to \eqref{kgs} of the form 
    \[
    (u(t, x), v(t, x)) = 
    (e^{i \omega t}\varphi_{\omega}(x), 
    \psi_{\omega}(x)) 
    \qquad (\omega > 0). 
    \]
Then, the pair $(\varphi_{\omega}, \psi_{\omega})$ 
satisfies the following system of 
elliptic equations: 
	\begin{equation} \label{sp}
	\begin{cases}
	- \Delta \varphi + \omega \varphi = 2 \varphi \psi 
	& \qquad \mbox{in $\R^{d}$}, \\
	- \Delta \psi + \psi = |\varphi|^{2} 
	& \qquad \mbox{in $\R^{d}$}. 
	\end{cases}
	\end{equation}
We put 
	\begin{equation*}
	\mathcal{J}(u, v) 
    := \|\nabla u\|_{L^{2}}^{2} 
	 + \|\nabla v\|_{L^{2}}^{2} 
	+ \|v\|_{L^{2}}^{2} - 2 \int_{\R^{d}} |u|^{2} v \, dx, 
    \end{equation*}
    \begin{equation*}
    \begin{split}
    \mathcal{S}_{\omega}(u, v) 
	:= \mathcal{J}(u, v) + \omega \|u\|_{L^{2}}^{2}  
    =  \|\nabla u\|_{L^{2}}^{2} + \omega \|u\|_{L^{2}}^{2}
	+ \|\nabla v\|_{L^{2}}^{2} 
	+ \|v\|_{L^{2}}^{2} - 2 \int_{\R^{d}} |u|^{2} v \, dx.  
    \end{split}
    \end{equation*}
A pair $(\varphi_{\omega}, 
\psi_{\omega}) \in H^{1}(\R^{d}) 
\times H^{1}(\R^{d})$ solves \eqref{sp} 
if and only if $\mathcal{S}_{\omega}^{\prime}
(\varphi_{\omega}, \psi_{\omega}) 
= 0$ in $H^{-1}(\R^{d}) 
\times H^{-1}(\R^{d})$. 
A \textit{ground state} 
is a non-trivial solution  $(\varphi_{\omega}, \psi_{\omega})$ 
of \eqref{sp} that 
minimizes 
$\mathcal{S}_{\omega}$ among all non-trivial solutions 
to \eqref{sp}. 
The existence of ground states is well known 
(see e.g. \cite[Section 2]{MR2474173}). 

To apply the potential well method 
and establish the strong instability of standing waves in four spatial dimensions, we seek solutions to \eqref{sp} through 
a constrained minimization problem. Before 
proceeding, we explain 
why the four-dimensional case is of 
particular interest. 
Roughly speaking, it corresponds to 
an $L^{2}$--critical setting. 
To explain this connection, we briefly review the stability theory of standing waves. 
We first recall the definitions of orbital stability and strong instability.
\begin{definition}\label{def-sta}
\begin{enumerate}
\item[\textrm{(i)}] 
The standing wave $(e^{i \omega t} \varphi_{\omega}, \psi_{\omega})$ is 
said to be \emph{orbitally stable} if, 
for any $\varepsilon > 0$, there exists $\delta > 0$ such that whenever 
\[
  \|(u_{0}, v_{0}, v_{1}) - (\varphi_{\omega}, \psi_{\omega}, 0)\|_{X} < \delta,
\]
then the corresponding solution $(u(t), v(t), \partial_{t} v(t))$ to \eqref{kgs} with $(u, v, \partial_{t} v)|_{t = 0} = (u_{0}, v_{0}, v_{1})$ exists globally in time and satisfies 
\[
  \sup_{t \ge 0} \inf_{(\theta, y) \in \mathbb{R} \times \mathbb{R}^{d}} \|(u(t), v(t), \partial_{t} v(t)) - (e^{i \theta} \varphi_{\omega}(\cdot + y), \psi_{\omega}(\cdot + y), 0)\|_{X} < \varepsilon.
\]
Otherwise, the standing wave $(e^{i \omega t} \varphi_{\omega}, \psi_{\omega})$ is said to be \emph{orbitally unstable}.
\item[\textrm{(ii)}]
The standing wave $(e^{i \omega t} \varphi_{\omega}, \psi_{\omega})$ is 
said to be \emph{strongly unstable} 
if for any $\varepsilon > 0$, there exists initial data 
$(u_{0}, v_{0}, v_{1}) \in X$ such that 
\[
  \|(u_{0}, v_{0}, v_{1}) - (\varphi_{\omega}, \psi_{\omega}, 0)\|_{X} < \varepsilon
\]
and the corresponding solution $(u(t), v(t), \partial_{t} v(t))$ to \eqref{kgs} with $(u, v, \partial_{t} v)|_{t = 0} = (u_{0}, v_{0}, v_{1})$ either blows up in finite time or exists globally in time and satisfies
\[
  \limsup_{t \to \infty} \|(u(t), v(t), \partial_{t} v(t))\|_{X} = \infty.
\]
\end{enumerate}
\end{definition}
Clearly, strong instability implies orbital instability.

The orbital stability and instability of standing waves for nonlinear Schr\"{o}dinger equations have been extensively studied 
(see \cite{MR646873, MR677997, 
MR691044} and references therein). 
We briefly recall some relevant results.
We first consider the  
nonlinear Schr\"{o}dinger equation with a power-type nonlinearity: 
\begin{equation}\label{NLS}
    i \partial_{t} u 
    + \Delta u  
    + |u|^{p-1}u 
    = 0 \quad 
    \text{in }\R \times \R^{d}, 
\end{equation}
where 
$d \ge 1$, $1 < p < 2^{*} - 1$, 
and 
\[
2^*=\begin{cases}
\dfrac{2d}{d - 2}, &d\ge3,\\
\infty, &d=1,2.
\end{cases}
\]
The equation \eqref{NLS} is invariant under the 
scaling 
\begin{equation} \label{eq-scale}
    u(t, x) 
    \mapsto u_{\lambda}(t, x) 
    = \lambda^{\frac{2}{p-1}} u(\lambda^{2} t, \lambda x). 
\end{equation}
Note that 
\begin{equation} \label{eq-scale2}
    \|u_{\lambda}(0, \cdot)\|_{L^{2}} 
    = \lambda^{\frac{2}{p - 1} - \frac{d}{2}}
    \|u(0, \cdot)\|_{L^{2}}.  
\end{equation}
Thus, the scaling \eqref{eq-scale} preserves 
the $L^{2}$-norm $\|\cdot\|_{L^{2}}$ 
precisely 
when $p = 1 + 4/d$.  
For this reason, the exponent 
$p = 1 + 4/d$
is referred to as ``$L^{2}$-critical''. 

When considering a standing wave solution of 
\eqref{NLS} of the form 
$u(t, x) = e^{i \omega t} \varphi_{\omega}(x)$, the 
function $\varphi_{\omega}$ 
satisfies the stationary equation
\begin{equation}\label{sp-0}
    - \Delta \varphi 
    + \omega \varphi 
    - |\varphi|^{p - 1} \varphi 
    = 0 \quad \text{in } \mathbb{R}^{d}. 
\end{equation}
Let $\phi_{\omega}$ be the ground state of \eqref{sp-0}.
Cazenave and P. -L.~Lions~\cite{MR677997} 
established the orbital 
stability of the standing waves 
$e^{i \omega t} \phi_{\omega}$ 
in the $L^{2}$-subcritical 
range $1 < p < 1 + 4/d$.
In contrast, 
Berestycki and Cazenave~\cite{MR646873} 
proved that the standing waves 
$e^{i \omega t} \phi_{\omega}$
are orbitally unstable 
in the $L^{2}$-supercritical range
$1 + 4/d < p < 2^{*} - 1$. 
In the $L^{2}$-critical case 
$p = 1 + 4/d$, 
Weinstein~\cite{MR691044} proved that, for 
any \(\omega>0\) and every nontrivial solution \(\phi_\omega\in H^1(\mathbb R^d)\) of (1.7), the standing wave \(e^{i\omega t}\phi_\omega\) is unstable, even when \(\phi_\omega\) is not a ground state.

We next recall stability results for 
the following Klein-Gordon equation: 
\begin{equation}
\partial_{t}^2 v - \Delta v + v = |v|^{p-1}v 
\quad  \text{in }\R \times \R^{d},  
\end{equation}
where $d \in \mathbb{N}$ and $1 < p < 2^{*} - 1$. 
Define 
    \[
    \omega_{c} := 
    \sqrt{\dfrac{p-1}{4 - (d-1)(p-1)}} 
    \]
and, let $\psi_{\omega}$ denote 
the unique ground state 
of the following semilinear elliptic equations:
\begin{equation}\label{sp-1}
    - \Delta \psi 
    + (1 - \omega^{2}) \psi 
    - |\psi|^{p - 1} \psi 
    = 0 \quad \mbox{in $\R^{d}$}. 
\end{equation}
Shatah~\cite{MR723756} established 
the orbital stability of 
the standing waves $e^{i \omega t} \psi_{\omega}$ for $\omega_{c} 
< |\omega| < 1$ in the range 
$1 < p < 1 + 4/d$. 
Subsequently, Shatah and Strauss~\cite{MR804458}
obtained the instability of 
the standing waves $e^{i \omega t} \psi_{\omega}$ for $|\omega| < \omega_{c}$ 
when $1 < p < 1 + 4/d$ or for $|\omega| < 1$
when $1 + 4/d \leq p < 2^{*} - 1$. 
Ohta and Todorova~\cite{MR2299435} 
further extended the result of Shatah and Strauss~\cite{MR804458} by proving that  
the standing waves 
$e^{i \omega t} \psi_{\omega}$ are 
strongly unstable for $|\omega| \leq \omega_{c}$~\footnote{Note that 
the case $|\omega| = \omega_{c}$ is included.
} 
when $1 < p < 1 + 4/d$ or for $|\omega| < 1$
when $1 + 4/d \leq p < 2^{*} - 1$ for $d 
\geq 2$. 
See Liu, Ohta, and Todorova~\cite{LiOhTo07} and Wu~\cite{MR4586806} for the case of $d = 1$. 

We next consider the following wave–Schr\"{o}dinger system:
\begin{equation}\label{ws}
\begin{cases}
i \partial_{t} u + \Delta u = -2uv & \qquad \mbox{in $\mathbb{R} \times 
\mathbb{R}^{d}$}, \\
\partial_{t}^{2} v - \Delta v = |u|^{2}
& \qquad \mbox{in $\mathbb{R} \times 
\mathbb{R}^{d}$}, 
\end{cases}
\end{equation}
Ohta~\cite{MR1393148} proved that, 
when $d = 3$, 
the standing waves $(e^{i \omega t} 
\Phi_{\omega}, \Psi_{\omega})$ are 
orbitally stable for any $\omega > 0$. 
Here $(\Phi_{\omega}, \Psi_{\omega})$ 
denotes the unique ground state of the following elliptic system:
\begin{equation}\label{lsp1}
\begin{cases}
- \Delta \Phi + \omega \Phi - 2 \Psi \Phi = 0
& \qquad \mbox{in $\mathbb{R}^{d}$}, 
\\
- \Delta \Psi - |\Phi|^{2} = 0
& \qquad \mbox{in $\mathbb{R}^{d}$}.
\end{cases}
\end{equation}
In contrast, Kikuchi~\cite{MR2600439} 
proved that, when $d = 4$, all standing 
waves are unstable, as in the 
$L^{2}$--critical case of \eqref{NLS}. 
This result suggests that 
the behavior of standing waves of 
\eqref{lsp1} in four spatial dimensions 
resembles that of standing waves for 
\eqref{NLS} in the $L^{2}$-critical case $p=1+4/d$. 
This observation motivates us to 
investigate \eqref{kgs} in four dimensions. 
For the Klein-Gordon-Schr\"{o}dinger 
system \eqref{kgs}, Kikuchi and Ohta~\cite{MR2474173, MR2585081} proved 
that, when $d=3$, 
there exist constants 
$0 < \omega_{S} < \omega_{L}$
such that the standing waves 
$(e^{i \omega t} 
\phi_{\omega}, \psi_{\omega})$
 are orbitally stable 
 for $\omega > \omega_{L}$
  and orbitally unstable for $0 < \omega < \omega_{S}$. 
Here $(\varphi_{\omega}, \psi_{\omega})$ 
denotes a ground state of 
  \eqref{sp}. 
 However, the arguments developed in \cite{MR2474173, MR2585081} do not 
 directly extend to the four-dimensional 
 case, 
 and the stability properties of 
 standing waves in this setting remain 
 incompletely understood. 

In this paper, we investigate the global dynamics of solutions to \eqref{kgs} 
in four spatial dimensions. Our approach is based on a potential well defined through a constrained minimization problem. This variational framework also enables us to establish the strong instability of the corresponding standing waves. 
The global dynamics of solutions to various dispersive equations have been extensively investigated in the literature (see, e.g., \cite{MR3102564, MR4502077, MR3940908, MR2470397, MR2491692, MR2662148, MR2838120, MR2257393,  MR3634031, MR2756065, MR2898769} and the references therein).

\subsection{Constrained minimization problem}
To formulate the minimization problem, we introduce some notation.
For $\lambda > 0$ and a 
function $u$ on $\R^{d}$, 	
we define 
	\begin{equation} \label{l2-sca}
	T_{\lambda} u(x) := \lambda^{\frac{d}{2}} 
    u(\lambda x). 
	\end{equation}

It is straightforward to verify that 
$\|T_{\lambda} u\|_{L^{2}} = \|u\|_{L^{2}}$ 
for every $u \in L^{2}(\R^{d})$ and 
$\lambda > 0$. 
We also introduce the functional $\K$ by 	
	\begin{equation}\label{K}
	\mathcal{K}(u, v) 
	:= \frac{1}{2}\frac{d}{d \lambda} 
	\mathcal{J}(T_{\lambda} u, 
	T_{\lambda} v)\bigl|_{\lambda = 1}  
	 = \|\nabla u\|_{L^{2}}^{2} 
	+ \|\nabla v\|_{L^{2}}^{2} 
	- \frac{d}{2} \int_{\R^{d}} |u|^{2} v \, dx .
	\end{equation}
In particular, when $d = 4$, 
we have
    \[
    \mathcal{K}(u, v)
    = \mathcal{J}(u, v) - \|v\|_{L^{2}}^{2}.
    \]
Then, for each $m > 0$, 
we put 
    \begin{equation} \label{Sm-set}
    S(m): = 
    \left\{
    (u, v) \in 
    H^{1}(\R^{d}) 
	\times H^{1}(\R^{d}): 
    \|u\|_{L^{2}}^{2} = m, \; 
    \K(u, v) = 0
    \right\}, 
    \end{equation}
and 
	\begin{equation} \label{m-value}
	m_{*} := \inf\left\{
	m >0 \colon  S(m) \neq \emptyset
	\right\}. 
	\end{equation}
When $d = 4$, we show that 
$m_{*} > 0$ and $S(m) 
\neq \emptyset$ for every 
$m \in (m_*, \infty)$ 
(see Lemmas \ref{lem2-1}--\ref{lem2-3} 
below).  
We then consider the constrained 
minimization problem: 
	\begin{equation} \label{mini}
	\begin{split}
	\mathcal{E}^{\mathcal{K}}(m) 
	:= 
    \inf\Big\{
	\mathcal{J}(u, v) \colon 
	& (u, v) \in S(m) 
	\Big\}. 
	\end{split}
	\end{equation}     
By the definition of $m_{*}$, 
we have 
$\E^{\K}(m) = \infty$ 
for $0 < m < m_*$. 
For $m > m_*$, we obtain the following 
existence result: 
\begin{theorem} \label{Thm:ground-sol}
Assume that $d = 4$. 
Let $m > m_{*}$. 
Then, $\E^{\K}(m)$ admits a positive minimizer 
$(\varphi_{\omega(m)}, \psi_{\omega(m)}) \in 
H_{\text{rad}}^{1}(\R^{d}) \times H_{\text{rad}}^{1}(\R^{d})$, 
which satisfies \eqref{sp} 
with $\omega = \omega(m)$ for some $\omega(m) > 0$, 
where 
    \[
    H_{\text{rad}}^{1}(\R^{d}) := 
    \left\{ 
    u \in H^{1}(\R^{d}) \colon u(x) = u(|x|)
    \right\}. 
    \]
Moreover, $\varphi_{\omega(m)}$ 
and $\psi_{\omega(m)}$ are 
radially decreasing functions of $r = |x|$. 
\end{theorem}
\begin{remark}
\begin{enumerate}
    \item[\textrm{(i)}] 
If $(\varphi_{\omega(m)}, \psi_{\omega(m)})$ is a minimizer of $\mathcal{E}^{\mathcal{K}}(m)$, then so is $(T_{\lambda} \varphi_{\omega(m)}, T_{\lambda} \psi_{\omega(m)})$ for each $\lambda>0$ (see Lemma~\ref{mini-l2} below). 
Thus, the minimizer of $\mathcal{E}^{\mathcal{K}}(m)$ is not unique. 
In fact, we exploit this property in the proof of Theorem~\ref{Thm:ground-sol}.
 \item[\textrm{(ii)}]
 It remains unclear whether 
 the frequency $\omega(m)$ 
 given in Theorem \ref{Thm:ground-sol} 
 is uniquely determined by $m$. 
\end{enumerate}
\end{remark}

Having established the existence of 
minimizers, we next investigate 
the asymptotic behavior and qualitative properties, which are of independent interest. In fact, they are essential for studying the global dynamics of solutions to \eqref{kgs} (see, e.g., \cite{MR2491692, MR2756065, MR2898769}). 
However, a complete analysis for 
\textit{all} $m > m_*$ appears difficult. 
We therefore focus on the 
two limiting regimes $m \to m_* + 0$ and 
$m \to \infty$.  
We first establish the following result: 
\begin{theorem}[Asymptotic 
behavior of minimizers as $m 
\downarrow m_*$] \label{thm-conv}
Assume that $d = 4$. 
Let 
$(\varphi_{\omega(m)}, \psi_{\omega(m)}) \in 
H_{\text{rad}}^{1}(\R^{d}) \times H_{\text{rad}}^{1}(\R^{d})$ 
be the positive minimizer of $\E^{\K}(m)$ and 
$\omega(m) > 0$ be the number given in 
Theorem \ref{Thm:ground-sol}. 
Then the following statements hold: 
\begin{enumerate}
\item[\textrm{(i)}]
$\lim_{m \to m_* + 0} \E^{\K}(m) = \infty$. 
\item[\textrm{(ii)}] 
$\lim_{m \to m_* + 0} \omega(m) = 
\infty$. 
\item[\textrm{(iii)}]
We put 
    \begin{equation} \label{rescale-1}
    \begin{split}
    (\widetilde{\varphi}_{\omega(m)}(\cdot), 
    \widetilde{\psi}_{\omega(m)}(\cdot))
&
= (T_{\omega^{-\frac{1}{2}}(m)} \varphi_{\omega(m)}, 
T_{\omega^{-\frac{1}{2}}(m)} \psi_{\omega(m)}
)\\    
& = (\omega^{-1}(m) \varphi_{\omega(m)}
(\omega^{-\frac{1}{2}}(m) \cdot), 
\omega^{-1}(m) \psi_{\omega(m)} 
(\omega^{-\frac{1}{2}}(m) \cdot)).
    \end{split}
    \end{equation}
Then, $(\widetilde{\varphi}_{\omega(m)}, 
\widetilde{\psi}_{\omega(m)})$ satisfies 
	\begin{equation}\label{msp1}
	\begin{cases}
	- \Delta \widetilde{\varphi} + \widetilde{\varphi} 
	- 2 \widetilde{\varphi} \widetilde{\psi} = 0 
	& \qquad \mbox{in $\R^{d}$}, \\
	- \Delta \widetilde{\psi} + \frac{1}{\omega(m)} \widetilde{\psi} 
	- |\widetilde{\varphi}|^{2} = 0 
	& \qquad \mbox{in $\R^{d}$}. 
	\end{cases}
	\end{equation}
As $m \downarrow m_*$, 
the rescaled pair 
$(\widetilde{\varphi}_{\omega(m)}, \widetilde{\psi}_{\omega(m)})$ converges 
strongly to $(\Phi_0, \Psi_0)$ in 
 $H^{1}(\R^{d}) 
 \times \dot{H}^{1}(\R^{d})$, 
 where $(\Phi_0, \Psi_0)$ is the unique ground state of \eqref{lsp1} with $\omega = 1$. 	
\item[\textrm{(iv)}]
There exists $\varepsilon_0>0$ such that, for any
$m\in(m_*,m_*+\varepsilon_0)$ and any fixed
$\omega>0$, there exists at most one minimizer
of $\E^{\K}(m)$, up to symmetries, that satisfies
\eqref{sp} with this frequency $\omega$.
\item[\textrm{(v)}]
There exists a sufficiently small $\varepsilon_1 > 0$ 
such that for any $m_* < m < m_* + \varepsilon_1$, 
the minimizer $(\varphi_{\omega(m)}, \psi_{\omega(m)})$ 
of $\mathcal{E}^{\mathcal{K}}(m)$ 
satisfying \eqref{sp} with $\omega = \omega(m)$ 
is
non-degenerate in 
$H_{\text{rad}}^{1}(\mathbb{R}^{d}) \times H_{\text{rad}}^{1}(\mathbb{R}^{d})$. 
More precisely, 
the linearized system at $(\varphi_{\omega(m)}, \psi_{\omega(m)})$, 
given by
\begin{equation*} 
\begin{cases}
- \Delta Z + \omega(m) Z - 2 \psi_{\omega(m)} Z - 2 \varphi_{\omega(m)} W = 0 
& \text{in } \mathbb{R}^{d}, \\
- \Delta W + W - 2 \varphi_{\omega(m)} Z = 0 
& \text{in } \mathbb{R}^{d}, 
\end{cases}   
\end{equation*}
admits only the trivial solution $(Z, W) 
= (0, 0)$ in 
$H_{\text{rad}}^{1}(\mathbb{R}^{d}) \times H_{\text{rad}}^{1}(\mathbb{R}^{d})$.
\end{enumerate}
    \end{theorem}

\begin{remark}
\begin{enumerate}
    \item[\textrm{(i)}] 
    The existence of a ground state 
    $(\Psi_0, \Phi_0)$ 
    follows from a 
    standard argument (see e.g. \cite{CaLi82}). 
    Moreover, 
    the uniqueness and non-degeneracy of 
    the ground state
    were established  
    by Krieger, Lenzmann and Raphael
    ~\cite[Theorems 3 and 4]{MR2557200}. 
    \item[\textrm{(ii)}] 
    It follows from Theorem \ref{thm-conv} 
    \textrm{(iii)} that 
    $m_* = \|\Phi_0\|_{L^{2}}^{2}$, 
    where $m_*$ is defined by \eqref{m-value}. 
\end{enumerate}  
\end{remark}
Next, we consider the case where 
$m>0$ is sufficiently large. 
\begin{theorem}[Asymptotic behavior 
of minimizers as $m 
\to \infty$] \label{thm-conv2}
Assume that $d = 4$. 
Let $(\varphi_{\omega(m)}, \psi_{\omega(m)}) \in 
H_{\text{rad}}^{1}(\R^{d}) \times H_{\text{rad}}^{1}(\R^{d})$ 
be the positive minimizer of $\E^{\K}(m)$ and 
$\omega(m) > 0$ be the number given in 
Theorem \ref{Thm:ground-sol}. 
Then, the following 
statements hold: 
\begin{enumerate}
\item[\textrm{(i)}] 
$\lim_{m \to \infty} \E^{\K}(m) = 
\frac{\sigma^{2}}{4}$, 
where 
	\begin{equation} \label{eq-sigma}
	\sigma := \inf\left\{ 
	\|\nabla u\|_{L^{2}}^{2} \colon 
	\mbox{$u \in \dot{H}^{1}(\R^{4})$ 
	with $\|u\|_{L^{4}} = 1$}
	\right\}.  
	\end{equation}
\item[\textrm{(ii)}] 
$\lim_{m \to \infty} \omega(m) = 0$. 
\item[\textrm{(iii)}]  
$\lim_{m \to \infty} M_{m} = 0$, 
where $M_{m} = \|\varphi_{\omega(m)}\|_{L^{\infty}}$.
\item[\textrm{(iv)}]
We set
\begin{equation} \label{eq13:large}
    (\widehat{\varphi}_{\omega(m)}(\cdot), \widehat{\psi}_{\omega(m)}(\cdot)) 
    = \left( M_{m}^{-1} \varphi_{\omega(m)}(M_{m}^{-1} \cdot), M_{m}^{-2}\psi_{\omega(m)}(M_{m}^{-1} \cdot) \right). 
\end{equation}
Then, $(\widehat{\varphi}_{\omega(m)}, \widehat{\psi}_{\omega(m)})$ satisfies
\begin{equation} \label{eq1:large}
\begin{cases}
- \Delta \widehat{\varphi} + \omega(m) M_{m}^{-2} \widehat{\varphi} - 2 \widehat{\varphi} \widehat{\psi} = 0
& \text{in } \mathbb{R}^{d}, \\
- M_{m}^2 \Delta \widehat{\psi} + \widehat{\psi} - |\widehat{\varphi}|^{2} = 0
& \text{in } \mathbb{R}^{d}.
\end{cases}
\end{equation}
As $m \to \infty$, 
we have $ 
\omega(m) M_{m}^{-2} \to 0$, 
and $\widehat{\varphi}_{\omega(m)}$ 
converges strongly to $V$ 
in $\dot{H}^{1}(\mathbb{R}^{d})$, 
where 
\begin{equation} \label{eq:talenti}
V(x) := \left(1 + \frac{|x|^{2}}{4}\right)^{-1}
\end{equation}
is the unique positive radial solution, up to symmetries, to the semilinear elliptic equation
\begin{equation} \label{lsp2}
- \Delta V - 2 V^{3} = 0 \quad \text{in } \mathbb{R}^{d}.
\end{equation}
\end{enumerate}
    \end{theorem}
\begin{remark}
\begin{enumerate}
    \item[\textrm{(i)}]
    The function 
$V$ is a suitable rescaling of the 
\textit{Aubin--Talenti function} 
and satisfies $\|V\|_{L^{4}}^{4} = \frac{\sigma^2}{4}$.
\item[\textrm{(ii)}]
Unlike the second equation of \eqref{msp1}, 
the second equation of \eqref{eq1:large}  
involves a singular perturbation. 
Consequently, obtaining a uniform bound 
for \(\|\nabla\widehat\psi_{\omega(m)}\|_{L^2}\) is more difficult, and the argument used to prove Theorem \ref{thm-conv} \textup{(iii)} cannot be directly applied to Theorem 
\ref{thm-conv2} \textup{(iv)}. To overcome this difficulty, we use the variational characterization of \(\sigma\).
\end{enumerate}
\end{remark}

\subsection{Unboundedness of solutions to 
\eqref{kgs}}
We next investigate the dynamics 
of solutions to 
\eqref{kgs} with initial date in a 
potential well. 
Assume that $d = 4$. 
For each $m > m_{*}$, we define 
the potential wells 
$\mathcal{R}_{+}(m)$ and 
$\mathcal{R}_{-}(m)$
by 
	\begin{equation}\label{R-+}
	\begin{split}
    & \mathcal{R}_{+}(m) := 
	\Big\{ 
	(u, v, w) \in H_{\text{rad}}^{1}(\R^{d}) \times H_{\text{rad}}^{1}(\R^{d}) \times 
    L_{\text{rad}}^{2}(\R^{d})
	\colon  
    \mathcal{E}(u, v, w) < \E^{\K}(m), \; 
	\|u\|_{L^{2}}^2 = m, \; \K(u, v) > 0
	\Big\}, 
    \end{split}
	\end{equation}
	\begin{equation}\label{R--}
	\begin{split}
    \mathcal{R}_{-}(m) := 
	\Big\{ 
	(u, v, w) \in H_{\text{rad}}^{1}(\R^{d}) \times H_{\text{rad}}^{1}(\R^{d}) \times 
    L_{\text{rad}}^{2}(\R^{d})
	\colon 
    \mathcal{E}(u, v, w) < \E^{\K}(m), \; 
	\|u\|_{L^{2}}^2 = m, \; \K(u, v) < 0
	\Big\}. 
    \end{split}
	\end{equation}
Both $\mathcal{R}_{+}(m)$ and 
$\mathcal{R}_{-}(m)$ 
are nonempty for each $m > m_*$  
(see Lemma \ref{lem1-invariant} below).   
We then establish the following result for 
solutions with initial data in 
$\mathcal{R}_{-}(m)$. 
\begin{theorem}[Unboundedness of solutions to \eqref{kgs}]\label{thm-grow-blow}
Assume that $d = 4$. 
Let $m > m_{*}$ and $(u_0, v_0, v_1) \in \mathcal{R}_{-}(m)$. 
Then the corresponding 
solution $(u, v, \p_{t}v)$ 
either blows up in finite time or exists globally in time and 
satisfies 
    \[
    \limsup_{t \to \infty} 
    \|(u(t), v(t), \partial_{t} v(t))\|_{X} = \infty.
    \]
\end{theorem}
Let $(\varphi_{\omega(m)}, \psi_{\omega(m)}) \in 
H_{\text{rad}}^{1}(\R^{d}) \times 
H_{\text{rad}}^{1}(\R^{d})$ be a  
minimizer of $\E^{\K}(m)$.  
Then we can find that 
    \[
    (\varphi_{\omega(m)}, \psi_{\omega(m)}, 0) 
\in \overline{\mathcal{R}_{-}(m)}^{
\|\cdot\|_{X}} \setminus \mathcal{R}_{-}(m).
    \]
See Lemma \ref{lem1-invariant} below. 
Combining this observation with  
Theorem \ref{thm-grow-blow}, 
we obtain the following strong instability 
result. 
\begin{theorem}[Strong instability of 
the standing waves]
\label{thm-instability}
Assume that $d = 4$. 
Let 
$m > m_*$ and 
$(\varphi_{\omega(m)}, \psi_{\omega(m)})$ be the minimizer of 
$\E^{\K}(m)$ satisfying \eqref{sp} 
with $\omega = \omega(m)$. 
Then, the standing wave 
$(e^{i \omega(m)t} \varphi_{\omega(m)}, \psi_{\omega(m)})$ 
is strongly unstable. 
\end{theorem}
The remainder of this paper is organized as follows. In Section \ref{sec-min}, 
we establish Theorem \ref{Thm:ground-sol}. 
Sections \ref{sec:m_*} 
and \ref{sec:large} are 
devoted to the proofs of Theorems \ref{thm-conv} and \ref{thm-conv2}, respectively. 
In Section \ref{sec:potential}, 
we prove Theorems \ref{thm-grow-blow} and 
\ref{thm-instability}. 
Appendix \ref{sec:cauchy} contains the proof of Theorem \ref{thm1:local-wp}, 
while Appendix \ref{sec:technical} 
provides several technical results used in the main text.

\section{Existence of a minimizer for $\E^{\K}(m)$}
\label{sec-min}
In this section, we shall show the existence of a minimizer 
for $\E^{\K}(m)$. 
To this end, we need several preparations. 
We set 
    \begin{equation}\label{eq-GN}
      C_{GN} := 
      \inf\left\{ 
      \frac{\|\nabla u\|_{L^{2}} \|u\|_{L^{2}}}{\|u\|_{L^{\frac{8}{3}}}^{2}} 
      \colon u \in H^{1}(\R^{4}) \setminus \{0\}
      \right\}. 
    \end{equation}
We first show the following: 
	\begin{lemma} \label{lem2-1}
	We have  
	$\mathcal{K} (u, v) > 0$
	for $0 < \|u\|_{L^{2}}^{2} < 
    C_{GN}^{2}\sigma^{2}$. 
	\end{lemma}
\begin{proof}
By the H\"{o}lder, Gagliardo-Nirenberg,
Sobolev and Young inequalities, 
we obtain 
	\begin{equation*}
	\begin{split}
	2 \int_{\R^{4}} |u|^{2} v \, dx 
	\leq 2 \|u\|_{L^{\frac{8}{3}}}^{2} \|v\|_{L^{4}} 
	\leq \frac{2}{C_{GN} \sigma} \|\nabla u\|_{L^{2}} \|u\|_{L^{2}} 
	\|\nabla v\|_{L^{2}} 
    \leq \frac{1}{C_{GN} \sigma}  \|u\|_{L^{2}}
    \left(\|\nabla u\|_{L^{2}}^{2} 
	+ \|\nabla v\|_{L^{2}}^{2} \right). 
	\end{split}
	\end{equation*}
This yields that 
    \[
    \mathcal{K}(u, v) = \|\nabla u\|_{L^{2}}^{2} 
	+ \|\nabla v\|_{L^{2}}^{2} -2 \int_{\R^{4}} |u|^{2} v \, dx > 0
    \]
for $0 < \|u\|_{L^{2}}^2 < C_{GN}^{2} \sigma^{2}$. 
This completes the proof. 
\end{proof}
It follows from Lemma \ref{lem2-1} that 
$m_{*} > 0$, where $m_*$ is defined by \eqref{m-value}. 
Let $Q$ be the unique positive solution to 
 	\begin{equation} \label{scf}
	- \Delta Q + Q - 2 |Q| Q = 0 
	\qquad \mbox{in $\R^{4}$}. 
	\end{equation}
See Berestycki and P. -L.~Lions~\cite{BeLi83} for the 
existence and Kwong~\cite{MR969899} for the uniqueness. 
We see that $m_{*}$ is finite in the following lemma:
\begin{lemma}\label{lem2-2}
Let $m_*$ be the one defined by \eqref{m-value}. 
We have $m_{*} \leq 2\|Q\|_{L^{2}}^{2}$.
\end{lemma}
\begin{proof}
We first claim that 
		\begin{equation} \label{poho}
		\|\nabla Q\|_{L^{2}}^{2} 
	= \frac{4}{3} \|Q\|_{L^{3}}^{3}, 
	\qquad 
	\|Q\|_{L^{2}}^{2} 
	= \frac{2}{3}\|Q\|_{L^{3}}^{3}. 
		\end{equation}
Observe that $Q$ satisfies 
	\begin{equation*} 
	\|\nabla Q\|_{L^{2}}^{2} + \|Q\|_{L^{2}}^{2} 
	- 2 \|Q\|_{L^{3}}^{3} = 0,  \qquad
	\|\nabla Q\|_{L^{2}}^{2} 
	+ 2 \|Q\|_{L^{2}}^{2} 
	- \frac{8}{3} \|Q\|_{L^{3}}^{3} = 0. 
	\end{equation*}
See e.g. \cite[Lemma 8.1.2]{Ca03}. 
These identities yield that \eqref{poho}. 

It follows from \eqref{poho} that 
$\|\nabla Q\|_{L^{2}}^{2} = \frac{4}{3} \|Q\|_{L^{3}}^{3}$. 
Thus, we have
	\[
	\K(\sqrt{2} Q, Q) 
	= 2 \|\nabla Q\|_{L^{2}}^{2} + \|\nabla Q\|_{L^{2}}^{2} 
    -4 \|Q\|_{L^{3}}^{3}
    = 3 \|\nabla Q\|_{L^{2}}^{2} - 4 \|Q\|_{L^{3}}^{3} 
    = 0. 
	\]
Thus, from the definition of $m_{*}$ (see \eqref{m-value}), 
we obtain $m_{*} \leq 2 \|Q\|_{L^{2}}^{2}$. 
\end{proof}
\begin{lemma}\label{lem2-3}
Let $S(m)$ 
be the subset of $H^{1}(\R^{4}) \times H^{1}(\R^{4})$
defined by \eqref{Sm-set}. 
Then, we have $S(m) \neq \emptyset$ 
for $m > m_{*}$. 
\end{lemma}
\begin{proof}
We claim that if $m_{1} < m_{2}$ and $S(m_{1}) \neq 
\emptyset$, we have $S(m_{2}) \neq \emptyset$.  
Since $S(m_{1}) \neq \emptyset$, there exists $(u_{1}, v_{1}) 
\in H^{1}(\R^{4}) \times H^{1}(\R^{4})$ such that 
$\|u_{1}\|_{L^{2}}^{2} = m_{1}$ and $\K(u_{1}, v_{1}) = 0$. 
Setting $u_{2} = (m_{2}/m_{1})^{1/2} u_{1}$, 
one has $\|u_{2}\|_{L^{2}}^{2} = m_{2}$. 
For $\lambda > 0$, we define 
    \[
    f(\lambda) = \K(u_{2}, \lambda v_{1}) 
    = \|\nabla u_{2}\|_{L^{2}}^{2} + 
    \lambda^{2} \|\nabla v_{1}\|_{L^{2}}^{2} 
    - 2 \lambda \int_{\R^{4}} |u_{2}^{2}| v_{1} \, dx.  
    \]
Then, since $u_{2} = (m_{2}/m_{1})^{1/2} u_{1}$ and 
 $\K(u_{1}, v_{1}) = 0$, we obtain 
    \begin{equation} 
    \label{eq-15}
    \begin{split}
    f(1) 
    & = \frac{m_{2}}{m_{1}} 
    \|\nabla u_{1}\|_{L^{2}}^{2} + 
     \|\nabla v_{1}\|_{L^{2}}^{2} 
    - 2 \frac{m_{2}}{m_{1}} 
    \int_{\R^{4}} |u_{1}^{2}| v_{1} \, dx \\
   & < \frac{m_{2}}{m_{1}} 
   \left(\|\nabla u_{1}\|_{L^{2}}^{2} + 
     \|\nabla v_{1}\|_{L^{2}}^{2} 
    - 2 \int_{\R^{4}} |u_{1}^{2}| v_{1} 
    \, dx \right) 
    = \frac{m_{2}}{m_{1}} \K(u_{1}, v_{1}) = 0. 
    \end{split}
    \end{equation}
In addition, we have 
$f(0) = \|\nabla u_{2}\|_{L^{2}}^{2}
> 0$. 
Thus, there exists $\lambda_{*} \in (0, 1)$ such that 
$\K(u_{2}, \lambda_{*} v_{1}) = 0$. 
Then, setting $v_{2} = \lambda_{*} v_{1}$, we see that 
$(u_{2}, v_{2}) \in S(m_{2})$, that is, $S(m_{2}) \neq 
\emptyset$. 

From the definition of $m_{*}$ and $m_{*} < \infty$, hence 
for any $\e>0$, there exists $m_{1} \in (m_*, m_{*} + \e)$ such 
that $S(m_{1}) \neq \emptyset$. 
This together with the above argument implies that 
$S(m) \neq \emptyset$ for $m > m_*$. 
This completes the proof. 

\end{proof}

 \begin{lemma} \label{vari-chara}
 Let $(u, v) \in H^{1}(\R^{4}) \times 
    H^{1}(\R^{4})$ satisfy 
$\mathcal{K}(u, v) < 0$ and $\|u\|_{L^{2}}^{2} = m$. 
Then, we obtain 
$\mathcal{E}^{\K}(m) < \|v\|_{L^{2}}^{2}$
 \end{lemma}

\begin{proof}
We argue as in the proof of Lemma \ref{lem2-3}. 
Let $(u, v) \in H^{1}(\R^{4}) \times H^{1}(\R^{4})$ satisfy 
$\mathcal{K}(u, v) < 0$ and $\|u\|_{L^{2}}^{2} = m$. 
Then, for $\lambda \geq 0$, 
we put 
\[
f(\lambda) := \mathcal{K}(u, \lambda v) 
= \|\nabla u\|_{L^{2}}^{2} 
+ \lambda^{2} \|\nabla v\|_{L^{2}}^{2} 
- 2 \lambda \int |u|^{2} v \, dx. 
\]
Then, we have 
$f(0) = \|\nabla u\|_{L^{2}}^{2} > 0$~\footnote{
Since $m > 0$, we have $u \neq 0$. 
} and $f(1) = 
\mathcal{K}(u, v) < 0$. 
Thus, there exists $\lambda_{m} \in 
(0, 1)$ such that $f(\lambda_{m}) = 0$. 
Namely, one has $\mathcal{K}(u, \lambda_{m} v) = 0$, 
$\|u\|_{L^{2}}^{2} = m$.
It follows from the definition $\E^{\K}(m)$ 
(see \eqref{mini}) that 
	\begin{equation} \label{po-lem-1}
	\begin{split}
	\mathcal{E}^{\mathcal{K}}(m)
 \leq  \mathcal{J}(u, \lambda_{m} v) 
	= \mathcal{J}(u, \lambda_{m} v) 
	- \mathcal{K} (u, \lambda_{m}v)  
    = \|\lambda_{m} v\|_{L^{2}}^{2} 
	= \lambda_{m}^{2} \|v\|_{L^{2}}^{2} 
    < \|v\|_{L^{2}}^{2}. 
 	\end{split}
	\end{equation}
This completes the proof. 
\end{proof}
 \begin{lemma} \label{vari-chara2}
The variational value 
$\E^{\K}(m)$ defined by \eqref{mini}
can be written by the following: 
	\begin{equation} \label{mini-1}
	\begin{split}
	\mathcal{E}^{\mathcal{K}}(m) 
	= \inf\big\{
	\|v\|_{L^{2}}^{2} 
	\colon 
    (u, v) \in H^{1}(\R^{4}) \times 
    H^{1}(\R^{4})
    \setminus \{(0, 0)\}, 
    \mathcal{K} (u, v) \leq 0, 
	\|u\|_{L^{2}}^{2} = m  
	\big\}. 
	\end{split}
	\end{equation}
\end{lemma} 
\begin{proof}
Let $(u, v) \in H^{1}(\R^{4}) \times H^{1}(\R^{4})$ 
satisfy  $\mathcal{K}(u, v) = 0$. 
Then, we obtain 
    \[
 \mathcal{J}(u, v) 
	= \mathcal{J}(u, v) 
	- \mathcal{K} (u, v) 
	= \|v\|_{L^{2}}^{2}.  
    \]
This implies that 
\begin{equation} \label{mini-2}
	\begin{split}
	\mathcal{E}^{\mathcal{K}}(m) 
	& = \inf\left\{
	\|v\|_{L^{2}}^{2} 
	\colon (u, v) \in H^{1}(\R^{4}) \times H^{1}(\R^{4}) \setminus \{(0, 0)\}, 
	\mathcal{K} (u, v) = 0, \; 
	\|u\|_{L^{2}}^{2} = m  
	\right\}. 
	\end{split}
	\end{equation}    
We see from this and Lemma \ref{vari-chara} that 
\eqref{mini-1} holds. 
\end{proof}


\begin{lemma}\label{mini-l2}
Let $m > m_{*}$. 
Assume that $\{(u_{n}, v_{n})\}$ is a 
minimizing sequence of $\mathcal{E}^{\mathcal{K}}(m)$. 
Then, for each $\lambda > 0$, 
$\{(T_{\lambda} u_{n}, T_{\lambda} v_{n})\}$ is also a 
minimizing sequence of $\mathcal{E}^{\mathcal{K}}(m)$.
\end{lemma}
\begin{proof}
Observe from \eqref{l2-sca} that for any $\lambda > 0$, we have 
	\[
	\|T_{\lambda} u_{n}\|_{L^{2}}^{2} = \|u_{n}\|_{L^{2}}^{2} = m, 
	\quad
	\|T_{\lambda} v_{n}\|_{L^{2}}^{2} = \|v_{n}\|_{L^{2}}^{2}, 
	\quad 
	\mathcal{K} (T_{\lambda} u_{n}, T_{\lambda} v_{n}) 
	= \lambda^{2} \mathcal{K}(u_{n}, v_{n}) = 0. 
	\]
Thus, we see from \eqref{mini-1} that 
$\{(T_{\lambda} u_{n}, T_{\lambda} v_{n})\}$ is also a 
minimizing sequence of $\mathcal{E}^{\mathcal{K}}(m)$. 
This completes the proof. 
\end{proof}

\begin{lemma}\label{mini-po}
For each $m > m_{*}$, we have 
$\mathcal{E}^{\mathcal{K}}(m) \geq  \frac{\sigma^{2}}{4}$,  
where $\sigma > 0$ is the constant given by 
\eqref{eq-sigma}. 
\end{lemma}
\begin{proof}
Let $\{(\varphi_{n}, \psi_{n})\}$ be a minimizing sequence 
of $\mathcal{E}^{\mathcal{K}}(m)$. 
We put 
	\begin{equation} \label{eq-14}
	\lambda_{n} = \left(2 
     \int_{\R^{4}} |\varphi_{n}|^2 \psi_{n} \, dx
	\right)^{- \frac{1}{2}}
	\end{equation}
and $(\widetilde{\varphi}_{n}, \widetilde{\psi}_{n}) 
= (T_{\lambda_{n}} \varphi_{n}, T_{\lambda_{n}} \psi_{n})$. 
Then, we see from Lemma \ref{mini-l2} that 
$(\widetilde{\varphi}_{n}, \widetilde{\psi}_{n})$ is also a minimizing 
sequence of $\mathcal{E}^{\mathcal{K}}(m)$.
It follows from $\mathcal{K} (\widetilde{\varphi}_{n}, \widetilde{\psi}_{n}) 
\leq 0$ and \eqref{eq-14} that 
	\begin{equation*} 
	1 =  2 \int_{\R^{4}} 
    |\widetilde{\varphi}_{n}|^2 
    \widetilde{\psi}_{n} \, dx
	\geq \|\nabla \widetilde{\varphi}_{n}\|_{L^{2}}^{2} 
	+ \|\nabla \widetilde{\psi}_{n}\|_{L^{2}}^{2} 
	\geq \|\nabla \widetilde{\varphi}_{n}\|_{L^{2}}^{2}. 
	\end{equation*}
This implies that $\sup_{n \in \mathbb{N}} 
\|\nabla \widetilde{\varphi}_{n}\|_{L^{2}}^{2} \leq 1$. 
Then, by the  H\"{o}lder  and Sobolev inequalities, 
we obtain 
	\[
	1 = 2 
     \int_{\R^{4}} |\widetilde{\varphi}_{n}|^2 
     \widetilde{\psi}_{n} \, dx
	\leq 2 \|\widetilde{\varphi}_{n}\|_{L^{4}}^{2} 
	\|\widetilde{\psi}_{n}\|_{L^{2}} 
	\leq \frac{2}{\sigma} \|\nabla \widetilde{\varphi}_{n}\|_{L^{2}}^{2} 
	\|\widetilde{\psi}_{n}\|_{L^{2}}  
	\leq \frac{2}{\sigma} \|\widetilde{\psi}_{n}\|_{L^{2}}. 
	\]
This yields that $\|\widetilde{\psi}_{n}\|_{L^{2}} \geq 
\sigma/2$. 
Therefore, we obtain 
	\[
	\mathcal{E}^{\mathcal{K}}(m) 
	= \lim_{n \to \infty} \|\psi_{n}\|_{L^{2}}^{2} 
	= \lim_{n \to \infty} \|\widetilde{\psi}_{n}\|_{L^{2}}^{2} 
	\geq \frac{\sigma^{2}}{4}. 
	\]
 This completes the proof. 
 \end{proof}
 \begin{lemma}[Monotonicity of 
$\E^{\K}(m)$]\label{monotone}
The variational value 
$\E^{\K}(m)$ is non-increasing 
on $m \in [m_{*}, \infty)$. 
\end{lemma}
\begin{proof}
Let $m_{1} < m_{2}$. 
Suppose that $(u, v) \in H^{1}(\R^{4}) \times H^{1}(\R^{4})$ 
satisfies $\K(u, v) = 0$ and $\|u\|_{L^{2}}^{2} = m_{1}$. 
Let $\theta_{*} := (m_{2}/m_{1})^{\frac{1}{2}} > 1$. 
As in \eqref{eq-15}, 
we see from $v \neq 0$ that $\|\theta_{*} u\|_{L^{2}}^{2} = m_{2}$ 
and 
\begin{equation*}
\begin{split}
\K(\theta_{*} u, v) 
= \theta_{*}^{2} 
\|\nabla u\|_{L^{2}}^{2} 
+ \|\nabla v\|_{L^{2}}^{2} 
- 2 \theta_{*}^{2}
    \int_{\R^{4}} |u|^{2} v \, dx 
< 
\theta_{*}^{2} \K(u, v) = 0. 
\end{split}
\end{equation*}
By Lemma \ref{vari-chara}, 
we have 
$\E^{\K}(m_{2}) < \|v\|_{L^{2}}^{2}$. 
Taking an infimum on $(u, v) \in H^{1}(\R^{4}) \times 
H^{1}(\R^{4})$ satisfying
$\mathcal{K}(u, v) \leq 0$ and 
$\|u\|_{L^{2}}^{2} = m_1$, 
we see that 
$\E^{\K}(m_{2}) \leq \E^{\K}(m_{1})$. 
Thus, we have obtained the desired result. 
\end{proof}
We obtain a convergence result, which 
is needed not only 
to establish the existence 
of a minimizer for $\E^{\K}(m)$ but also 
for investigating the asymptotic of 
$\E^{\K}(m)$ as $m$ goes to $m_*$. 
\begin{lemma}\label{lem-nonc}
Let $\{u_{n}\}$  be bounded sequence 
in $H_{\text{rad}}^{1}(\R^{4})$ and 
$\{v_{n}\}$ be a bounded one in 
$\dot{H}_{\text{rad}}^{1}(\R^{4})$
satisfying 
\begin{equation*}
\lim_{n\to \infty}u_{n} =u_{\infty}
\qquad 
\mbox{weakly in $H^{1}(\R^{4})$}, 
\qquad 
\lim_{n\to \infty}v_{n} =v_{\infty}
\qquad 
\mbox{weakly in $\dot{H}^{1}(\R^{4})$}
\end{equation*}
for some functions $u_{\infty}\in H^{1}(\R^{4})$ 
and $v_{\infty} \in \dot{H}^{1}(\R^{4})$. 
Then, up to some subsequence, we have
\begin{equation*}
\lim_{n \to \infty} 
 \int_{\R^{4}} |u_{n}|^2 v_{n} \, dx 
 = 
\int_{\R^{4}} |u_\infty|^{2} v_\infty \, dx. 
\end{equation*}
\end{lemma}
\begin{proof}
For each $R > 0$, we have by H\"{o}lder's inequality that 
    \begin{equation} \label{eq-13}
    \begin{split}
    \biggl| 
    \int_{\R^{4}} |u_{n}|^2 
v_{n} \, dx  
- \int_{\R^{4}} |u_{\infty}|^2 v_{\infty} \, dx 
    \biggl| 
    & \leq 
    \biggl| 
    \int_{\R^{4}} |u_{n}|^2 
v_{n} \, dx  
- \int_{\R^{4}} |u_{\infty}|^2 v_{n} \, dx
\biggl| + 
\biggl| 
    \int_{\R^{4}} |u_{\infty}|^2 
v_{n} \, dx  
- \int_{\R^{4}} |u_{\infty}|^2 v_{\infty} \, dx
\biggl| \\
& \leq \|u_{n} 
+ u_{\infty}\|_{L^{\frac{8}{3}}} 
\|u_{n} 
- u_{\infty}\|_{L^{\frac{8}{3}}} 
\|v_{n} 
\|_{L^{4}} 
+ \|u_{\infty}\|_{L^{3}(B_{R})}^2 
\|v_{n} 
- v_{\infty}\|_{L^{3}(B_{R})} 
\\
& \quad + \|u_{\infty}\|_{L^{\frac{8}{3}}(\R^{4} \setminus B_{R})}^2 
\|v_{n} - v_{\infty}
\|_{L^{4}(\R^{4} \setminus B_{R})}. 
    \end{split}
    \end{equation}
Since $\{v_{n}\}$ is bounded 
in $\dot{H}^{1}(\R^{4})$ and $u_{\infty} \in H^{1}(\R^{4})$, 
for any $\e > 0$, 
we can take $R > 0$ sufficiently large so that
    \[
    \|u_{\infty}\|_{L^{\frac{8}{3}}(\R^{4} \setminus B_{R})}^2 
\|v_{n} - v_{\infty}
\|_{L^{4}(\R^{4} \setminus B_{R})} 
\leq \|u_{\infty}\|_{L^{\frac{8}{3}}(\R^{4} \setminus B_{R})}^2 
(\|v_{n}\|_{L^{4}(\R^{4} \setminus B_{R})}
+ \|v_{\infty}\|_{L^{4}(\R^{4} \setminus B_{R})}) 
< \frac{\e}{3}. 
    \]
By the compact embedding 
$H^{1}_{\text{rad}}(\R^{4}) 
\hookrightarrow L^{\frac{8}{3}} (\R^{4})$ 
(see e.g. Berestycki and P.-L.~Lions~\cite[Theorem A.I']{BeLi83}) 
and the 
Rellich-Kondrachov theorem, 
there exists $n_{0} \in \mathbb{N}$ such that 
we have 
  \[
\|u_{n} 
+ u_{\infty}\|_{L^{\frac{8}{3}}} 
\|u_{n} 
- u_{\infty}\|_{L^{\frac{8}{3}}} 
\|v_{n} 
\|_{L^{4}} < \frac{\e}{3}, \qquad  
 \|u_{\infty}\|_{L^{3}(B_{R})}^2 
\|v_{n} 
- v_{\infty}\|_{L^{3}(B_{R})} < \frac{\e}{3} 
  \] 
for all $n \geq n_0$. 
Therefore, by \eqref{eq-13}, we obtain 
    \[
    \biggl| 
    \int_{\R^{4}} |u_{n}|^2 
v_{n} \, dx  
- \int_{\R^{4}} |u_{\infty}|^2 v_{\infty} \, dx 
    \biggl| < \e 
    \]
for all $n \geq n_0$. 
Therefore, we deduce that 
$\lim_{n \to \infty} 
\int_{\R^{4}} |u_{n}|^2 
v_{n} \, dx 
= \int_{\R^{4}} |u_{\infty}|^2 v_{\infty} \, dx$. 
This completes the proof. 
\end{proof}
\begin{proposition}[Existence of a minimizer for $\E^{\K}(m)$]\label{ex-mini}
For each $m > m_{*}$, 
there exists a minimizer 
$(\varphi_{\omega(m)}, \psi_{\omega(m)}) 
\in H_{\text{rad}}^{1}(\R^{4}) \times 
H_{\text{rad}}^{1}(\R^{4})$ 
of $\E^{\K}(m)$. 
In addition, both functions $\varphi_{\omega(m)}$ and $\psi_{\omega(m)}$ are non-negative and 
decreasing in $r = |x|$. 
\end{proposition}
\begin{proof}
Let $\{(\varphi_{n}, \psi_{n})\} \subset H^{1}(\R^{4}) \times H^{1}(\R^{4})$ 
be a minimizing sequence of $\E^{\K} (m)$. 
Let $\varphi_{n}^{*}$ and $\psi_{n}^{*}$ 
be Schwartz rearrangements of 
$|\varphi_{n}|$ and $|\psi_{n}|$, 
respectively. 
Then, since 
$ \int_{\R^{4}} |\varphi_{n}|^2 |\psi_{n}| \, dx
\leq 
\int_{\R^{4}} |\varphi_{n}^{*}|^2 |\psi_{n}^{*}| \, dx$ 
(see Theorem 3.4 of Lieb and 
Loss~\cite{LiebLoss2001}), 
it follows from \eqref{mini-1} that 
$\{(\varphi_{n}^{*}, \psi_{n}^{*})\} \subset 
 H_{\text{rad}}^{1}(\R^{4}) \times 
H_{\text{rad}}^{1}(\R^{4})$ 
is also a minimizing sequence of $\E^{\K} (m)$.
Thus, we may assume that $\varphi_{n}$ and $\psi_{n}$ are non-negative, radially symmetric 
and monotone decreasing in $r = |x| > 0$.

We put 
	\[
	(\widetilde{\varphi}_{n}, \widetilde{\psi}_{n}) 
	:= (T_{\lambda_{n}} \varphi_{n}, T_{\lambda_{n}} \psi_{n}), 
	\qquad 
	\lambda_{n} := \left(2
    \int_{\R^{4}} |\varphi_{n}|^2 |\psi_{n}| \, dx 
    \right)^{-\frac{1}{2}}.  
	\]
Then,  
we obtain  
	\[
	\|\nabla \widetilde{\varphi}_{n}\|_{L^{2}}^{2} + 
	\|\nabla \widetilde{\psi}_{n}\|_{L^{2}}^{2} 
    = \lambda_{n}^2 
    \left(\|\nabla \varphi_{n}\|_{L^{2}}^{2} + 
	\|\nabla \psi_{n}\|_{L^{2}}^{2}\right)
    \leq 2 \lambda_{n}^{2} 
    \int_{\R^{4}} |\varphi_{n}|^2 |\psi_{n}| \, dx
    = 1. 
	\]
Hence, 
$\{(\widetilde{\varphi}_{n}, \widetilde{\psi}_{n})\} \subset 
H_{\text{rad}}^{1}(\R^{4}) \times 
H_{\text{rad}}^{1}(\R^{4})$ 
is a bounded minimizing sequence of $\E^{\K}(m)$.  
Thus, there exists $(\Phi_m, 
\Psi_{m}) \in H_{\text{rad}}^{1}(\R^{4}) \times 
H_{\text{rad}}^{1}(\R^{4})$ and 
a subsequence of $\{(\widetilde{\varphi}_{n}, \widetilde{\psi}_{n})\}$ 
(we still denote it by the same letter) such that 
	\[
	\lim_{n \to \infty} (\widetilde{\varphi}_{n}, \widetilde{\psi}_{n}) 
= (\Phi_m, \Psi_m) \qquad 
\mbox{weakly in $H^{1}(\R^{4}) \times H^{1}(\R^{4})$}.
	\] 
By the compact embedding 
$H^{1}_{\text{rad}}(\R^{4}) 
\hookrightarrow L^{q}(\R^{4})\; (2 < q < 4)$ (see e.g. Berestycki and P.-L.~Lions~\cite[Theorem A.I']{BeLi83}), 
we have
	\[
	\lim_{n \to \infty} (\widetilde{\varphi}_{n}, \widetilde{\psi}_{n}) 
= (\Phi_m, \Psi_m) \qquad 
\mbox{strongly in $L^{q}(\R^{4}) \times L^{q}(\R^{4})$ 
for $2 < q < 4$}.
	\] 
By H\"{o}lder's inequality, we obtain 
	\begin{equation} \label{ex-mini-eq1}
	1 = 2 
    \int_{\R^{4}} |\widetilde{\varphi}_{n}|^2 
    \widetilde{\psi}_{n} \, dx
	\leq 2 \|\widetilde{\varphi}_{n}\|_{L^{3}}^{2} 
	\|\widetilde{\psi}_{n}\|_{L^{3}} = 
	2 \|\Phi_m\|_{L^{3}}^{2} \|\Psi_m\|_{L^{3}}.  
	\end{equation}
This yields that $\Phi_m \neq 0$ and $\Psi_m \neq 0$. 
It follows from Lemma \ref{lem-nonc} and weak 
lower semicontinuity of $H^{1}$-norm that 
\begin{align}
& \mathcal{K} (\Phi_m, \Psi_m) \leq 
\liminf_{n \to \infty} 
\mathcal{K} (\widetilde{\varphi}_{n}, \widetilde{\psi}_{n})  
= \liminf_{n \to \infty}
\lambda_{n}^2
\mathcal{K} (\varphi_{n}, \psi_{n})
\leq 0.
\label{BL-1} 
\end{align}
This together with \eqref{mini-1} 
implies that 
$\|\Phi_{m}\|_{L^{2}}^{2} \geq m_{*}$. 

From the weak lower semicontinuity, we have 
	\[
	\|\Phi_m\|_{L^{2}}^{2} \leq \liminf_{n \to \infty} \|\widetilde{\varphi}_{n}\|_{L^{2}}^{2} 
	= \liminf_{n \to \infty} \|\varphi_{n}\|_{L^{2}}^{2} = m. 
	\]
This together with 
the monotonicity of $\E^{\K}(m)$ 
(see Lemma \ref{monotone}), \eqref{mini-1}, 
and the weak lower semicontinuity yields that 
\[
\E^{\K} (m) \leq \E^{\K} (\|\Phi_m\|_{L^{2}}^{2}) \leq 
\|\Psi_m\|_{L^{2}}^{2} 
\leq \liminf_{n \to \infty} 
\|\widetilde{\psi}_{n}\|_{L^{2}}^{2} 
= \liminf_{n \to \infty} 
\|\psi_{n}\|_{L^{2}}^{2} 
= \E^{\K}(m). 
\]
Thus, one has $\|\Psi_m\|_{L^{2}}^{2} = \E^{\K}(m)$. 

We claim that $\|\Phi_{m}\|_{L^{2}}^{2} = m$. 
Suppose to the contrary that $\|\Phi_{m}\|_{L^{2}}^{2} < m$. 
Let $\theta_{*} = (m/\|\Phi_{m}\|_{L^{2}}^{2})^{\frac{1}{2}} 
> 1$. 
By a similar argument to \eqref{eq-15}, 
we see from $\Psi_{m} \neq 0$ that $\|
\theta_{*} \Phi_m\|_{L^{2}}^{2} = m$ 
and 
\begin{equation*}
\begin{split}
\K(\theta_{*} \Phi_m, \Psi_m) 
 < 
\theta_{*}^{2} \K(\Phi_m, \Psi_m) \leq 0. 
\end{split}
\end{equation*}
By Lemma \ref{vari-chara}, 
we have 
$\E^{\K}(m) < \|\Psi_{m}\|_{L^{2}}^{2}$, 
which is a contradiction. 

By \eqref{BL-1} and Lemma \ref{vari-chara}, 
we obtain $\mathcal{K}(\Phi_m, \Psi_m) = 0$. 
Therefore, we see that 
$(\Phi_m, \Psi_m)$ is a minimizer of $\E^{\K}(m)$. 
Since $\{(\widetilde{\varphi}_{n}, \widetilde{\psi}_{n})\}$ is non-negative, 
we see that both of functions 
$\Phi_m$ and $\Psi_m$ 
are non-negative. 
This completes the proof. 
\end{proof}
We are now in the position to prove Theorem 
\ref{Thm:ground-sol}. 
\begin{proof}[Proof of Theorem \ref{Thm:ground-sol}]
Let $(\varphi_{\omega(m)}, \psi_{\omega(m)}) 
\in H_{\text{rad}}^{1}(\R^{4}) \times 
H_{\text{rad}}^{1}(\R^{4})$ be a non-negative minimizer 
of $\E^{\K}(m)$~\footnote{Existence of a minimizer is 
proved in Proposition \ref{ex-mini}.}. 
We shall show that 
$\mathcal{K}^{\prime}(\varphi_{\omega(m)}, \psi_{\omega(m)})$ and 
$(\|\varphi_{\omega(m)}\|_{L^{2}}^{2})^{\prime}$ are linearly independent 
by contradiction. 
Suppose to the contrary that there exists 
$C_{m} \in \R$ such that 
$\mathcal{K}^{\prime}(\varphi_{\omega(m)}, \psi_{\omega(m)}) = C_{m} 
(\|\varphi_{\omega(m)}\|_{L^{2}}^{2})^{\prime}$. 
Then, we have 
	\[
	\begin{cases}
	- \Delta \varphi_{\omega(m)} - 2 \varphi_{\omega(m)} \psi_{\omega(m)} = - C_{m} \varphi_{\omega(m)}
	& \qquad \mbox{in $\R^{4}$},  \\
	- \Delta \psi_{\omega(m)} - |\varphi_{\omega(m)}|^{2} = 0 
	& \qquad \mbox{in $\R^{4}$}.
	\end{cases}
	\]
We claim that the following inequality holds: 
	\begin{equation} \label{eq-b1-1}
	\psi_{\omega(m)}(x) \geq \frac{\psi_{\omega(m)}(1)}{2 |x|^{2}} 
	\qquad \mbox{for $x \in \R^{4}$ satisfying $|x| \geq 1$}. 
	\end{equation}
Indeed, putting $f_{m}(x) = \frac{\psi_{\omega(m)}(1)}{2 |x|^{2}}$, 
we have $- \Delta f_{m} = 0$ in $\R^{4} 
\setminus \{0\}$. 
 In addition, we have $\psi_{\omega(m)}(x) \geq f_{m}(x)$ for $x \in \R^{4}$ 
 with $|x| = 1$. 
Since 
    \[
    - \Delta \psi_{\omega(m)} = |\varphi_{\omega(m)}|^{2} > 0,
    \]
we see from the maximum27
principle that 
\eqref{eq-b1-1} holds. 

One has by \eqref{eq-b1-1} that 
	\[
	\begin{split}
	\E^{\K}(m) = 
	\|\psi_{\omega(m)}\|_{L^{2}}^{2} \geq \int_{|x| \leq R} |\psi_{\omega(m)}(x)|^{2} \, dx 
	\geq C \int_{1}^{R} r^{-1} dr 
	=  C \log R \to \infty \quad \mbox{as $R \to \infty$},  
	\end{split}
	\]		
which is a contradiction. 
Hence, $\mathcal{K}^{\prime}(\varphi_{\omega(m)}, \psi_{\omega(m)})$ and 
$(\|\varphi_{\omega(m)}\|_{L^{2}}^{2})^{\prime}$ are linearly independent. 

Then, there exists Lagrange multiplier $(\lambda, \mu) \in \R^{2}$ 
such that 
	\begin{equation} \label{eq-4}
	\begin{cases}
	- \Delta \varphi_{\omega(m)} - 2 \varphi_{\omega(m)} \psi_{\omega(m)} = 
	\lambda (- \Delta \varphi_{\omega(m)} - 2 \varphi_{\omega(m)} \psi_{\omega(m)})
	 - \mu \varphi_{\omega(m)}
	& \qquad \mbox{in $\R^{4}$},  \\
	- \Delta \psi_{\omega(m)} + \psi_{\omega(m)} - |\varphi_{\omega(m)}|^{2} = 
	\lambda (- \Delta \psi_{\omega(m)} - |\varphi_{\omega(m)}|^{2}) 
	& \qquad \mbox{in $\R^{4}$}.
	\end{cases}
	\end{equation}
We shall show that $\lambda \neq 1$ by contradiction. 
Suppose to the contrary that $\lambda = 1$. 
Then, from the second equation of \eqref{eq-4}, we have 
$\psi_{\omega(m)} \equiv 0$, which is absurd. 
Hence, $\lambda \neq 1$. 

Then, \eqref{eq-4} can be written by 
	\begin{equation} \label{eq-4-1}
	\begin{cases}
	- \Delta \varphi_{\omega(m)} + \frac{\mu}{1 - \lambda} \varphi_{\omega(m)} 
	- 2 \varphi_{\omega(m)} \psi_{\omega(m)} 
	 = 0
	& \qquad \mbox{in $\R^{4}$},  \\
	- \Delta \psi_{\omega(m)} + \frac{1}{1 - \lambda} \psi_{\omega(m)} - |\varphi_{\omega(m)}|^{2} 
	= 0
	& \qquad \mbox{in $\R^{4}$}.
	\end{cases}
	\end{equation}
This yields that $\frac{\mu}{1 - \lambda} > 0$. 
Indeed, suppose to the contrary that $\frac{\mu}{1 - \lambda} \leq 0$. 
Then, since $\varphi_{\omega(m)}$ and $\psi_{\omega(m)}$ are non-negative, we have 
by \eqref{eq-4-1} that 
	\[
	- \Delta \varphi_{\omega(m)}= - \frac{\mu}{1 - \lambda} \varphi_{\omega(m)}  + 2 \varphi_{\omega(m)} \psi_{\omega(m)} 
	 \geq 0. 
	\] 
As in the previous argument, using the maximum principle, we have 
	\[
	\varphi_{\omega(m)}(x) \geq \frac{\varphi_{\omega(m)}(1)}{2 |x|^{2}} 
	\qquad \mbox{for $x \in \R^{4}$ satisfying 
	$|x| \geq 1$}.  
	\] 
This implies that $\varphi_{\omega(m)} \not\in L^{2}(\R^{4})$, 
which contradicts the fact that  $\|\varphi_{\omega(m)}\|_{L^{2}}^{2} = m$. 
Therefore, we obtain 
$\frac{\mu}{1 - \lambda} > 0$. 
Using a similar argument to the second equation 
of \eqref{eq-4-1}, we deduce that $\frac{1}{1 - \lambda} > 0$, 
which implies $\mu > 0$. 

We put $(\widetilde{\varphi}_{m}, \widetilde{\psi}_{m}) 
= (T_{\nu} \varphi_{\omega(m)}, T_{\nu} \psi_{\omega(m)})$. 
Then we see from Lemma \ref{mini-l2} 
that $(\widetilde{\varphi}_{m}, \widetilde{\psi}_{m})$ 
is also a minimizer for $\E^{\K}(m)$ and 
satisfies 
\begin{equation*} 
	\begin{cases}
	- \Delta \widetilde{\varphi}_{\omega(m)} + \frac{\mu}{1 - \lambda} \nu^{2} \widetilde{\varphi}_{\omega(m)} 
	- 2 \widetilde{\varphi}_{\omega(m)} \widetilde{\psi}_{\omega(m)} 
	 = 0
	& \qquad \mbox{in $\R^{4}$},  \\
	- \Delta \widetilde{\psi}_{\omega(m)} + \frac{1}{1 - \lambda} \nu^{2} \widetilde{\psi}_{\omega(m)} - |\widetilde{\varphi}_{\omega(m)}|^{2} 
	= 0
	& \qquad \mbox{in $\R^{4}$}.
	\end{cases}
\end{equation*}
Then, putting $\nu = \sqrt{1 - \lambda} 
\; (> 0)$, we obtain
$(\widetilde{\varphi}_{m}, \widetilde{\psi}_{m})$ satisfies \eqref{sp} 
with $\omega = \mu$. 
This completes the proof. 	
\end{proof}

\section{Asymptotic of the minimum value $\E^{\K}(m)$ and its 
minimizer as $m \to m_{*} + 0$}
\label{sec:m_*}
In this section, we consider the case where 
$m$ is sufficiently close to $m_*$ and 
shall show Theorem \ref{thm-conv}. 
\subsection{Asymptotic of the minimum value 
$\E^{\K}(m)$ and the frequency 
$\omega(m)$}
In this subsection, 
we show Theorem 
\ref{thm-conv} \textrm{(i)} and \textrm{(ii)}.  
\begin{proof}[Proof of Theorem 
\ref{thm-conv} \textrm{(i)}]
Suppose to the contrary that $\limsup_{m \to m_{*} + 0} \E^{\K}(m) < \infty$. 
Let $\{m_{n}\} \subset (0, \infty)$ be a sequence satisfying 
$\lim_{n \to \infty} m_{n} = m_{*}$ and 
$m_{n} > m_{*}$ for all $n \in \mathbb{N}$. 
Then, for each $n \in \N$, there exists a 
positive 
minimizer $(\varphi_{\omega(m_n)}, \psi_{\omega(m_n)}) 
\in H^{1}_{\text{rad}}(\R^{4}) \times H^{1}_{\text{rad}}(\R^{4})$ 
for $\E^{\K}(m_{n})$, 
which is a solution to \eqref{sp} with 
$\omega = \omega(m_n)$ 
(see Theorem \ref{Thm:ground-sol}). 
We set 
	\[
	\lambda_{n} := \left(2 
    \int_{\R^{4}} |\varphi_{\omega(m_{n})}|^2 
    \psi_{\omega(m_{n})} \, dx
	\right)^{- \frac{1}{2}}
	\]
and $(\widetilde{\varphi}_{n}, \widetilde{\psi}_{n}) 
:= (T_{\lambda_{n}} \varphi_{\omega(m_n)}, T_{\lambda_{n}} \psi_{\omega(m_n)})$. 
Then, by Lemma \ref{mini-l2}, we see that 
$(\widetilde{\varphi}_{n}, \widetilde{\psi}_{n})$ is also a minimizer of $\mathcal{E}^{\mathcal{K}}(m_n)$.
It follows from $\mathcal{K} (\widetilde{\varphi}_{n}, \widetilde{\psi}_{n}) 
= 0$ that 
	\begin{equation*} 
	1 = 2 
    \int_{\R^{4}} |\widetilde{\varphi}_{n}|^2 
    \widetilde{\psi}_{n} \, dx
	= \|\nabla \widetilde{\varphi}_{n}\|_{L^{2}}^{2} 
	+ \|\nabla \widetilde{\psi}_{n}\|_{L^{2}}^{2} 
	\geq \|\nabla \widetilde{\varphi}_{n}\|_{L^{2}}^{2}, 
    \; 
    \|\nabla \widetilde{\psi}_{n}\|_{L^{2}}^{2}. 
	\end{equation*}
This implies that $\sup_{n \in \mathbb{N}} 
\|\nabla \widetilde{\varphi}_{n}\|_{L^{2}}^{2}, 
\; \sup_{n \in \mathbb{N}} 
\|\nabla \widetilde{\psi}_{n}\|_{L^{2}}^{2}\leq 1$. 
Moreover, $\lim_{n \to \infty} 
\|\widetilde{\varphi}_{n}\|_{L^{2}}^{2} 
= \lim_{n \to \infty} m_{n} = m_{*}$ and 
from our assumption, we obtain 
    \[
    \limsup_{n \to \infty} 
    \|\widetilde{\psi}_{n}\|_{L^{2}}^{2} 
    = \limsup_{n \to \infty}\mathcal{E}^{\K}(m_{n}) 
    < \infty. 
    \]
Then, up to a subsequence of $\{(\widetilde{\varphi}_{n}, 
\widetilde{\psi}_{n})\}$ (we still denote it by the 
same letter), 
there exists $(\Phi_{m_*}, \Psi_{m_*}) \in 
H^{1}_{\text{rad}}(\R^{4}) \times 
H^{1}_{\text{rad}}(\R^{4})$ such that 
$(\widetilde{\varphi}_{n}, \widetilde{\psi}_{n}) 
\rightharpoonup (\Phi_{m_*}, \Psi_{m_*})$ 
weakly in $H^{1}_{\text{rad}}(\R^{4}) \times 
H^{1}_{\text{rad}}(\R^{4})$ as $n \to \infty$. 

It follows from the weak lower semicontinuity 
of norm and 
Lemma \ref{lem-nonc} that 
\begin{equation} \label{ex-mini-eq22}
\K(\Phi_{m_*}, \Psi_{m_*}) 
\leq \liminf_{n \to \infty} \K(\widetilde{\varphi}_{n}, \widetilde{\psi}_{n}) = 0. 
\end{equation}
As in the proof of 
Lemma \ref{lem2-3}, we put $f(\lambda) 
= \K(\Phi_{m_*}, \lambda \Psi_{m_*})$. 
Then, it follows from \eqref{ex-mini-eq22} 
that $f(1) \leq 0$ and $f(0) = 
\|\nabla \Phi_{m_*}\|_{L^{2}}^{2} > 0$. 
Hence, there exists $\lambda_{*} \in (0, 1]$ 
such that $f(\lambda_*) = 0$, that is, 
    \begin{equation} \label{eq-mini-e33}
    \mathcal{K}(\Phi_{m_*}, 
    \lambda_* \Psi_{m_*}) = 0.
    \end{equation}

By \eqref{eq-mini-e33} and the definition of $m_*$, we have
\[
m_* \leq \|\Phi_{m_*}\|_{L^2}^2.
\]
On the other hand, the weak lower semicontinuity
of the $L^2$-norm yields
\[
\|\Phi_{m_*}\|_{L^2}^2
\leq \liminf_{n\to\infty}
\|\widetilde{\varphi}_n\|_{L^2}^2
=m_*.
\]
Consequently,
\[
\|\Phi_{m_*}\|_{L^2}^2=m_*.
\]
Since $\K(\Phi_{m_*},\Psi_{m_*})\leq 0$,
it follows from \eqref{mini-1} and the weak lower
semicontinuity of the $L^2$-norm that
\begin{equation} \label{eq1-energy}
\E^{\K}(m_*)
\leq \|\Psi_{m_*}\|_{L^2}^2
\leq \liminf_{n\to\infty}
\|\widetilde{\psi}_n\|_{L^2}^2
=\lim_{n\to\infty}\E^{\K}(m_n).
\end{equation}
We next show that
\[
\E^{\K}(m_n)\leq \E^{\K}(m_*)
\]
for every $n\in\mathbb{N}$.
Indeed, by \eqref{eq-mini-e33} and
$\|\Phi_{m_*}\|_{L^2}^2=m_*$,
we have 
$(\Phi_{m_*}, 
    \lambda_* \Psi_{m_*}) \in S(m_*)$, so that 
$S(m_*)\neq\emptyset$.
For any $(u,v)\in S(m_*)$, we set
\[
\theta_n:=\left(\frac{m_n}{m_*}\right)^{1/2}>1.
\]
Then $\|\theta_n u\|_{L^2}^2=m_n$ and
\[
\begin{aligned}
\K(\theta_n u,v)
=\theta_n^2 \K(u,v)
 +(1-\theta_n^2)\|\nabla v\|_{L^2}^2 
 =(1-\theta_n^2)\|\nabla v\|_{L^2}^2
<0.
\end{aligned}
\]
Thus, by \eqref{mini-1}, we obtain
\[
\E^{\K}(m_n)\leq\|v\|_{L^2}^2=J(u,v).
\]
Taking the infimum over $(u,v)\in S(m_*)$,
we conclude that $\E^{\K}(m_n)\leq \E^{\K}(m_*)$.
Combining this with the preceding inequality 
\eqref{eq1-energy},
we obtain
\[
\E^{\K}(m_*)
=\|\Psi_{m_*}\|_{L^2}^2
=\lim_{n\to\infty}\E^{\K}(m_n).
\]
Next, we claim that $\K(\Phi_{m_*},\Psi_{m_*})=0$.
Suppose to the contrary that
$\K(\Phi_{m_*},\Psi_{m_*})<0$.
Then, by Lemma \ref{vari-chara}, we have
\[
\E^{\K}(m_*)<\|\Psi_{m_*}\|_{L^2}^2,
\]
which contradicts the equality above.
Therefore, $\K(\Phi_{m_*},\Psi_{m_*})=0$.
Hence, $(\Phi_{m_*},\Psi_{m_*})$ is a
minimizer of $\E^{\K}(m_*)$.


By the same argument as in the proof of
Theorem \ref{Thm:ground-sol}, there exists $\nu>0$ such that
$(T_\nu\Phi_{m_*},T_\nu\Psi_{m_*})$
satisfies \eqref{sp} for some $\omega>0$.
Since this scaling preserves both the
$L^2$-norm and the condition $\K=0$,
the rescaled pair is also a minimizer
of $\E^{\K}(m_*)$.
For simplicity, we continue to denote
the rescaled pair by
$(\Phi_{m_*},\Psi_{m_*})$.

Observe from the definition of $m_*$ (see \eqref{m-value}) that 
 $\K(\Phi_{m_*}, v) \geq 0$ for all $v \in H^{1}(\R^{4})$.  
This yields that 
  \[
  0 \leq \K(\Phi_{m_*}, \lambda \Psi_{m_*}) 
  = \|\nabla \Phi_{m_*}\|_{L^{2}}^{2} + \lambda^{2} 
  \|\nabla \Psi_{m_*}\|_{L^{2}}^{2} 
  - 2 \lambda 
  \int_{\R^{4}} |\Phi_{m_{*}}|^2 \Psi_{m_{*}} \, dx
  \qquad \mbox{for all $\lambda \geq 0$}. 
  \]
This together with the discriminant of the 
quadratic equation in $\lambda$
implies that 
   \begin{equation} \label{eq-5}
  0 \geq 
    \left( \int_{\R^{4}} |\Phi_{m_{*}}|^2 \Psi_{m_{*}} \, 
    dx \right)^{2} - 
  \|\nabla \Phi_{m_*}\|_{L^{2}}^{2} 
  \|\nabla \Psi_{m_*}\|_{L^{2}}^{2}.  
  \end{equation}
On the other hand, it follows from 
$\K(\Phi_{m_*}, \Psi_{m_*}) = 0$ that 
  \begin{equation} \label{eq-6}
  0 = \K(\Phi_{m_*}, \Psi_{m_*}) 
  = \|\nabla \Phi_{m_*}\|_{L^{2}}^{2} +  
  \|\nabla \Psi_{m_*}\|_{L^{2}}^{2} - 
  2 \int_{\R^{4}} |\Phi_{m_{*}}|^2 \Psi_{m_{*}} \, 
    dx.  
  \end{equation}
By \eqref{eq-5} and \eqref{eq-6}, we have 
  \[
  \begin{split}
  0 
  \geq 
  \frac{1}{4}\left(
  \|\nabla \Phi_{m_*}\|_{L^{2}}^{2} +  
  \|\nabla \Psi_{m_*}\|_{L^{2}}^{2}
  \right)^{2} - 
  \|\nabla \Phi_{m_*}\|_{L^{2}}^{2} 
  \|\nabla \Psi_{m_*}\|_{L^{2}}^{2} 
  = \frac{1}{4}\left(
  \|\nabla \Phi_{m_*}\|_{L^{2}}^{2} -  
  \|\nabla \Psi_{m_*}\|_{L^{2}}^{2}
  \right)^{2} \geq 0. 
  \end{split}
  \]
Hence, 
$\|\nabla \Phi_{m_*}\|_{L^{2}}^{2} =  
  \|\nabla \Psi_{m_*}\|_{L^{2}}^{2}$. 
  This together with 
  $\K(\Phi_{m_*}, \Psi_{m_*}) = 0$ yields that 
    \begin{equation} \label{eq-16}
    \|\nabla \Phi_{m_*}\|_{L^{2}}^{2} =  
  \|\nabla \Psi_{m_*}\|_{L^{2}}^{2} = 
 \int_{\R^{4}} |\Phi_{m_{*}}|^2 \Psi_{m_{*}} \, dx
    \end{equation}
On the other hand, since 
$(\Phi_{m_*}, \Psi_{m_*})$ satisfies \eqref{sp} for some 
$\omega = \omega(m_{*})>0$, we see from the second 
equation of \eqref{sp} that 
  \[
 -\Delta \Psi_{m_*} + \Psi_{m_*} - |\Phi_{m_*}|^{2} = 0.  
  \]
Then, multiplying the above by $\Psi_{m_*}$ and 
integrating the resulting equation, we have 
by \eqref{eq-16} that 
  \[
  0 < \|\Psi_{m_*}\|_{L^{2}}^{2} = 
\int_{\R^{4}} |\Phi_{m_{*}}|^2 \Psi_{m_{*}} \, 
    dx
  - \|\nabla \Psi_{m_*}\|_{L^{2}}^{2} 
  = 0, 
  \]
which is a contradiction. 
Thus, we infer that 
$\lim_{m \to m_{*} + 0} \E^{\K}(m) = \infty$.
\end{proof}

\begin{proof}[Proof of Theorem 
\ref{thm-conv} \textrm{(ii)}]
We shall show Theorem \ref{thm-conv} 
\textrm{(ii)} by contradiction. 
Suppose to the contrary that there exists a sequence 
$\{m_{n}\}$ with $\lim_{n \to \infty} m_{n} = m_{*}$ such that 
$\limsup_{n \to \infty} \omega(m_{n}) < \infty$. 
Let $(\varphi_{\omega(m_n)}, \psi_{\omega(m_n)}) \in H^{1}(\R^{4}) 
\times H^{1}(\R^{4})$ be the minimizer for 
$\E^{\K}(m_{n})$ satisfying \eqref{sp} 
for some $\omega = \omega(m_n)$, that is, 
\begin{equation} \label{sp-n}
	\begin{cases}
	- \Delta \varphi_{\omega(m_n)} + \omega(m_{n}) \varphi_{\omega(m_n)} 
    = 2 \varphi_{\omega(m_n)} \psi_{\omega(m_n)} 
	& \qquad \mbox{in $\R^{4}$}, \\
	- \Delta \psi_{\omega(m_n)} + \psi_{\omega(m_n)} = |\varphi_{\omega(m_n)}|^{2} 
	& \qquad \mbox{in $\R^{4}$}. 
	\end{cases}
	\end{equation}
From \eqref{sp-n}, we obtain 
  \begin{align}
  & \|\nabla \varphi_{\omega(m_n)}\|_{L^{2}}^{2} + 
  \omega(m_{n}) \|\varphi_{\omega(m_n)}\|_{L^{2}}^{2} = 
  2 \int_{\R^{4}} |\varphi_{\omega(m_{n})}|^2 \psi_{\omega(m_n)} \, dx, 
  \label{eq-7} \\
   & \|\nabla \psi_{\omega(m_n)}\|_{L^{2}}^{2} + 
  \|\psi_{\omega(m_n)}\|_{L^{2}}^{2} = 
 \int_{\R^{4}} |\varphi_{\omega(m_{n})}|^2 \psi_{\omega(m_n)} \, dx. 
   \label{eq-8}
  \end{align}
In addition, 
it follows from \eqref{mini-1} and 
Theorem \ref{thm-conv} \textrm{(i)}
that 
  \[
  \lim_{n \to \infty} \|\psi_{\omega(m_n)}\|_{L^{2}}^{2} 
  = \lim_{n \to \infty} \E^{\K}(m_{n}) = \infty. 
  \]
This together with \eqref{eq-8} yields that 
  \begin{equation} \label{eq-9}
   \lim_{n \to \infty} 
   \int_{\R^{4}} |\varphi_{\omega(m_{n})}|^2 \psi_{\omega(m_n)} \, dx
   = \infty. 
  \end{equation}
Observe from $\sup_{n \in \N} \|\varphi_{\omega(m_n)}\|_{L^{2}}^{2} 
= \sup_{n \in \N} m_{n} < \infty$, 
the assumption that $\sup_{n \in \N} \omega(m_{n}) < \infty$,  
\eqref{eq-7} and \eqref{eq-9} 
that $\lim_{n \to \infty} \|\nabla \varphi_{\omega(m_n)}\|_{L^{2}}^{2} 
= \infty$ and 
  \begin{equation} \label{eq-10}
  \lim_{n \to \infty} \dfrac{\|\nabla \varphi_{\omega(m_n)}\|_{L^{2}}^{2}}
  {\int_{\R^{4}} |\varphi_{\omega(m_{n})}|^2 \psi_{\omega(m_n)} \, dx} = 2.
  \end{equation}
Moreover, 
we see from $\K(\varphi_{\omega(m_n)}, \psi_{\omega(m_n)}) = 0$ 
that 
  \[
  0 = \dfrac{\K(\varphi_{\omega(m_n)}, \psi_{\omega(m_n)})}{\int_{\R^{4}} |\varphi_{\omega(m_{n})}|^2 \psi_{\omega(m_n)} \, dx} = 
  \left(\dfrac{\|\nabla \varphi_{\omega(m_n)}\|_{L^{2}}^{2}}
  {\int_{\R^{4}} |\varphi_{\omega(m_{n})}|^2 \psi_{\omega(m_n)} \, dx} 
  + \dfrac{\|\nabla \psi_{\omega(m_n)}\|_{L^{2}}^{2}}
  {\int_{\R^{4}} |\varphi_{\omega(m_{n})}|^2 \psi_{\omega(m_n)} \, dx} - 2\right)
  \]
This together with \eqref{eq-10} yields that 
  \begin{equation} \label{eq-11}
  \lim_{n \to \infty} \dfrac{\|\nabla \psi_{\omega(m_n)}\|_{L^{2}}^{2}}
  {\int_{\R^{4}} |\varphi_{\omega(m_{n})}|^2 \psi_{\omega(m_n)} \, dx} 
  =  \lim_{n \to \infty} \dfrac{\|\nabla \psi_{\omega(m_n)}\|_{L^{2}}^{2}}
  {\|\nabla \varphi_{\omega(m_n)}\|_{L^{2}}^{2}} 
  = 0. 
  \end{equation}
Then, by the H\"{o}lder, Gagliardo-Nirenberg and 
Sobolev inequalities, we obtain 
  \begin{equation} \label{eq-17}
	\begin{split}
	\int_{\R^{4}} |\varphi_{\omega(m_{n})}|^2 \psi_{\omega(m_n)} \, dx
	\leq  \|\varphi_{\omega(m_n)}\|_{L^{\frac{8}{3}}}^{2} 
    \|\psi_{\omega(m_n)}\|_{L^{4}} 
	& \leq \frac{1}{C_{GN}\sigma} 
    \|\nabla \varphi_{\omega(m_n)}\|_{L^{2}} \|\varphi_{\omega(m_n)}\|_{L^{2}} 
	\|\nabla \psi_{\omega(m_n)}\|_{L^{2}} \\
	& = \frac{\sqrt{m_{n}}}{C_{GN}\sigma}   
    \|\nabla \varphi_{\omega(m_n)}\|_{L^{2}} \|\nabla \psi_{\omega(m_n)}\|_{L^{2}}. 
	\end{split}
	\end{equation}
This together with \eqref{eq-10} 
and \eqref{eq-11} implies that 
  \[
  \begin{split}
  \frac{C_{GN}\sigma}{\sqrt{m_{n}}} 
  \leq 
  \dfrac{\|\nabla \varphi_{\omega(m_n)}\|_{L^{2}} \|\nabla \psi_{\omega(m_n)}\|_{L^{2}}}
  {\int_{\R^{4}} |\varphi_{\omega(m_{n})}|^2 \psi_{\omega(m_n)} \, dx} 
  = \dfrac{\|\nabla \varphi_{\omega(m_n)}\|_{L^{2}}}
  {\left(\int_{\R^{4}} |\varphi_{\omega(m_{n})}|^2 \psi_{\omega(m_n)} \, dx\right)^{\frac{1}{2}}} 
  \times 
  \dfrac{\|\nabla \psi_{\omega(m_n)}\|_{L^{2}}}
  {\left(\int_{\R^{4}} |\varphi_{\omega(m_{n})}|^2 \psi_{\omega(m_n)} \, dx\right)^{\frac{1}{2}}} 
  \to 0
  \end{split}
  \]
as $n \to \infty$, which is a contradiction. 
Thus, we conclude that $\lim_{m \to m_{*} + 0} \omega(m) = \infty$. 
\end{proof}

\subsection{Convergence to the massless ground state}
In this subsection, we shall show Theorem 
\ref{thm-conv} \textrm{(iii)}. 
Let $(\widetilde{\varphi}_{\omega(m)}), 
\widetilde{\psi}_{\omega(m)})$ be the function 
defined by \eqref{rescale-1}. 
Then, we have 
  \[
  \mathcal{J}(\varphi_{\omega(m)}, \psi_{\omega(m)}) 
  = \omega(m) \widetilde{\mathcal{J}}(\widetilde{\varphi}_{\omega(m)}, 
  \widetilde{\psi}_{\omega(m)}), \qquad 
  \mathcal{K}(\varphi_{\omega(m)}, \psi_{\omega(m)}) 
  = \omega(m) \K(\widetilde{\varphi}_{\omega(m)}, 
  \widetilde{\psi}_{\omega(m)})
  \]
where
  \begin{equation*} 
   \widetilde{\mathcal{J}}(u, v) 
   = \|\nabla u\|_{L^{2}}^{2} 
   + \|\nabla v\|_{L^{2}}^{2} 
   + \omega^{-1}(m) \|v\|_{L^{2}}^{2} 
   - 2 \int_{\R^{4}} |u|^{2} v \, dx.
  \end{equation*}
We first show the following: 
\begin{lemma} \label{lem-bounded}
$\limsup_{m \to m_{*} + 0} (\|\nabla \widetilde{\varphi}_{\omega(m)}
\|_{L^{2}}^{2} + \|\widetilde{\varphi}_{\omega(m)}
\|_{L^{2}}^{2} + \|\nabla \widetilde{\psi}_{\omega(m)}
\|_{L^{2}}^{2}) < + \infty$. 
\end{lemma}
\begin{proof}
We claim that 
  \begin{equation}\label{eq-conv6}
    \limsup_{m \to m^{*} + 0} 
    \int_{\R^{4}} |\widetilde{\varphi}_{\omega(m)}|^{2} 
    \widetilde{\psi}_{\omega(m)} \, dx < + \infty. 
  \end{equation}
To this end, we argue as in the proof 
of Theorem \ref{thm-conv} \textrm{(ii)}. 
Suppose to the contrary that 
  \[
  \lim_{m \to m_{*} + 0}  \int_{\R^{4}}|\widetilde{\varphi}_{\omega(m)}|^{2} 
  \widetilde{\psi}_{\omega(m)} \, dx = + \infty.
  \]
Observe that the pair $(\widetilde{\varphi}_{\omega(m)}, 
\widetilde{\psi}_{\omega(m)})$ satisfies 
  \begin{equation} \label{eq-conv13}
  \|\nabla \widetilde{\varphi}_{\omega(m)}\|_{L^{2}}^{2} 
  + \|\widetilde{\varphi}_{\omega(m)}\|_{L^{2}}^{2} 
  = 2 
 \int_{\R^{4}} |\widetilde{\varphi}_{\omega(m)}|^2  \widetilde{\psi}_{\omega(m)} \, dx. 
  \end{equation}
Since $\|\widetilde{\varphi}_{\omega(m)}\|_{L^{2}}^{2} 
= \|\varphi_{\omega(m)}\|_{L^{2}}^{2} = m$ 
(see \eqref{rescale-1}), we obtain 
  \begin{equation} \label{eq-conv11}
  \lim_{m \to m^{*} + 0} \frac
  { \int_{\R^{4}} |\widetilde{\varphi}_{\omega(m)}|^{2} 
    \widetilde{\psi}_{\omega(m)} \, dx}{\|\nabla 
\widetilde{\varphi}_{\omega(m)}\|_{L^{2}}^{2}} 
  = \frac{1}{2}. 
  \end{equation}
Since $\K(\widetilde{\varphi}_{\omega(m)}, 
\widetilde{\psi}_{\omega(m)}) = 0$, 
we have 
  \begin{equation} \label{eq-conv7}
  0 = \K(\widetilde{\varphi}_{\omega(m)}, \widetilde{\psi}_{\omega(m)}) 
  = \|\nabla \widetilde{\varphi}_{\omega(m)} 
  \|_{L^{2}}^{2} 
  + \|\nabla \widetilde{\psi}_{\omega(m)}\|_{L^{2}}^{2} 
  - 2\int_{\R^{4}} |\widetilde{\varphi}_{\omega(m)}|^2  \widetilde{\psi}_{\omega(m)} \, dx. 
  \end{equation}
Dividing the above equality by 
$\|\nabla \widetilde{\varphi}_{\omega(m)} 
  \|_{L^{2}}^{2}$, one has by \eqref{eq-conv11} that 
    \begin{equation} \label{eq-conv12}
    \frac{\|\nabla \widetilde{\psi}_{\omega(m)}\|_{L^{2}}^{2}}
    {\|\nabla \widetilde{\varphi}_{\omega(m)} 
  \|_{L^{2}}^{2}} 
  = -1 + 2 \dfrac{\int_{\R^{4}} |\widetilde{\varphi}_{\omega(m)}|^{2} 
    \widetilde{\psi}_{\omega(m)} \, dx}{\|\nabla \widetilde{\varphi}_{\omega(m)} \|_{L^{2}}^{2}}
    \to 0 \qquad 
    \mbox{as $m \to m^{*} + 0$}. 
    \end{equation}
As in \eqref{eq-17}, 
one has by the H\"{o}lder, 
the Gagliardo-Nirenberg and the Sobolev inequalities that 
  \[
  \begin{split}
 \int_{\R^{4}} |\widetilde{\varphi}_{\omega(m)}|^2  \widetilde{\psi}_{\omega(m)} \, dx
  \leq \|\widetilde{\varphi}_{\omega(m)}
  \|_{L^{\frac{8}{3}}}^{2} 
  \|\widetilde{\psi}_{\omega(m)}\|_{L^{4}} 
  & 
  \leq \frac{1}{C_{GN}} \|\nabla \widetilde{\varphi}_{\omega(m)}\|_{L^{2}}
  \|\widetilde{\varphi}_{\omega(m)}\|_{L^{2}} 
  \|\widetilde{\psi}_{\omega(m)}\|_{L^{4}} \\
  & \leq \frac{\sqrt{m}}{C_{GN}\sigma} \|\nabla \widetilde{\varphi}_{\omega(m)}\|_{L^{2}}
  \|\nabla \widetilde{\psi}_{\omega(m)}\|_{L^{2}}. 
  \end{split}
  \]
Dividing the above inequality by 
$\|\nabla 
\widetilde{\varphi}_{\omega(m)}\|_{L^{2}}^{2}$, 
we obtain 
  \[
  \frac{\int_{\R^{4}} |\widetilde{\varphi}_{\omega(m)}|^{2} 
    \widetilde{\psi}_{\omega(m)} \, dx}{\|\nabla 
\widetilde{\varphi}_{\omega(m)}\|_{L^{2}}^{2}}
  \leq \frac{\sqrt{m}}{C_{GN}\sigma}
  \frac{\|\nabla 
  \widetilde{\psi}_{\omega(m)}\|_{L^{2}}}
  {\|\nabla \widetilde{\varphi}_{\omega(m)}\|_{L^{2}}}
  \]
  This together with \eqref{eq-conv11} and 
  \eqref{eq-conv12} implies that 
    \[
    \frac{1}{2} 
    = \lim_{m \to m^{*} + 0}\frac{
 \int_{\R^{4}} |\widetilde{\varphi}_{\omega(m)}|^{2} 
    \widetilde{\psi}_{\omega(m)} \, dx}{\|\nabla 
\widetilde{\varphi}_{\omega(m)}\|_{L^{2}}^{2}}
  \leq  \lim_{m \to m^{*} + 0} 
 \frac{\sqrt{m}}{C_{GN}\sigma}
  \frac{\|\nabla 
  \widetilde{\psi}_{\omega(m)}\|_{L^{2}}}
  {\|\nabla \widetilde{\varphi}_{\omega(m)}\|_{L^{2}}} 
  = 0, 
    \]
which is a contradiction. 
Therefore, \eqref{eq-conv6} holds. 

We see from \eqref{eq-conv6}, 
\eqref{eq-conv13} and \eqref{eq-conv7} that 
  \[
  \limsup_{m \to m^{*} + 0} 
  \left\{ 
  \|\nabla 
  \widetilde{\varphi}_{\omega(m)}\|_{L^{2}}^{2} 
 +  \|
  \widetilde{\varphi}_{\omega(m)}\|_{L^{2}}^{2} 
  + \|\nabla 
  \widetilde{\psi}_{\omega(m)}\|_{L^{2}}
  \right\} < + \infty. 
  \]
This completes the proof.
\end{proof}

\begin{proof}[Proof of Theorem \ref{thm-conv} \textrm{(iii)}]
Let $\{m_{n}\}$ be a sequence in 
$(m_*, \infty)$ with $\lim_{n \to \infty} 
m_{n} = m_*$. 
It follows from Lemma \ref{lem-bounded} that 
the sequence 
$(\widetilde{\varphi}_{\omega(m_n)}, \widetilde{\psi}_{\omega(m_n)})$ is bounded in 
$H^{1}(\R^{4}) \times \dot{H}^{1}(\R^4)$. 
Then, there exist $(\Phi_0, \Psi_{0}) 
\in H^{1}(\R^{4}) \times \dot{H}^{1}(\R^4)$ 
such that $\lim_{n \to \infty} 
(\widetilde{\varphi}_{\omega(m_n)}, \widetilde{\psi}_{\omega(m_n)}) = 
(\Phi_0, \Psi_{0})$ weakly in 
$H^{1}(\R^{4}) \times \dot{H}^{1}(\R^4)$. 
Then, we see that $(\Phi_0, \Psi_{0})$ satisfies 
\eqref{lsp1} with $\omega = 1$. 
We shall show that $\Phi_0 \neq 0$ and $\Psi_0 \neq 0$. 
It follows from \eqref{eq-conv7} and 
the Sobolev inequality that 
    \begin{equation} \label{eq-conv9}
    \begin{split}
    \|\widetilde{\varphi}_{\omega(m_n)}
    \|_{L^{\frac{8}{3}}}^2 
    \leq \|\nabla \widetilde{\varphi}_{\omega(m_n)}
    \|_{L^{2}}^{2} + 
     \|\nabla \widetilde{\psi}_{\omega(m_n)}\|_{L^{2}}^2
    = 2 \int_{\R^{4}} 
    |\widetilde{\varphi}_{\omega(m_n)}|^{2} 
    \widetilde{\psi}_{\omega(m_n)} \, dx 
    \lesssim \|\widetilde{\varphi}_{\omega(m_n)}
    \|_{L^{\frac{8}{3}}}^2 
    \|\nabla \widetilde{\psi}_{\omega(m_n)}
    \|_{L^{2}}, 
    \end{split}
    \end{equation}
    \begin{equation} \label{eq-conv10}
    \begin{split}
    \|\nabla \widetilde{\psi}_{\omega(m_n)}\|_{L^{2}}^2 
    \leq \|\nabla \widetilde{\varphi}_{\omega(m_n)}
    \|_{L^{2}}^{2} + 
     \|\nabla \widetilde{\psi}_{\omega(m_n)}\|_{L^{2}}^2
    = 2 \int_{\R^{4}} 
|\widetilde{\varphi}_{\omega(m_n)}|^2 
    \widetilde{\psi}_{\omega(m_n)} \, dx 
    \lesssim \|\widetilde{\varphi}_{\omega(m_n)}
    \|_{L^{\frac{8}{3}}}^2 
    \|\nabla \widetilde{\psi}_{\omega(m_n)}\|_{L^{2}}.
     \end{split}
     \end{equation}
By \eqref{eq-conv10}, we have 
$\|\nabla \widetilde{\psi}_{\omega(m_n)}
\|_{L^{2}} \lesssim 
\|\widetilde{\varphi}_{\omega(m_n)}
    \|_{L^{\frac{8}{3}}}^2$. 
This together with \eqref{eq-conv9} 
yields that 
    \[
    \|\widetilde{\varphi}_{\omega(m_n)}
    \|_{L^{\frac{8}{3}}}^2 
    \lesssim |\widetilde{\varphi}_{\omega(m_n)}
    \|_{L^{\frac{8}{3}}}^2 
    \|\nabla \widetilde{\psi}_{\omega(m_n)}
    \|_{L^{2}} 
    \lesssim \|\widetilde{\varphi}_{\omega(m_n)}
    \|_{L^{\frac{8}{3}}}^4.  
    \]
Therefore, we obtain 
$\|\widetilde{\varphi}_{\omega(m_n)}
    \|_{L^{\frac{8}{3}}} \gtrsim 1$. 
Then, by the radial compactness (See Berestycki and 
P. -L.~Lions~\cite[Appendix A.II]{BeLi83}), we have 
    \[
\|\Phi_0\|_{L^{\frac{8}{3}}} 
= \lim_{n \to \infty} 
\|\widetilde{\varphi}_{\omega(m_n)}
\|_{L^{\frac{8}{3}}} \gtrsim 1.
    \]
Thus, we obtain $\Phi_0 \neq 0$. 
Since $(\Phi_0, \Psi_{0})$ satisfies 
\eqref{lsp1}, it follows from 
the second equation of 
\eqref{lsp1} that $\Psi_{0} \neq 0$. 
Thus, our claim holds. 

We see from Lemma \ref{lem-nonc} that 
$\lim_{n \to \infty} 
    \int_{\R^{4}} 
    |\widetilde{\varphi}_{\omega(m_n)}|^2 
    \widetilde{\psi}_{\omega(m_n)} \, dx = \int_{\R^{4}} |\Phi_0|^2 \Psi_0 \, dx$. 
Then, we have by \eqref{msp1} and \eqref{lsp1}
that 
    \[
    \begin{split}
    \lim_{n \to \infty} 
    \left(\|\nabla \widetilde{\varphi}_{\omega(m_n)}\|
    _{L^{2}}^{2} + \|\widetilde{\varphi}_{\omega(m_n)}
    \|_{L^{2}}^{2} \right) 
    = \lim_{n \to \infty} 
    2 \int_{\R^{4}} 
    |\widetilde{\varphi}_{\omega(m_n)}|^2 
    \widetilde{\psi}_{\omega(m_n)} \, dx 
    = 2 \int_{\R^{4}} |\Phi_0|^2 \Psi_0 \, dx 
    = \|\nabla \Phi_0\|_{L^{2}}^{2} 
    + \|\Phi_0\|_{L^{2}}^{2},  
    \end{split}
    \]
 \[
    \begin{split}
    \limsup_{n \to \infty} 
    \|\nabla \widetilde{\psi}_{\omega(m_n)}\|
    _{L^{2}}^{2} 
    = \limsup_{n \to \infty}
    \left(
     \int_{\R^{4}} 
    |\widetilde{\varphi}_{\omega(m_n)}|^2 
    \widetilde{\psi}_{\omega(m_n)} \, dx 
    - \frac{1}{\omega(m_{n})} 
    \|\widetilde{\psi}_{\omega(m_n)}
    \|_{L^{2}}^{2}
    \right)
    \leq  \int_{\R^{4}} |\Phi_0|^2 \Psi_0 \, dx 
    = \|\nabla \Psi_{0}\|_{L^{2}}^{2}.  
    \end{split}
    \]
By the weak lower semicontinuuity, we have 
    \[
    \liminf_{n \to \infty} 
    \|\nabla \widetilde{\psi}_{\omega(m_n)}\|
    _{L^{2}}^{2} \geq  \|\nabla \Psi_{0}\|_{L^{2}}^{2}.   
    \]
Therefore, we deduce that 
    \[
    \lim_{n \to \infty} 
    \|\nabla \widetilde{\psi}_{\omega(m_n)}\|
    _{L^{2}}^{2} =  \|\nabla \Psi_{0}\|_{L^{2}}^{2}, 
    \qquad 
    \lim_{n \to \infty} 
    \|\widetilde{\varphi}_{\omega(m_n)}\|
    _{H^{1}}^{2} =  \|\Phi_{0}\|_{H^{1}}^{2}.
    \]
This together with 
$\lim_{n \to \infty} 
(\widetilde{\varphi}_{\omega(m_n)}, 
\widetilde{\psi}_{\omega(m_n)}) 
= (\Phi_0, \Psi_0)$ weakly in 
$H^{1}(\R^{4}) \times \dot{H}^{1}(\R^{4})$ 
yields that 
    \[
    \lim_{n \to \infty} 
    (\widetilde{\varphi}_{\omega(m_n)}, 
    \widetilde{\psi}_{\omega(m_n)})
    = (\Phi_0, \Psi_0)
    \]
strongly in $H^{1}(\R^{4}) \times 
\dot{H}^{1}(\R^4)$. 
This concludes the proof. 
\end{proof}

\subsection{Uniqueness and non-degeneracy of the 
minimizer for $|m - m_{*}| \ll 1$}
In this subsection, we 
will prove Theorem \ref{thm-conv} \textrm{(iv)} and \textrm{(v)}.

\begin{proof}[Proof of 
Theorem \ref{thm-conv} \textrm(iv)]
 Suppose to the contrary that 
 there exist 
 sequence $\{m_n\}$ with 
 $\lim_{n \to \infty} m_{n} = m_*, m_n > m_*$ 
and  $(\widetilde{\varphi}_{\omega(m_n)}^{(i)}, 
\widetilde{\psi}_{\omega(m_n)}^{(i)}
 \; (i = 1, 2)$ such that 
 $(\widetilde{\varphi}_{\omega(m_n)}^{(1)}, \widetilde{\psi}_{\omega(m_n)}^{(1)}) \neq 
 (\widetilde{\varphi}_{\omega(m_{n})}^{(2)}, \widetilde{\psi}_{\omega(m_{n})}^{(2)})$ and 
 $(\widetilde{\varphi}_{\omega(m_n)}^{(i)}, \widetilde{\psi}_{\omega(m_n)}^{(i)})
 \; (i = 1, 2)$ is the ground state of 
 $\mathcal{E}^{\mathcal{K}}(m_n)$. 
 Then, we put 
    \[
    Z_{n} = 
    \frac{\widetilde{\varphi}_{\omega(m_{n})}^{(1)} - \widetilde{\varphi}_{\omega(m_{n})}^{(2)}}
    {\|\widetilde{\varphi}_{\omega(m_{n})}^{(1)} 
    - \widetilde{\varphi}_{\omega(m_{n})}^{(2)}\|_{H^{1}}}, 
    \qquad 
    W_{n} = 
    \frac{\widetilde{\psi}_{\omega(m_{n})}^{(1)} - \widetilde{\psi}_{\omega(m_{n})}^{(2)}}
    {\|\widetilde{\varphi}_{\omega(m_{n})}^{(1)} 
    - \widetilde{\varphi}_{\omega(m_{n})}^{(2)}
    \|_{H^{1}}}. 
    \]
Then, we see that 
$(Z_{n}, W_{n})$ satisfies $\|Z_{n}\|_{H^{1}} = 1\; 
(n \in \N)$ and 
 \begin{equation}\label{eq1:unique}
      \begin{cases}
      - \Delta Z_{n} + 
      Z_{n} 
      -2 \widetilde{\varphi}_{\omega(m_{n})}^{(1)} 
        W_{n}
      - 2 \widetilde{\psi}_{\omega(m_{n})}^{(2)} 
      Z_{n} = 0, & \\
        - \Delta W_{n} + 
        \frac{1}{\omega(m_{n})}W_{n} 
        - (\widetilde{\varphi}_{\omega(m_{n})}^{(1)} 
        + \widetilde{\varphi}_{\omega(m_{n})}^{(2)}) 
        Z_{n} = 
        0. & 
      \end{cases}
    \end{equation}
It follows from the second equation of 
\eqref{eq1:unique} that 
    \[
    \|\nabla W_{n}\|_{L^{2}}^{2} 
    + \frac{1}{\omega(m_n)} \|W_{n}\|_{L^{2}}^2
    - \int_{\R^{4}} 
    (\widetilde{\varphi}_{\omega(m_{n})}^{(1)} + \widetilde{\varphi}_{\omega(m_{n})}^{(2)})
        Z_{n}W_{n} \, dx = 0. 
    \]
By the H\"{o}lder, the Sobolev 
inequalities and $\|Z_{n}\|_{H^{1}} = 1$, 
we obtain 
    \[
    \begin{split}
    \|\nabla W_{n}\|_{L^{2}}^{2} 
    \leq \|\widetilde{\varphi}_{\omega(m_{n})}^{(1)} + \widetilde{\varphi}_{\omega(m_{n})}^{(2)}\|_{L^{\frac{8}{3}}}
    \|Z_{n}\|_{L^{\frac{8}{3}}} 
    \|W_{n}\|_{L^{4}} 
    & 
    \lesssim \|\widetilde{\varphi}_{\omega(m_{n})}^{(1)} + \widetilde{\varphi}_{\omega(m_{n})}^{(2)}\|_{L^{\frac{8}{3}}}
    \|Z_{n}\|_{H^{1}} 
    \|\nabla W_{n}\|_{L^{2}} \\
    & \leq  \|\widetilde{\varphi}_{\omega(m_{n})}^{(1)} + \widetilde{\varphi}_{\omega(m_{n})}^{(2)}\|_{L^{\frac{8}{3}}}
    \|\nabla W_{n}\|_{L^{2}} . 
    \end{split}
    \]
  This yields that 
    \begin{equation*} 
        \|\nabla W_{n}\|_{L^{2}} 
        \lesssim \|\widetilde{\varphi}_{\omega(m_{n})}^{(1)} + \widetilde{\varphi}_{\omega(m_{n})}^{(2)}\|_{L^{\frac{8}{3}}}. 
    \end{equation*}
Since $\{\widetilde{\varphi}_{\omega(m_n)}^{(i)}\}\; 
(i = 1, 2)$ is bounded in $H^{1}(\R^{4})$ (see Lemma \ref{lem-bounded}), 
we see that $\sup_{n \in \N} 
\|\nabla W_{n}\|_{L^{2}} < + \infty$. 

Since $(Z_{n}, W_{n})$ is bounded 
in $H^{1}(\R^{4}) \times \dot{H}^{1}(\R^{4})$, 
there exists $(Z_{\infty}, W_{\infty}) 
\in H^{1}(\R^{4}) \times \dot{H}^{1}(\R^{4})$ 
and a subsequence of $\{(Z_{n}, W_{n})\}$ 
(we still denote it by the same letter)
such that $\lim_{n \to \infty} 
(Z_{n}, W_{n}) = (Z_{\infty}, W_{\infty})
$ weakly in $H^{1}(\R^{4}) \times \dot{H}^{1}(\R^{4})$. 
We see from Theorem \ref{thm-conv} \textrm{(ii)} 
that $(Z_{\infty}, W_{\infty})$ satisfies 
    \begin{equation} \label{eq3:unique}
    \begin{cases}
    - \Delta Z + Z - 2 \Phi_0 W 
    - 2 \Psi_0 Z = 0, & \\
    - \Delta W  - 2 \Phi_0 Z = 0.  &
    \end{cases}    
    \end{equation}
     It follows from the non-degeneracy of 
$(\Phi_0, \Psi_0)$ (see Krieger-
Lenzmann-Raphael~\cite[Theorem 4]{MR2557200}) 
that $(Z_{\infty}, W_{\infty}) 
= (0, 0)$. 
On the other hand, we have by \eqref{eq1:unique} and $\|Z_{n}\|_{H^{1}} 
= 1$ that 
    \begin{equation} \label{eq4:unique}
    \begin{split}
    0  &= 
    \|\nabla Z_{n}\|_{L^{2}}^{2} 
    + \|Z_{n}\|_{L^{2}}^{2}
    - 2 \int_{\R^{4}} Z_{n}^2 
    \widetilde{\psi}_{\omega(m_{n})}^{(2)} \, dx 
    - 2 \int_{\R^{4}} Z_{n} W_{n} 
    \widetilde{\varphi}_{\omega(m_{n})}^{(1)} \, dx \\
    & = 1 - 2 \int_{\R^{4}} Z_{n}^2 
    \widetilde{\psi}_{\omega(m_{n})}^{(2)} \, dx 
    - 2 \int_{\R^{4}} Z_{n} W_{n} 
    \widetilde{\varphi}_{\omega(m_{n})}^{(1)} \, dx. 
    \end{split}
    \end{equation}
This together with $(Z_{\infty}, W_{\infty}) 
= (0, 0)$ and radial compactness yields that
    \[
    \biggl|\int_{\R^{4}} Z_{n}^2 
    \widetilde{\psi}_{\omega(m_{n})}^{(2)} \, dx \biggl|
    \leq \|Z_{n}\|_{L^{\frac{8}{3}}}^{2} 
    \|\widetilde{\psi}_{\omega(m_{n})}^{(2)}\|_{L^{4}} 
    \lesssim \|Z_{n}\|_{L^{\frac{8}{3}}}^{2} 
    \to 0 \quad \mbox{as $n \to \infty$}, 
    \]
    \[
    \biggl|\int_{\R^{4}} Z_{n} W_{n} 
    \widetilde{\varphi}_{\omega(m_{n})}^{(1)} \, dx \biggl|
    \leq \|Z_{n}\|_{L^{\frac{8}{3}}}
    \|W_{n}\|_{L^{4}}
    \|\widetilde{\varphi}_{\omega(m_{n})}^{(1)}\|_{L^{\frac{8}{3}}} 
    \lesssim \|Z_{n}\|_{L^{\frac{8}{3}}}
    \to 0 \quad \mbox{as $n \to \infty$}.  
    \]
These imply
    \[
    \lim_{n \to \infty} 
\int_{\R^{4}} Z_{n}^2 
    \widetilde{\psi}_{\omega(m_{n})}^{(2)} \, dx = \lim_{n \to \infty} 
    \int_{\R^{4}} Z_{n} W_{n} 
    \widetilde{\varphi}_{\omega(m_{n})}^{(1)} \, dx = 0.
    \]
Then, by \eqref{eq4:unique}, we obtain 
    \[
    0 = 1 - 
    \lim_{n \to \infty} \left\{
    2 \int_{\R^{4}} Z_{n}^2 
    \widetilde{\psi}_{\omega(m_{n})}^{(2)} \, dx 
    + 2 \int_{\R^{4}} Z_{n} W_{n} 
    \widetilde{\varphi}_{\omega(m_{n})}^{(1)} \, dx
    \right\} = 1,  
    \]
which is absurd. 
Thus, we obtain the desired result. 
\end{proof}

\begin{proof}[Proof of 
Theorem \ref{thm-conv} \textrm{(v)}]
We use an argument similar to 
that in the proof of 
Theorem \ref{thm-conv} \textrm{(iv)}.
 Suppose to the contrary that 
 there exist 
 sequence $\{m_n\}$ with 
 $\lim_{n \to \infty} m_{n} = m_*, m_n > m_*$  such that 
the minimizer 
 $(\widetilde{\varphi}_{\omega(m_n)}, \widetilde{\psi}_{\omega(m_n)})$ 
 is degenerate. 
 Then, there exists 
 $(Z_{n}, W_{n}) \neq (0, 0)$ such that 
 \begin{equation}\label{eq5:unique}
      \begin{cases}
      - \Delta Z_{n} + Z_{n} 
      - 2 Z_{n} \widetilde{\psi}_{\omega(m_n)} 
      -2 \widetilde{\varphi}_{\omega(m_n)} 
      W_{n} = 0, & \\
        - \Delta W_{n} + 
        \frac{1}{\omega(m_n)} W_{n} 
        - 2 \widetilde{\varphi}_{\omega(m_n)}
        Z_{n} = 0. & 
      \end{cases}
    \end{equation}
Here, we may assume that $\|Z_{n}\|_{H^{1}} = 1$.
Indeed, if $(Z_{n}, W_n)$ 
satisfies \eqref{eq5:unique} , so does $(\kappa Z_{n}, \kappa W_n)$ 
for any $\kappa \in \R$. In addition, if  $(Z_{n}, W_n) 
\neq (0, 0)$, we see that $Z_{n} \neq 0$ and $W_{n} \neq 0$. 
Thus, it suffices to take $\kappa_{n} = \|Z_{n}\|_{H^{1}}^{-1} (>0)$ 
and we replace$(\kappa_{n} Z_{n}, \kappa_{n} W_n)$ by 
$(Z_{n}, W_{n})$. 

The rest of the proof is similar to that of Theorem \ref{thm-conv} (iv) replaced $(\widetilde{\varphi}_{\omega(m_n)}^{(i)},\widetilde{\psi}_{\omega(m_n)}^{(i)})$ $(i = 1,2)$ with $(\widetilde{\varphi}_{\omega(m_n)},\widetilde{\psi}_{\omega(m_n)})$ and hence, we omit it.
\end{proof}

\section{Asymptotic of the minimum value $\E^{\K}(m)$ and 
its minimizer as $m \to \infty$}
\label{sec:large}
In this section, we study the case where $m>0$ is large and  
shall show Theorem \ref{thm-conv2}.
\subsection{Asymptotic of 
the minimum value $\E^{\K}(m)$ and 
the frequency 
$\omega(m)$ as $m \to \infty$}

\begin{proof}[Proof of Theorem \ref{thm-conv2} \textrm{(i)}]
\textrm{(i)}\; 
By Lemma \ref{mini-po}, we have
\[
    \E^{\K} (m)\ge \frac{\sigma^2}{4}
\]
for every $m>m^*$. Thus, it suffices to show that
\[
    \limsup_{m\to\infty}\E^{\K} (m)
    \le \frac{\sigma^2}{4}.
\]

Let $f\in C_{0}^{\infty}(\R^{4})$ be a nonnegative,
nontrivial function. For $a,b>0$, we have
\begin{align*}
    \K(af,bf^{2})
    &=a^2\left(
        \|\nabla f\|_{L^2}^2
        -2b\|f\|_{L^4}^4
      \right)
      +b^2\|\nabla(f^2)\|_{L^2}^2.
\end{align*}
Fix
\[
    b>\frac{\|\nabla f\|_{L^2}^2}
            {2\|f\|_{L^4}^4}.
\]
Since
\[
    \|\nabla f\|_{L^2}^2
    -2b\|f\|_{L^4}^4<0,
\]
we have $\K(af,bf^2)<0$ for all sufficiently
large $a>0$.

For each $m>0$, choose
\[
    a=\frac{\sqrt m}{\|f\|_{L^2}}.
\]
Then $\|af\|_{L^2}^2=m$, and $\K(af,b f^{2})<0$
whenever $m$ is sufficiently large.
Therefore, by Lemma \ref{vari-chara2}, 
\[
    \E^{\K} (m)\le\|b f^{2}\|_{L^2}^2
            =b^2\|f\|_{L^4}^4.
\]
Consequently,
\[
    \limsup_{m\to\infty}\E^{\K} (m)
    \le b^2\|f\|_{L^4}^4.
\]
Letting
\[
    b\downarrow
    \frac{\|\nabla f\|_{L^2}^2}
         {2\|f\|_{L^4}^4},
\]
we obtain
\[
    \limsup_{m\to\infty}\E^{\K} (m)
    \le
    \frac{\|\nabla f\|_{L^2}^4}
         {4\|f\|_{L^4}^4}.
\]
Since $f\in C_{0}^{\infty}(\R^{4})$ is arbitrary
and the best Sobolev constant satisfies
\[
    \sigma
    =\inf_{\substack{
        f\in C_{0}^{\infty}(\R^{4})\\
        f\ne0
    }}
    \frac{\|\nabla f\|_{L^2}^2}
         {\|f\|_{L^4}^2},
\]
it follows that
\[
    \limsup_{m\to\infty}\E^{\K} (m)
    \le\frac{\sigma^2}{4}.
\]
Combining this with the lower bound from
Lemma \ref{mini-po}, we conclude that
\[
    \lim_{m\to\infty}\E^{\K} (m)
    =\frac{\sigma^2}{4}.
\]
\end{proof}

\begin{lemma}\label{lem12:mlarge}
Let $(\varphi_{\omega(m)},\psi_{\omega(m)})$
be a minimizer for $\E^{\K}(m)$ satisfying
$\K(\varphi_{\omega(m)},\psi_{\omega(m)})=0$,
Then
\begin{align} 
& \lim_{m\to\infty}
\|\nabla\psi_{\omega(m)}\|_{L^2}^2=0, 
\label{eq16:large} \\
& \lim_{m \to \infty} 
\int_{\R^{4}}|\varphi_{\omega(m)}|^2\psi_{\omega(m)}\,dx = \frac{\sigma^{2}}{4}.
\label{eq10:mlarge}
\end{align}
\end{lemma}

\begin{proof}
Since $\K(\varphi_{\omega(m)},\psi_{\omega(m)})=0$, we have
\[
\|\nabla\varphi_{\omega(m)}\|_{L^2}^2
+\|\nabla\psi_{\omega(m)}\|_{L^2}^2
=2\int_{\R^{4}}|\varphi_{\omega(m)}|^2\psi_{\omega(m)}\,dx.
\]
On the other hand, multiplying the second
equation of the elliptic system \eqref{sp}
by $\psi_{\omega(m)}$
and integrating over $\R^{4}$, we obtain
\begin{equation} \label{eq44:mlarge}
\|\nabla\psi_{\omega(m)}\|_{L^2}^2
+\|\psi_{\omega(m)}\|_{L^2}^2
= \int_{\R^{4}}|\varphi_{\omega(m)}|^2\psi_{\omega(m)}\,dx.
\end{equation}
Consequently,
\[
\|\nabla\varphi_{\omega(m)}\|_{L^2}^2
= \|\nabla \psi_{\omega(m)}\|_{L^2}^2 + 
2\|\psi_{\omega(m)}\|_{L^2}^2.
\]
By the H\"{o}lder and Sobolev inequalities,
\begin{align*}
\|\nabla\varphi_{\omega(m)}\|_{L^2}^2
 +\|\nabla\psi_{\omega(m)}\|_{L^2}^2
=2\int_{\R^{4}}
       |\varphi_{\omega(m)}|^2\psi_{\omega(m)}\,dx
       \le
2\|\varphi_{\omega(m)}\|_{L^4}^2\|\psi_{\omega(m)}\|_{L^2} \le
 \frac{2}{\sigma}
 \|\nabla\varphi_{\omega(m)}\|_{L^2}^2
 \|\psi_{\omega(m)}\|_{L^2}.
\end{align*}
It follows from the identity that 
\begin{align*}
2\|\nabla\psi_m\|_{L^2}^2
+2\|\psi_{\omega(m)}\|_{L^2}^2 
\le
\frac{2}{\sigma}
\left(
\|\nabla\psi_{\omega(m)}\|_{L^2}^2
+2\|\psi_{\omega(m)}\|_{L^2}^2
\right)
\|\psi_{\omega(m)}\|_{L^2}.
\end{align*}
Hence,
\begin{align*}
\left(
1-\frac{\|\psi_{\omega(m)}\|_{L^2}}{\sigma}
\right)
\|\nabla\psi_{\omega(m)}\|_{L^2}^2 
\le
\|\psi_{\omega(m)}\|_{L^2}^2
\left(
\frac{2\|\psi_{\omega(m)}\|_{L^2}}{\sigma}-1
\right).
\end{align*}
Since
\begin{equation} \label{eq45:mlarge}
\|\psi_{\omega(m)}\|_{L^2}^2
=\mathcal{E}^{\mathcal{K}}(m)\longrightarrow\frac{\sigma^2}{4}
\qquad\text{as }m\to\infty,
\end{equation}
we have
\[
\frac{\|\psi_{\omega(m)}\|_{L^2}}{\sigma}
\longrightarrow\frac{1}{2}, \qquad 
\|\psi_{\omega(m)}\|_{L^2}^2
\left(
\frac{2\|\psi_{\omega(m)}\|_{L^2}}{\sigma}-1
\right)
\longrightarrow0 
\qquad \mbox{as $m \to \infty$}.
\]
Therefore,
\[
\lim_{m\to\infty}
\|\nabla\psi_{\omega(m)}\|_{L^2}^2=0.
\]
It follows from 
\eqref{lem12:mlarge}, 
\eqref{eq44:mlarge} and 
\eqref{eq45:mlarge} that 
\begin{equation*} 
\lim_{m \to \infty} 
\int_{\R^{4}}|\varphi_{\omega(m)}|^2\psi_{\omega(m)}\,dx
=
\lim_{m \to \infty} \left\{
\|\nabla\psi_{\omega(m)}\|_{L^2}^2+\E^{\K}(m) 
\right\} = \frac{\sigma^{2}}{4}.
\end{equation*}
This completes the proof.
\end{proof}

\begin{proof}[Proof of Theorem \ref{thm-conv2} \textrm{(ii)}]
We argue by contradiction.  
Suppose to the contrary that 
$\limsup_{m \to \infty} \omega(m) > 0$. 
Then, there exists a subsequence of $\{\omega(m_n)\}$ 
with $\lim_{n \to \infty} \omega(m_n) > 0$. 
It follows from the first equation of 
\eqref{sp} that 
    \begin{equation} 
    \label{eq14-large}
    \|\nabla \varphi_{\omega(m_n)}\|_{L^{2}}^{2}
    + \omega(m_m) \|\varphi_{\omega(m_n)}\|_{L^{2}}^{2}
    = 2 \int_{\R^{4}} |\varphi_{\omega(m_n)}|^2 
    \psi_{\omega(m_n)} \, dx. 
    \end{equation}
Since $m_n = 
\|\varphi_{\omega(m_n)}\|_{L^{2}}^{2}$ and 
$\liminf_{n \to \infty} \omega(m_n) > 0$, 
we have by \eqref{eq14-large} 
that $\lim_{n \to \infty} \int_{\R^{4}} |\varphi_{\omega(m_n)}|^2 
    \psi_{\omega(m_n)} \, dx = \infty$, 
    which contradicts \eqref{eq10:mlarge}. 
\end{proof} 

\subsection{Asymptotic of $L^{\infty}$-norm of  
the minimizer for $\E^{\K}(m)$ as $m \to \infty$}
To prove Theorem \ref{thm-conv2} \textrm{(iii)}, 
we need several preparations. 

It follows from \eqref{eq10:mlarge} and 
$\K(\varphi_{\omega(m)}, \psi_{\omega(m)}) 
= 0$ that 
  \begin{equation} \label{eq9:large}
     \limsup_{m \to \infty} \; 
     (\|\nabla \varphi_{\omega(m)}\|_{L^{2}}^{2} +
    \|\nabla \psi_{\omega(m)}\|_{L^{2}}^{2}) 
    < + \infty. 
    \end{equation}
Observe from \eqref{mini-1} and Lemma 
\ref{monotone} that 
\begin{equation} \label{eq36:large}
    \limsup_{m \to \infty} 
\|\psi_{\omega(m)}\|_{L^{2}}^{2} = 
\limsup_{m \to \infty} \E^{\K}(m) < \infty, 
    \end{equation}
Therefore, once we obtain \eqref{eq9:large}, 
we see that $\limsup_{m \to \infty} 
\|\psi_{\omega(m)}\|_{H^{1}} < \infty$. 

In addition, 
since we see from \eqref{eq9:large} 
that $|\varphi|^2 \in L^{2}(\R^{4})$, 
we have by the second equation for \eqref{sp} that $\psi \in H^{2, 2}
(\R^{4})$. 
Then, from the Sobolev embedding,  
there exists $C > 0$ such that 
\[
\limsup_{m \to \infty} 
\|\psi_{\omega(m)}\|_{L^{q}} \leq C
\]
for $2\leq q < \infty$.     

It follows from  
\eqref{eq9:large} that 
there exists $(\Phi, \Psi) \in 
\dot{H}^1(\R^{4}) \times H^1(\R^{4})$ such that 
$\lim_{m \to \infty} 
(\varphi_{\omega(m)}, \psi_{\omega(m)}) 
= (\Phi, \Psi)$ weakly in 
$\dot{H}^1(\R^{4}) \times H^1(\R^{4})$. 
This together with Theorem 
\ref{thm-conv2} \textrm{(ii) }yields that 
$(\Phi, \Psi)$ satisfies 
the following: 
    \begin{equation}\label{eq22:large}
    \begin{cases}
    - \Delta \Phi - 2 \Phi \Psi = 0
    & \qquad \mbox{in $\R^{4}$}, \\
    - \Delta \Psi + \Psi - |\Phi|^2 = 0 
    & \qquad \mbox{in $\R^{4}$}. 
    \end{cases}
    \end{equation}

We now in position to prove Theorem \ref{thm-conv2} (iii).

\begin{proof}[Proof of Theorem \ref{thm-conv2} (iii)]
We claim that \eqref{eq22:large} has 
no non-trivial solution 
in $\dot{H}^{1}(\R^{4}) \times H^{1}(\R^{4})$.   
It follows from Lemma \ref{lem2-pohozaev} below that 
   \begin{align} 
& \|\nabla \Phi\|_{L^{2}}^{2} 
     - 4 \int_{\R^{4}} |\Phi|^2\Psi \, dx 
    - \int_{\R^{4}} (x \cdot \nabla \Psi) |\Phi|^2 
     \, dx 
    = 0, \label{eq25:large} \\
    & \|\nabla \Psi\|_{L^{2}}^{2} 
    + 2 \|\Psi\|_{L^{2}}^2
    + \int_{\R^{4}} 
    |\Phi|^{2} (x \cdot \nabla \Psi) \, dx 
    = 0.  \label{eq24:large}
    \end{align}
By \eqref{eq25:large} and 
\eqref{eq24:large}, we obtain 
    \begin{equation} \label{eq37:large}
    \|\nabla \Phi\|_{L^{2}}^{2} 
   + \|\nabla \Psi\|_{L^{2}}^{2}
    + 2 \|\Psi \|_{L^{2}}^{2}
    - 4 \int_{\R^{4}} |\Phi|^2 \Psi \, dx = 0.
    \end{equation}
On the other hand, it follows from 
\eqref{eq22:large} that 
    \begin{align}
    & \|\nabla \Phi\|_{L^{2}}^{2} 
    - 2 \int_{\R^{4}} |\Phi|^2 \Psi \, dx
    = 0, \label{eq26:large}\\
    & \|\nabla \Psi\|_{L^{2}}^{2} 
    + \|\Psi\|_{L^{2}}^2 
    - \int_{\R^{4}} |\Phi|^2 \Psi \, dx 
    = 0. \label{eq27:large}
    \end{align}
It follows from 
\eqref{eq26:large} and \eqref{eq27:large} that 
    \[
    \|\nabla \Phi\|_{L^{2}}^{2} 
    + \|\nabla \Psi\|_{L^{2}}^{2} 
    + \|\Psi\|_{L^{2}}^2
    - 3 \int_{\R^{4}} |\Phi|^2 \Psi \, dx
    = 0.
    \]
This together with \eqref{eq37:large} 
yields that 
    \[
    \|\Psi\|_{L^{2}}^{2} -  
    \int_{\R^{4}} |\Phi|^2 \Psi \, dx = 0. 
    \]
Using \eqref{eq27:large} again, one has  
$\|\nabla \Psi\|_{L^{2}} = 0$. 
Then, we see from \eqref{eq26:large} that  
$\|\nabla \Phi\|_{L^{2}} = 0$.

We next show that
\[
\lim_{m\to\infty}
\|\varphi_{\omega(m)}\|_{L^\infty}=0.
\]
By \eqref{eq9:large}, 
\eqref{eq36:large}, and the Sobolev inequality,
we have
\[
\sup_{m\gg1}
\left(
\|\varphi_{\omega(m)}\|_{L^4}
+\|\psi_{\omega(m)}\|_{L^2}
\right)<\infty.
\]
Since
\[
-\Delta\psi_{\omega(m)}+\psi_{\omega(m)}=|\varphi_{\omega(m)}|^2,
\]
the elliptic estimate yields
\[
\|\psi_{\omega(m)}\|_{H^2}
\le C\|\varphi_{\omega(m)}\|_{L^4}^2
\le C.
\]
By the Sobolev embedding in dimension four,
we obtain
\[
\sup_{m\gg1}\|\psi_{\omega(m)}\|_{L^8}<\infty.
\]
On the other hand, the first equation of
\eqref{sp} gives
\[
-\Delta\varphi_{\omega(m)}
=(2\psi_{\omega(m)}-\omega(m))\varphi_{\omega(m)}.
\]
For any fixed $R>0$, H\"{o}lder's inequality
and Theorem \ref{thm-conv2} 
\textrm{(ii)} imply that
\begin{align*}
\|(2\psi_{\omega(m)}-\omega(m))\varphi_{\omega(m)}
  \|_{L^{8/3}(B_{2R})} 
  \le
2\|\psi_{\omega(m)}\|_{L^8}
 \|\varphi_{\omega(m)}\|_{L^4}
+\omega(m)\|\varphi_{\omega(m)}\|_{L^{8/3}(B_{2R})} 
\le C_R.
\end{align*}
Hence, by the interior elliptic
$W^{2,p}$ estimate,
\begin{align*}
\|\varphi_{\omega(m)}\|_{W^{2,8/3}(B_R)}
\le C_R\bigl(
\|\Delta\varphi_{\omega(m)}\|_{L^{8/3}(B_{2R})}
+\|\varphi_{\omega(m)}\|_{L^{8/3}(B_{2R})}
\bigr) \le C_R,
\end{align*}
where $C_R$ is independent of $m$.
From $8/3> 2$, the embedding
\[
W^{2,8/3}(B_R)
\hookrightarrow C^0(\overline{B_R})
\]
is compact. 
Therefore,
$\{\varphi_{\omega(m)}\}$ is relatively compact
in $C^0_{\mathrm{loc}}(\mathbb R^4)$.

As shown above, every weak limit of
$\{\varphi_{\omega(m)}\}$ in $\dot H^1(\mathbb R^4)$
is identically zero.
Consequently, for every sequence
$m_n\to\infty$, there exists a subsequence
such that
\[
\varphi_{m_n}\longrightarrow0
\quad\text{strongly in }
C^0_{\mathrm{loc}}(\mathbb R^4).
\]
In particular,
\[
\varphi_{m_n}(0)\longrightarrow0.
\]
Since this holds for every sequence
$m_n\to\infty$, we conclude that
\[
\lim_{m\to\infty}\varphi_{\omega(m)}(0)=0.
\]
Finally, the radial monotonicity of
$\varphi_{\omega(m)}$ implies
\[
\|\varphi_{\omega(m)}\|_{L^\infty}
=\varphi_{\omega(m)}(0)\longrightarrow0.
\]
This completes the proof.
\end{proof}

\subsection{Asymptotic of the minimum value $\E^{\K}(m)$ and 
its minimizer as $m \to \infty$}

\begin{lemma} \label{lem13:mlarge}
Let $(\varphi_{\omega(m)}, \psi_{\omega(m)}) 
\in H^{1}(\R^{4}) \times H^{1}(\R^{4})$ be a minimizer 
for $\E^{\K}(m)$ that  
satisfies \eqref{sp} with $\omega = \omega(m)$
and is radially symmetric and 
non-increasing in $r = |x|$. 
Furthermore, let 
$(\widehat{\varphi}_{\omega(m)}, 
\widehat{\psi}_{\omega(m)}) 
\in H^{1}(\R^{4}) \times H^{1}(\R^{4})$ 
be the pair of functions defined by \eqref{eq13:large}. 
Then, we obtain 
\begin{equation}\label{eq6:large}
\lim_{m \to \infty} 
\|\nabla \widehat{\varphi}_{\omega(m)}\|_{L^2}^2 = 
\frac{\sigma^{2}}{2}. 
\end{equation}
\end{lemma}
\begin{proof}
Since $\K(\varphi_{\omega(m)},\psi_{\omega(m)})=0$,
we have
\[
\|\nabla\varphi_{\omega(m)}\|_{L^2}^2
+\|\nabla\psi_{\omega(m)}\|_{L^2}^2
=
2\int_{\R^{4}}
|\varphi_{\omega(m)}|^2\psi_{\omega(m)}\,dx.
\]
On the other hand, testing the second equation
of \eqref{sp} with $\psi_{\omega(m)}$ yields
\[
\|\nabla\psi_{\omega(m)}\|_{L^2}^2
+ \|\psi_{\omega(m)}\|_{L^2}^2
=
\int_{\R^{4}}
|\varphi_{\omega(m)}|^2\psi_{\omega(m)}\,dx.
\]
Consequently,
\[
\|\nabla\varphi_{\omega(m)}\|_{L^2}^2
=\|\nabla \psi_{\omega(m)}\|_{L^2}^2 + 2 \|\psi_{\omega(m)}\|_{L^2}^2.
\]
Since the rescaling (1.20) preserves
the $\dot H^1(\R^{4})$ norm, it follows from 
Theorem \ref{thm-conv2} \textrm{(i)} and 
\eqref{eq16:large} that 
\[
\|\nabla\widehat\varphi_{\omega(m)}\|_{L^2}^2
= \|\nabla\varphi_{\omega(m)}\|_{L^2}^2
= \|\nabla \psi_{\omega(m)}\|_{L^2}^2 + 2 \|\psi_{\omega(m)}\|_{L^2}^2
= \|\nabla \psi_{\omega(m)}\|_{L^2}^2 + 
2\mathcal{E}^{\mathcal{K}}(m) 
\longrightarrow\frac{\sigma^2}{2}
\qquad\text{as }m\to\infty.
\]
\end{proof}

It follows from \eqref{eq6:large} 
that there exists 
$\widehat{\varphi}_{\infty} 
\in \dot{H}^{1}(\mathbb{R}^{4})$ 
and a subsequence 
$\{\widehat{\varphi}_{\omega(m_n)}\}$ 
(we still denote it by the same 
symbol) such that 
$\lim_{n \to \infty} 
\widehat{\varphi}_{\omega(m_n)} = 
\widehat{\varphi}_{\infty}$ 
weakly in $\dot{H^{1}}(\R^{4})$. 
Then, we obtain the following:
\begin{lemma}\label{lem9:mlarge}
$\widehat{\varphi}_{\infty} = V$, 
where $V$ is defined by 
\eqref{eq:talenti}. 
\end{lemma}
\begin{proof}
We first show that $\widehat{\varphi}_{\infty}(0)=1$.
Let $(\widehat{\varphi}_{\omega(m)}, \widehat{\psi}_{\omega(m)})$ 
be the pair of functions defined by \eqref{eq13:large}. 
By the radial monotonicity of the minimizers, 
both functions 
$\widehat{\varphi}_{\omega(m)}$ and 
$\widehat{\psi}_{\omega(m)}$ take a 
maximum at $x = 0$
and
$\Delta \widehat{\psi}_{\omega(m)}(0) \le 0$.
In addition, it follows from the definition of the rescaling that
\[
\widehat\varphi_{\omega(m)}(0)=1, \qquad
0<\widehat\varphi_{\omega(m)}(x)\le 1,\qquad
0<\widehat\psi_{\omega(m)}(x)\le 1
\]
for all $x\in\mathbb R^4$.
Moreover, evaluating the first equation of
\eqref{eq1:large} at $x=0$ and using
$\Delta \widehat\varphi_{\omega(m)}(0)\le 0$, we obtain
\begin{equation}
\label{Ineq1}
0\le \omega(m)M_m^{-2}\le 2\widehat\psi_{\omega(m)}(0)\le 2.
\end{equation}
Consequently,
\[
\|(2\widehat\psi_{\omega(m)}-\omega(m)M_m^{-2})\widehat\varphi_{\omega(m)}\|_{L^\infty(\mathbb R^4)}
\le 2\|\widehat\psi_{\omega(m)}\|_{L^\infty(\mathbb R^4)} + \omega(m)M_m^{-2}\|\widehat\varphi_{\omega(m)}\|_{L^\infty(\mathbb R^4)}
\le 4.
\]
The first equation of \eqref{eq1:large} can be written as
\[
-\Delta \widehat\varphi_{\omega(m)}=(2\widehat\psi_{\omega(m)}-\omega(m)M_m^{-2})\widehat\varphi_{\omega(m)}.
\]
Fix $R>0$ and $1<p<\infty$.
By the interior $W^{2,p}$ estimate for the
Laplace operator, we have
\begin{align*}
\|\widehat\varphi_{\omega(m)}\|_{W^{2,p}(B_R)}
\le C_{R,p}
\left(
\|\Delta \widehat\varphi_{\omega(m)}\|_{L^p(B_{2R})}
+\|\widehat\varphi_{\omega(m)}\|_{L^p(B_{2R})}
\right) 
\le C_{R,p},
\end{align*}
where $C_{R,p}$ is independent of $m$.

Choosing $p>2$, the compact Sobolev embedding
\[
W^{2,p}(B_R)\hookrightarrow C^0(\overline{B_R})
\]
implies that, after passing to a subsequence,
$\widehat\varphi_{\omega(m)}$ converges strongly in
$C^0(\overline{B_R})$.
By a diagonal argument, there exists a function
$U\in C^0(\mathbb R^4)$ such that
\[
\widehat\varphi_{\omega(m)}\longrightarrow U
\quad\text{strongly in }
C^0_{\mathrm{loc}}(\mathbb R^4).
\]
Since $\widehat\varphi_{\omega(m)}\rightharpoonup 
\widehat{\varphi}_{\infty}$ weakly in
$\dot H^1(\mathbb R^4)$, the local uniform
limit $U$ agrees with $\widehat{\varphi}_{\infty}$ almost everywhere.
Thus, identifying $\widehat{\varphi}_{\infty}$ with its continuous
representative, we obtain
\[
\widehat{\varphi}_{\infty}(0)
=\lim_{m\to\infty}\widehat\varphi_{\omega(m)}(0)
=1.
\]
In particular, $\widehat{\varphi}_{\infty}$ is nontrivial.


We take any sequence $\{m_n\}$ with 
$\lim_{n \to \infty} m_{n} = \infty$. 
Then, 
since one has $\limsup_{m \to \infty} \omega(m) 
M_{m}^{-2} < \infty$ from \eqref{Ineq1},
there exists a subsequence of 
$\{\omega(m_n) M_{m_n}^{-2}\}$ such that 
$\lim_{n \to \infty} \omega(m_n) M_{m_n}^{-2} 
    = C_0$ for some $C_0 \geq 0$ and 
$\lim_{n \to \infty} 
\widehat{\varphi}_{\omega(m_n)} 
= \widehat{\varphi}_{\infty}$ weakly in $\dot{H}^{1}
(\R^{4})$. 

We shall show that 
$\widehat{\varphi}_{\infty}$ 
is a non-trivial 
weak solution to \eqref{lsp2}. 
We take $\zeta \in C_{0}^{\infty} 
(\R^{4})$ and fix it. 
Multiplying the first and the second equation 
of \eqref{eq1:large} 
by $\zeta$  and $\widehat{\varphi}
_{\omega(m_{n})} \zeta$ and integrating 
the resulting equation, we obtain 
    \begin{equation}
    \label{eq30:large}
    (\nabla \widehat{\varphi}_{\omega(m_{n})}, 
    \nabla \zeta) + 
    \omega(m_n) M_{m_{n}}^{-2} 
    (\widehat{\varphi}_{\omega(m_{n})}, 
    \zeta) - 2 (\widehat{\varphi}_{\omega(m_n)} 
    \widehat{\psi}_{\omega(m_n)}, 
    \zeta) = 0, 
    \end{equation}
    \begin{equation}
    \label{eq31:large}
    M_{m}^2 
    (\nabla  \widehat{\psi}_{\omega(m_n)}, 
    \nabla (\widehat{\varphi}
_{\omega(m_{n})} \zeta)) 
+ (\widehat{\psi}_{\omega(m_n)} 
\widehat{\varphi}_{\omega(m_{n})}, 
\zeta) - 
(|\widehat{\varphi}_{\omega(m_{n})}|^2 
\widehat{\varphi}_{\omega(m_{n})},  
\zeta) = 0. 
    \end{equation}
Then, by 
\eqref{eq6:large}, 
\eqref{eq16:large}, \eqref{eq31:large},  $\|\widehat{\varphi}_{\omega(m_{n})}
\|_{L^{\infty}} = 1$ and the H\"{o}lder inequality, 
we obtain 
    \begin{equation*}
    \begin{split}
    & \quad \biggl| 
    (\widehat{\psi}_{\omega(m_n)} 
\widehat{\varphi}_{\omega(m_{n})}, 
\zeta) - 
(|\widehat{\varphi}_{\omega(m_n)}|^2 
\widehat{\varphi}_{\omega(m_n)},  
\zeta)\biggl| \\
& = \biggl|  M_{m_n}^2 
    (\nabla  \widehat{\psi}_{\omega(m_n)}, 
    \nabla (\widehat{\varphi}
_{\omega(m_n)} \zeta)) \biggl| \\
& \leq M_{m_n}^2
 \|\nabla 
 \widehat{\psi}_{\omega(m_n)}\|_{L^2}
 \|\nabla \widehat{\varphi}
_{\omega(m_n)}\|_{L^{2}} 
\|\zeta\|_{L^{\infty}} 
+ M_{m_n}^2
 \|\nabla 
 \widehat{\psi}_{\omega(m_n)}\|_{L^2}
 \|\widehat{\varphi}
_{\omega(m_n)}\|_{L^{\infty}} 
\|\nabla \zeta\|_{L^{2}} \\
& \lesssim 
M_{m_n}
\left(M_{m_n}^2
 \|\nabla 
 \widehat{\psi}_{\omega(m_n)}\|_{L^2}^2 
 \right)^{\frac{1}{2}}
 \to 0 \quad \mbox{as $n \to \infty$}. 
  \end{split}
     \end{equation*}
By the Rellich-Kondrachov theorem, one 
has 
    \[
    \biggl| 
    (
\widehat{\varphi}_{\infty}^3, 
\zeta) - 
(|\widehat{\varphi}_{\omega(m)}|^2 
\widehat{\varphi}_{\omega(m)},  
\zeta)\biggl| \to 0 \quad 
    \mbox{as $n \to \infty$}. 
    \]
These imply that 
     \[
    \biggl| 
    (\widehat{\varphi}_{\infty}^3, 
\zeta) - (\widehat{\psi}_{\omega(m_n)} 
\widehat{\varphi}_{\omega(m)}, 
\zeta) \biggl| \to 0 \quad 
    \mbox{as $n \to \infty$}. 
    \]
Then, 
we see that $\widehat{\varphi}_{\infty}$ is a weak solution to 
the following elliptic equation: 
    \begin{equation} \label{eq45:large}
    - \Delta \widehat{\varphi}_{\infty} + C_0 
    \widehat{\varphi}_{\infty} - 2 |\widehat{\varphi}_{\infty}|^2 
    \widehat{\varphi}_{\infty} = 0 
    \qquad \mbox{in $\R^{4}$}. 
    \end{equation}
Since $\widehat{\varphi}_{\infty} \neq 0$, we have $C_0 = 0$. 
Indeed, suppose that $C_0>0$. Since
$\widehat\varphi_\infty\in\dot H^1(\mathbb R^4)$,
a standard cutoff argument applied to \eqref{eq45:large}
shows that $\widehat\varphi_\infty\in L^2(\mathbb R^4)$.
The Pohozaev identity and the identity obtained
by testing \eqref{eq45:large} with $\widehat\varphi_\infty$
then yield
\[
C_0\|\widehat\varphi_\infty\|_{L^2}^2=0,
\]
contradicting $\widehat\varphi_\infty(0)=1$.
Hence $C_0=0$.

Therefore, we see that 
$\widehat{\varphi}_{\infty}$ satisfies 
\eqref{lsp2}. 
Namely, $\widehat{\varphi}_{\infty}$ is a 
non-trivial weak positive solution to \eqref{lsp2} 
with $\widehat{\varphi}_{\infty}(0) = 1$. 
From the uniqueness of positive solution to \eqref{lsp2} 
with $V(0) = 1$, we 
have obtained the desired result. 
\end{proof}

\begin{proof}[Proof of Theorem \ref{thm-conv2} 
\textrm{(iv)}]
By the preceding result and \eqref{eq16:large}, 
we have
\[
\|\nabla\widehat\varphi_{\omega(m)}\|_{L^2}^2
= \|\nabla\widehat \psi_{\omega(m)}\|_{L^2}^2 
+ 2\E^{\K}(m)
\longrightarrow \frac{\sigma^2}{2}
=\|\nabla V\|_{L^2}^2
\qquad\text{as }m\to\infty.
\]
On the other hand, Lemma \ref{lem9:mlarge} implies that,
for any sequence $m_n\to\infty$, there exists
a subsequence, still denoted by $\{m_n\}$,
such that
\[
\widehat\varphi_{\omega(m_n)}
\rightharpoonup V
\quad\text{weakly in }\dot H^1(\R^{4}).
\]
Consequently,
\begin{align*}
\|\nabla(\widehat\varphi_{\omega(m_n)}-V)\|_{L^2}^2 =\|\nabla\widehat\varphi_{\omega(m_n)}\|_{L^2}^2
+\|\nabla V\|_{L^2}^2 
-2\int_{\R^{4}}
\nabla\widehat\varphi_{\omega(m_n)}
\cdot\nabla V\,dx
\longrightarrow0.
\end{align*}
Thus,
\[
\widehat\varphi_{\omega(m_n)}
\longrightarrow V
\quad\text{strongly in }\dot H^1(\R^{4}).
\]
Since the weak limit $V$ is uniquely determined
and the above argument applies to every sequence
$m_n\to\infty$, the convergence holds for the
whole family as $m\to\infty$.
This completes the proof. 
\end{proof}

\section{Analysis on 
the potential well $\mathcal{R}_{-}(m)$}
\label{sec:potential}

\begin{lemma} \label{lem1-invariant}
Let $m > m_{*}$. 
For any $\varepsilon > 0$, 
there exist 
$\lambda_{\e, \pm} \in (1 - \e, 1 + \e)$ and $\theta_{\e, \pm} \in (1- \e, 1 + \e)$ such that 
$(T_{\lambda_{\e, \pm}}
     \varphi_{\omega(m)}, T_{\lambda_{\e, \pm}}
     \theta_{\e, \pm} \psi_{\omega(m)}, 0) \in 
    \mathcal{R}_{\pm}(m)$. 
As a result, we have 
$\mathcal{R}_{\pm}(m) \neq \emptyset$. 
\end{lemma}
\begin{proof}
We first show that $\mathcal{R}_{-}(m) \neq \emptyset$. 
We take a minimizer 
$(\varphi_{\omega(m)}, 
\psi_{\omega(m)})$ of $\E^{\K}(m)$ 
satisfying \eqref{sp} with $\omega = \omega(m)$. 
For $\theta > 0$, we put 
    \begin{equation} \label{eq4:inv}
    \begin{split}
     f_m(\theta) := 
     \K(
     \varphi_{\omega(m)}, 
     \theta \psi_{\omega(m)}) 
=  
\|\nabla \varphi_{\omega(m)}
\|_{L^{2}}^2 + \theta^2 \|\nabla \psi_{\omega(m)}\|_{L^{2}}^2 
- 2 \theta \int_{\R^{4}} 
|\varphi_{\omega(m)}|^2 
\psi_{\omega(m)} \, dx. 
    \end{split}
    \end{equation}
We have $f_m(1) =  \K(
     \varphi_{\omega(m)}, 
     \psi_{\omega(m)})= 0$. 
It follows from the second equation of 
\eqref{sp} that  
    \begin{equation} \label{eq5:inv} 
    f_{m}^{\prime}(1)
    = 2\|\nabla \psi_{\omega(m)}\|_{L^{2}}^2 
    -2 \int_{\R^{4}} 
|\varphi_{\omega(m)}|^2 
\psi_{\omega(m)} \, dx 
= -2 \|\psi_{\omega(m)}\|_{L^{2}}^{2}. 
    \end{equation}
Since $f_{m}(1) = 0, f_{m}'(1) < 0$ and 
$f$ is a quadratic function, 
there exists $\alpha_{m} > 1$ such that 
$f_{m}(\alpha_{m}) = 0$ and 
\begin{equation} \label{eq6:inv}
f_{m}(\theta) < 0 \quad  
\mbox{for $\theta \in (1, \alpha_{m})$}.
\end{equation}
By $f_{m}(1) = 0$, \eqref{eq4:inv} and 
\eqref{eq5:inv}, we obtain 
 \begin{equation} \label{eq1:inv}
 \begin{split}
f_{m}(\theta) 
& = f_{m}(\theta) - f_{m}(1) = 
-2 (\theta - 1) \int_{\R^{4}} |\varphi_{\omega(m)}|^{2} 
\psi_{\omega(m)} \, dx 
+ (\theta^{2} - 1) \|\nabla \psi_{\omega(m)}\|_{L^{2}}^{2} \\
& = (\theta - 1) \left( -2\int_{\R^{4}} |\varphi_{\omega(m)}|^{2} 
\psi_{\omega(m)} \, dx 
+ (\theta + 1) \|\nabla \psi_{\omega(m)}\|_{L^{2}}^{2} 
\right) \\
& = (\theta - 1) \left( -2\int_{\R^{4}} 
|\varphi_{\omega(m)}|^{2} \psi \, dx 
+ 2 \|\nabla \psi_{\omega(m)}\|_{L^{2}}^{2} 
+ (\theta - 1) \|\nabla \psi_{\omega(m)}\|_{L^{2}}^{2}
\right) \\
& = -2 (\theta - 1) \|\psi_{\omega(m)}\|_{L^{2}}^{2} 
+ (\theta - 1)^{2} \|\nabla \psi_{\omega(m)}\|_{L^{2}}^{2}. 
 \end{split}
\end{equation}
    It follows from $\mathcal{J}(u, v) 
    = \mathcal{K}(u, v) + \|v\|_{L^{2}}^{2}$ that 
    \begin{equation} \label{eq2:inv}
    \begin{split}
    \mathcal{J}(T_{\lambda}
     \varphi_{\omega(m)}, T_{\lambda}
     \theta \psi_{\omega(m)}) 
    & = \mathcal{K}(T_{\lambda}
     \varphi_{\omega(m)}, T_{\lambda}
     \theta \psi_{\omega(m)})
    + \theta^2 
    \|T_{\lambda} \psi_{\omega(m)}\|_{L^{2}}^{2} \\
   & = 
   \lambda^{2} \mathcal{K}(
     \varphi_{\omega(m)}, 
     \theta \psi_{\omega(m)})
   + \theta^2 
    \|\psi_{\omega(m)}\|_{L^{2}}^{2} \\
    & = \lambda^{2} f_{m}(\theta)
   + (\theta^2 - 1) 
    \|\psi_{\omega(m)}\|_{L^{2}}^{2} 
    + \|\psi_{\omega(m)}\|_{L^{2}}^{2} 
     \end{split}
     \end{equation}
for any $\lambda > 0$. 
By \eqref{eq1:inv}, we obtain 
    \begin{equation} \label{eq3:inv}
    \begin{split}
    \lambda^{2} f_{m}(\theta)
   + (\theta^2 - 1) 
    \|\psi_{\omega(m)}\|_{L^{2}}^{2} 
    & 
    = -2 (\theta -1) \lambda^{2} \|\psi_{\omega(m)}\|_{L^{2}}^{2}
+ \lambda^{2} (\theta - 1)^2
 \|\nabla \psi_{\omega(m)}\|_{L^{2}}^{2}
+ (\theta^2 - 1) 
    \|\psi_{\omega(m)}\|_{L^{2}}^{2} \\
    & = 
    (\theta - 1) \left( 
    (-2 \lambda^{2} + \theta + 1)
    \|\psi_{\omega(m)}\|_{L^{2}}^{2}
    + \lambda^{2} (\theta - 1)
 \|\nabla \psi_{\omega(m)}\|_{L^{2}}^{2} \right). 
    \end{split}
    \end{equation}
We take $\lambda_{\e, -} \in (1, 1 + \e)$ and fix it. 
For such $\lambda_{\e, -} > 1$, 
we set
    \[
    F_{\lambda_{\e, -}}(\theta) 
    = \lambda_{\e, -}^2 
    f_{m}(\theta) + (\theta^{2} - 1) 
    \|\psi_{\omega(m)}\|_{L^{2}}^{2}. 
    \]
Then, we have $F_{\lambda_{\e, -}}(1) = 0$
and $F_{\lambda_{\e, -}}^{\prime}(1) 
= 2(1- \lambda_{\e, -}^2) 
\|\psi_{\omega(m)}\|_{L^{2}}^2 < 0$. 
Thus, 
we can take 
$\theta_{\e, -} \in (1, \min\{1 + \e, \alpha_m\})$ 
so that  
$F_{\lambda_{\e, -}}(\theta_{\e, -}) < 0$. 
Then, it follows from \eqref{eq3:inv} that 
    \[
    \lambda_{\e, -}^{2} f_m(\theta_{\e, -})
   + (\theta_{\e, -}^2 - 1) 
    \|\psi_{\omega(m)}\|_{L^{2}}^{2} < 0.
    \]
 This together with $\E^{\K}(m) 
= \|\psi_{\omega(m)}\|_{L^{2}}^{2}$ 
(see \eqref{mini-1}) and \eqref{eq2:inv} 
yields that 
    \[
    \mathcal{J}(T_{\lambda_{\e, -}}
     \varphi_{\omega(m)}, T_{\lambda_{\e, -}}
     \theta_{\e, -} \psi_{\omega(m)}) 
      = \lambda_{\e, -} ^{2}
    f_m(\theta_{\e, -})
   + (\theta_{\e, -}^2 - 1) 
    \|\psi_{\omega(m)}\|_{L^{2}}^{2} 
    + \|\psi_{\omega(m)}\|_{L^{2}}^{2} 
     < \E^{\K}(m).
     \]
In addition, observe from \eqref{eq6:inv} and 
$\theta_{\e, -} \in (1, 1 + \e)$ that $\|T_{\lambda_{\e, -}}
     \varphi_{\omega(m)}\|_{L^{2}}^{2} 
     = \|
     \varphi_{\omega(m)}\|_{L^{2}}^{2}
    = m$ and 
    \[
    \begin{split} 
    \mathcal{K}(T_{\lambda_{\e, -}}
     \varphi_{\omega(m)}, T_{\lambda_{\e, -}}
     \theta_{\e, -} \psi_{\omega(m)})
   = \lambda_{\e, -}^{2} \mathcal{K}(
     \varphi_{\omega(m)}, 
     \theta_{\e, -} \psi_{\omega(m)})
   = \lambda_{\e, -}^{2} f_m(\theta_{\e, -}) < 0.  
    \end{split}
    \]
This implies $(T_{\lambda_{\e, -}}
     \varphi_{\omega(m)}, T_{\lambda_{\e, -}}
     \theta_{\e, -} \psi_{\omega(m)}, 0) \in 
    \mathcal{R}_{-}(m)$. 

Similarly, we will show that  
$\mathcal{R}_{+}(m) \neq \emptyset$. 
Observe that 
$f_{m}(\theta) > 0$ for $\theta \in (0, 1)$. 
We take $\lambda_{\e, +} \in 
(1 - \e, 1)$ and fix it. 
For such $\lambda_{\e, +} \in 
(1 - \e, 1)$, 
\[
    F_{\lambda_{\e, +}}(\theta) 
    = \lambda_{\e, +}^2 
    f_{m}(\theta) + (\theta^{2} - 1) 
    \|\psi_{\omega(m)}\|_{L^{2}}^{2}. 
    \]
It follows that $  F_{\lambda_{\e, +}}(1) = 0$
and $F_{\lambda_{\e, +}}^{\prime}(1) 
= 2(1 - \lambda_{\e, +}^2)\|\psi_{\omega(m)}\|_{L^{2}}^{2} > 0$. 
Thus, 
we can take $\theta_{\e, +} \in (1 - \e, 1)$ 
such that 
$F_{\lambda_{\e, +}}(\theta_{\e, +}) < 0$ 
and $f_{m}(\theta_{\e, +}) > 0$.
    It follows from \eqref{eq2:inv} and 
    $\E^{\K}(m) 
= \|\psi_{\omega(m)}\|_{L^{2}}^{2}$ that 
    \[
    \begin{split}
    \mathcal{J}(T_{\lambda_{\e, +}}
     \varphi_{\omega(m)}, T_{\lambda_{\e, +}}
     \theta_{\e, +} \psi_{\omega(m)}) 
    = \lambda_{\e, +} ^{2}
    f_m(\theta_{\e, +})
   + (\theta_{\e, +}^2 - 1) 
    \|\psi_{\omega(m)}\|_{L^{2}}^{2} 
    + \|\psi_{\omega(m)}\|_{L^{2}}^{2} 
    < \E^{\K}(m). 
     \end{split}
     \]
In addition, observe that $\|T_{\lambda_{\e, +}}
     \varphi_{\omega(m)}\|_{L^{2}}^{2} 
     = \|
     \varphi_{\omega(m)}\|_{L^{2}}^{2}
    = m$ and 
    \[
    \begin{split} 
    \mathcal{K}(T_{\lambda_{\e, +}}
     \varphi_{\omega(m)}, T_{\lambda_{\e, +}}
     \theta_{\e, +} \psi_{\omega(m)})
   = \lambda_{\e, +}^{2} \mathcal{K}(
     \varphi_{\omega(m)}, 
     \theta_{\e, +} \psi_{\omega(m)})
   = \lambda_{\e, +}^{2} f_m(\theta_{\e, +}) > 0.  
    \end{split}
    \]
Hence, 
$(T_{\lambda_{\e, +}}
     \varphi_{\omega(m)}, T_{\lambda_{\e, +}}
     \theta_{\e, +} \psi_{\omega(m)}, 0) \in 
    \mathcal{R}_{+}(m)$.     
This completes the proof. 
\end{proof}
We see that the sets 
$\mathcal{R}_{\pm}(m)$ are invariant under the flow \eqref{kgs}. 
Namely, we have the following: 
\begin{lemma}\label{thm-inv}
Let $(u_{0}, v_{0}, w_{0}) \in \mathcal{R}_{\pm}(m)$. 
Then, we have $(u(t), v(t), \partial_{t} v(t)) 
\in \mathcal{R}_{\pm}(m)$ 
for all $t \in [0, T_{\max})$, 
where $(u(t), v(t))$ is a solution to \eqref{kgs} 
for $t \in [0, T_{\max})$ with $(u(t), v(t), \p_{t} v(t))|_{t = 0} 
= (u_{0}, v_{0}, w_{0})$
and $T_{\max}$ denotes the 
maximal existence time of the solution.
\end{lemma}
\begin{proof}
Let $(u_{0}, v_{0}, w_{0}) \in \mathcal{R}_{+}(m)$. 
It follows from the conservation laws \eqref{eq:conserv} 
that 
	\[
	\M(u(t)) = \M(u_{0}) = m, 
	\]  
\begin{equation}\label{eq-1}
\E^{\K} (m) 
> \E(u_{0}, v_{0}, w_0) 
= \E(u(t), v(t), \partial_{t} v(t)) 
\geq \mathcal{J}(u(t), v(t))
\end{equation}
for all $t \in [0, T_{\max})$. 
We shall show that $\K(u(t), v(t)) > 0$ 
for all $t \in [0, T_{\max})$. 
Suppose that there exists $t_{1} \in (0, T_{\max})$ such that 
$\K(u(t_{1}), v(t_{1})) \leq 0$. 
Then, from the continuity of the function $\K(u(t), v(t))$, 
it follows that there exists $t_{*} \in (0, t_{1}]$ such that 
$\K(u(t_{*}), v(t_{*})) = 0$. 
From the definition of $\E^{\K}(m)$, we have 
$\E^{\K}(m) \leq \mathcal{J}(u(t_{*}), v(t_{*}))$, 
which contradicts \eqref{eq-1}. 
Thus, we see that $\mathcal{R}_{+}(m)$ 
is an invariant set. 
Similarly, we can prove that 
$\mathcal{R}_{-}(m)$ 
is also invariant under the flow of \eqref{kgs}. 
\end{proof}
Next, we show the following: 
\begin{lemma} \label{lem-uni-neg}
Let $(u,v,\partial_t v)\in C([0, T_{\max}), X)$ 
be a solution to \eqref{kgs} with 
$(u,v,\partial_t v)|_{t = 0} \in \mathcal{R}_{-}(m)$.  
There exists $\delta >0$ such that 
    \[
    \sup_{t \in [0,T_{\max})}
\mathcal{K}(u(t),v(t))\leq - \delta.
    \]
\end{lemma}
\begin{proof}
We use a contradiction argument.
Suppose to the contrary that 
there exists $t_n\in [0,T_{\max})$ such that we get
\begin{align*}
-\frac{1}{n} < \mathcal{K}(u(t_n),v(t_n)) < 0.
\end{align*}
Therefore, there exists $\{t_n\}_{n\in\mathbb N} 
\subset [0, T_{\max})$ such that $\lim_{n\to\infty} 
\mathcal{K}(u(t_n),v(t_n))=0$.
We take
\begin{equation} \label{eq5:non-ub}
\lambda_n := \frac{\|\nabla u(t_n)\|_{L^2}^2 + \|\nabla v(t_n)\|_{L^2}^2}{2 \int_{\R^{4}} |u(t_n)|^2 v(t_n) \, dx}.
\end{equation}
Then, it follows that 
\begin{align*}
\mathcal{K}
(\lambda_n u(t_n), \lambda_n v(t_n)) 
= \lambda_{n}^{2} 
\left( 
\|\nabla u(t_n)\|_{L^2}^2 + \|\nabla v(t_n)\|_{L^2}^2 
- 2 \lambda_{n}  \int_{\R^{4}} |u(t_n)|^2 v(t_n) \, dx
\right) 
= 0.
\end{align*}
Moreover, since 
$\mathcal{K}(u(t_n), v(t_n)) < 0$, 
we have $0 < \lambda_n < 1$ for each 
$n \in \mathbb{N}$.
On the other hand, 
by using $\mathcal{K}(u(t_n),v(t_n)) < 0$, 
we see that 
\begin{equation} \label{eq4:non-ub}
0 \leq \|\nabla u(t_n)\|_{L^2}^2 + \|\nabla v(t_n)\|_{L^2}^2 < 2 \int_{\R^{4}} |u(t_n)|^2 v(t_n) \, dx. 
\end{equation}

We claim that there exists $C_0 > 0$, 
which is independent of $n \in \mathbb{N}$, 
such that 
    \begin{equation} \label{eq3:non-ub}
    \liminf_{n \to \infty} 
\int_{\R^{4}} |u(t_n)|^2 v(t_n) \, dx 
\geq C_0. 
    \end{equation}   
Suppose to the contrary that 
there exists a subsequence of $\{(u(t_{n}), v(t_{n}))\}$ 
(we still denote it by the same symbol) such that 
$\lim_{n \to \infty} 
\int_{\R^{4}} |u(t_n)|^2 v(t_n) \, dx 
= 0$. 
By \eqref{eq4:non-ub}, we obtain 
$\lim_{n \to \infty} 
(\|\nabla u(t_n)\|_{L^2}^2 + \|\nabla v(t_n)\|_{L^2}^2) = 0$. 
This together with the energy conservation laws 
\eqref{eq:conserv} yields that 
\begin{equation}\label{eq1:non-ub}
\mathcal{E}^{\mathcal{K}}(m) 
> \E(u_0, v_0, w_0) 
= \lim_{n \to \infty} 
\mathcal{E}(u(t_n), v(t_{n}), \partial_{t} v(t_n)) 
= \lim_{n \to \infty} 
\left( 
\|v(t_{n})\|_{L^{2}}^{2} + 
\|\partial_{t} v(t_n)\|_{L^{2}}^{2}
\right). 
\end{equation}
On the other hand, since 
$\mathcal{K}(\lambda_{n}u(t_n), \lambda_n
v(t_n)) = 0$ and $0 < \lambda_{n} < 1$, 
we have $\lambda_{n}^2 m \geq m_*$.  
In addition, Lemmas \ref{monotone} and
the definition of $\mathcal{E}^{\mathcal{K}}(m)$ 
yield
\begin{equation}\label{eq2:non-ub}
\begin{split}
\mathcal{E^{\mathcal{K}}}(m)
\leq \mathcal{E^{\mathcal{K}}}(\lambda_n^2 m)  
\leq \mathcal{J}(\lambda_n u(t_n), 
\lambda_n v(t_n)). 
\end{split}
\end{equation}
Using $0 < \lambda_{n} < 1$,  
$\lim_{n \to \infty} 
\|\nabla u(t_n)\|_{L^2}^2 + \|\nabla v(t_n)\|_{L^2}^2 = 0$ and 
\eqref{eq1:non-ub}, we obtain 
    \[
    \limsup_{n \to \infty} 
    \mathcal{J}(\lambda_n u(t_n), 
    \lambda_n v(t_n))
    \leq \limsup_{n \to \infty} \lambda_{n}^2
    \|v(t_{n})\|_{L^{2}}^{2} 
    \leq \limsup_{n \to \infty} 
    \|v(t_{n})\|_{L^{2}}^{2} 
    < \mathcal{E}^{\mathcal{K}}(m), 
    \]
which contradicts \eqref{eq2:non-ub}. 
Thus, \eqref{eq3:non-ub} holds. 


On the other hand, 
since $\lim_{n \to \infty}
\mathcal{K}(u(t_{n}), v(t_{n})) = 0$
and $0 < \lambda_{n} < 1$, we obtain
\begin{align*}
& \quad  
\liminf_{n \to \infty}
\mathcal{E^{\mathcal{K}}}(\lambda_n^2 m) 
\leq \liminf_{n \to \infty} 
\mathcal{J}(\lambda_n u(t_n), 
\lambda_n v(t_n)) 
= \liminf_{n \to \infty} 
\lambda_{n}^2 \|v(t_n)\|_{L^{2}}^{2}
 \\
&  
\leq \liminf_{n \to \infty} 
\|v(t_n)\|_{L^{2}}^{2}
= \liminf_{n \to \infty} 
\mathcal{J}(u(t_{n}), v(t_n)) 
\leq 
\liminf_{n \to \infty}
\mathcal{E}(u(t_n), v(t_{n}), \partial_t 
v(t_n)) \\
& = \mathcal{E}(u_0, v_0, w_0)
< \mathcal{E}^{\mathcal{K}}(m),
\end{align*}
which contradicts \eqref{eq2:non-ub}.
Thus, we have obtained the desired result. 
\end{proof}

We define $\eta \in C^{\infty}([0,\infty))$ as a non-negative function satisfying
\begin{align*}
\eta(r)=
\begin{cases}
4 & \mbox{for } 0 \leq r \leq 1, \\ 
0 & \mbox{for } r \geq 2, 
\end{cases}
\quad 
\eta'(r) \leq 0 
\quad \mbox{for } 1 \leq r \leq 2. 
\end{align*}
For each $n > 0$, we define $\eta_{n}: \R^{4} \to \R$ and $\zeta_{n}: \R^{4} \to \R$ by
\begin{align*}
\eta_{n}(r) = \eta\left(\frac{r}{n}\right), 
\quad 
\zeta_{n}(r) = \frac{1}{r^{4}} \int_{0}^{r} s^{3} \eta_{n}(s) \, ds. 
\end{align*}
The functions 
$\eta_{n}$ and $\zeta_{n}$ satisfy the 
following properties: 
\begin{lemma}[c.f. Lemma 7 of Ohta and Todorova~\cite{MR2299435}]\label{lem-cutoff}
For each $n > 0$, 
the following assertions hold:
\begin{enumerate}
    \item[\textrm{(i)}] 
    $\eta_{n}(r) = 4$ $(0 \leq r \leq n)$, \; 
    $0 \leq \eta_{n}(r) \leq 4 \; (r \geq 0)$, 
    \item[\textrm{(ii)}]
    $\zeta_{n}(r) = 1$ $(0 \leq r \leq n)$, \; 
    $0 \leq \zeta_{n}(r) \leq 1 \; (r \geq 0)$
     \item[\textrm{(iii)}]
    \[
    r \zeta_{n}^{\prime}(r) + 
    4 \zeta_{n}(r) = \eta_{n}(r), 
    \]
    \item[\textrm{(iv)}]
    $\zeta_{n}^{\prime}(r) \leq 0, \quad
    x \cdot \nabla \zeta_{n}(x) \leq 0 \quad
    \quad (r \geq 0)$, 
    \item[\textrm{(v)}]
    $|\zeta_{n}(r)| \leq 
    \min\{1, \frac{16n^4}{r^{4}}\}$, 
    \item[\textrm{(vi)}]
    $|\Delta \eta_{n}(x)| \leq \frac{C}{n^{2}}, \quad
    |\Delta \zeta_n(x)| \leq \frac{C}{n^2}, \quad  
    |[x\cdot\nabla(\Delta \zeta_n)](x)| \leq \frac{C}{n^2} 
    \quad (r \geq 0)$. 
\end{enumerate}
\end{lemma}
\begin{proof}
 Assertions \textrm{(i)}--\textrm{(iii)}, and \textrm{(vi)} 
  follows directly from the definitions. 
We now prove \textrm{(iv)}. 
Integrating by parts, we have by $\eta_{n}^{\prime}(r) 
\leq 0$ that 
     \[
    \int_{0}^{r} 4 s^{3} \eta_{n}(s) \, ds 
    = \int_{0}^{r} (s^{4})^{\prime} \eta_{n}(s) \, ds 
    = r^{4} \eta_{n}(r) - \int_{0}^{r} 
    s^{4} \eta_{n}^{\prime}(s) \, ds 
    \geq r^{4} \eta_{n}(r). 
    \]
This implies that 
    \[
    4 \zeta_{n}(r) = \frac{1}{r^{4}} 
     \int_{0}^{r} 4 s^{3} \eta_{n}(s) \, ds 
     \geq \eta_{n}(r). 
    \]
Thus, we obtain $r \zeta_{n}^{\prime}(r) 
    = \eta_{n}(r) 
    - 4 \zeta_{n}(r) 
    \leq 0$. 
    
Next, we prove \textrm{(v)}.     
For $0 \leq r \leq 2n$, one has 
    \[
    \zeta_{n}(r) = \frac{1}{r^{4}} 
     \int_{0}^{r} s^{3} \eta_{n}(s) \, ds
     \leq \frac{4}{r^{4}} 
     \int_{0}^{r} s^{3} \, ds = 1. 
    \]
Finally, 
for $2n \leq r$, we obtain 
    \[
    \zeta_{n}(r) = \frac{1}{r^{4}} 
     \int_{0}^{r} s^{3} \eta_{n}(s) \, ds
     = \frac{1}{r^{4}} 
     \int_{0}^{2n} s^{3} \eta_{n}(s) \, ds
     + \frac{1}{r^{4}} 
     \int_{2n}^{r} s^{3} \eta_{n}(s) \, ds
     \leq \frac{4}{r^{4}} \int_{0}^{2n} s^{3}  \, ds
     = \frac{16 n^{4}}{r^{4}}. 
    \]
This completes the proof. 
\end{proof}
We set
\[
I_{n}\left(u(t), v(t), \p_{t} v(t)\right) 
:= \text{Im} \int_{\R^{4}} (x \cdot \nabla u) \overline{u} \zeta_{n} \, dx 
- 2 \int_{\R^{4}} \p_{t} v (x \cdot \nabla v) \zeta_{n} \, dx 
- \int_{\R^{4}} \p_{t} v v \eta_{n} \, dx.   
\]

\begin{lemma}\label{lem-lvirial}
The following identity holds:
\begin{equation}\label{eq0:virial}
\frac{d}{dt} I_{n} \left(u(t), v(t), \p_{t} v(t)\right) 
= 2\mathcal{K}(u(t), v(t)) + \mathcal{N}_{n}\left(u(t), v(t), \p_{t} v(t)\right), 
\end{equation}
where
\begin{align*}
\mathcal{N}_n\left(u, v, \p_{t} v\right)
&:= \int_{\R^{4}} |\nabla u|^2 (2 \zeta_{n} - 2) \, dx 
+ \int_{\R^{4}} |\nabla v|^2 (\eta_n - 2 \zeta_{n} - x \cdot \nabla \zeta_{n} -2) \, dx \\
& \quad + \int_{\R^{4}} |v|^2 (-4\zeta_{n} + \eta_n - x \cdot \nabla \zeta_{n}) \, dx 
+ \int_{\R^{4}} |u|^2 v (4 - \eta_n) \, dx \\
& \quad + \int_{\R^{4}} \left|\p_{t} v\right|^2 (x \cdot \nabla \zeta_{n} + 4 \zeta_{n} - \eta_n) \, dx 
+ 2 \mathrm{Re} \int_{\R^{4}} (\nabla u \cdot \nabla\zeta_{n}) (x \cdot \nabla\overline{u}) \, dx \\
& \quad - 3 \int_{\R^{4}} |u|^{2} \Delta \zeta_{n} \, dx 
+ 2 \int_{\R^{4}} (x \cdot \nabla v) (\nabla v \cdot \nabla\zeta_{n}) \, dx 
- \frac{1}{2} \int_{\R^{4}} |v|^2 \Delta \eta_n \, dx 
- \frac{1}{2} \int_{\R^{4}} [x \cdot \nabla (\Delta \zeta_{n})] |u|^{2} \, dx.
\end{align*}
\end{lemma}
\begin{proof}
By a direct computation with a little bit effort, 
we obtain 
    \begin{equation} \label{eq1-virial}
\begin{split}
\frac{d}{dt} \text{Im} \int_{\R^{4}} (x \cdot \nabla u) \overline{u} \zeta_{n} \, dx
& = 2 \int_{\R^{4}} |\nabla u|^2 \zeta_{n} \, dx 
+ 2 \int_{\R^{4}} |u|^2 (x \cdot \nabla v) \zeta_{n} \, dx   \\
& \quad 
+ 2 \mathrm{Re} \int_{\R^{4}} (\nabla u \cdot \nabla\zeta_{n}) (x \cdot \nabla\overline{u}) \, dx  
- 3 \int_{\R^{4}} |u|^{2} \Delta \zeta_{n} \, dx \\
& \quad  
- \frac{1}{2} \int_{\R^{4}} [x \cdot \nabla (\Delta \zeta_{n})] |u|^{2} \, dx, 
\end{split}
\end{equation}
\begin{equation} \label{v-eq17}
\begin{split}
\frac{d}{dt} \int_{\R^{4}} \p_{t} v (x \cdot \nabla v) \zeta_{n} \, dx
&= \int_{\R^{4}} |\nabla v|^2 \zeta_{n} \, dx 
+ 2 \int_{\R^{4}} |v|^2 \zeta_{n} \, dx 
+ \int_{\R^{4}} (x \cdot \nabla v) |u|^2 \zeta_{n} \, dx  \\
& \quad - 2 \int_{\R^{4}} 
\left|\p_{t} v\right|^2 \zeta_{n} \, dx	
- \frac{1}{2} \int_{\R^{4}} \left|\p_{t} v\right|^2 (x \cdot \nabla\zeta_{n}) \, dx   \\
& \quad + \frac{1}{2} \int_{\R^{4}} (x \cdot \nabla\zeta_{n}) |\nabla v|^2 \, dx
- \int_{\R^{4}} (x \cdot \nabla v) (\nabla v \cdot \nabla\zeta_{n}) \, dx  
+ \frac{1}{2} \int_{\R^{4}} |v|^2 (x \cdot \nabla\zeta_{n}) \, dx,  
\end{split}
\end{equation}
\begin{equation} \label{v-eq18-2}
\begin{split}
	\frac{d}{dt} \int_{\R^{4}} v 
    \p_{t} v \eta_n \, dx 
	= - \int_{\R^{4}} |\nabla v|^2 \eta_n \, dx
    + \frac{1}{2} \int_{\R^{4}} |v|^2 \Delta \eta_n \, dx 
    - \int_{\R^{4}} |v|^2 \eta_n \, dx 
	+ \int_{\R^{4}} |u|^2 v \eta_n \, dx 
    + \int_{\R^{4}} |\p_{t} v|^2 \eta_n \, dx.
\end{split}
\end{equation}
Then, \eqref{eq0:virial} follows from 
\eqref{eq1-virial}--\eqref{v-eq18-2}. 
\end{proof}

\begin{lemma}
Let $(u(t), v(t))$ be a radially symmetric solution to 
\eqref{kgs} with initial data in $X_{\text{rad}}$, where 
$X_{\text{rad}} = \left\{u \in X 
\colon u(x) = u(|x|) \right\}$. 
Then, there exists a positive constant $C_{0}$ which is 
independent of $n \in \mathbb{N}$ such that 
    \begin{equation} \label{v-eq30}
    \begin{split}
    \frac{d}{d t} I_{n} (u(t), v(t), \p_{t} v(t))
\leq 2 \mathcal{K}(u(t), v(t)) 
    + 4 \int_{|x| \geq n} |u|^{2} |v|(t, x) \, dx 
    + \frac{C_{0}^2}{n^{2}} 
    (\|u(t)\|_{L^{2}}^{2} + \|v(t)\|_{L^{2}}^{2}). 
    \end{split}
    \end{equation}
\end{lemma}
\begin{proof}
Using Lemma \ref{lem-cutoff}, we can estimate 
several terms in $\mathcal{N}_n$ as follows:
\begin{align*}
\int_{\R^{4}} |\nabla u|^2 (2 \zeta_{n} - 2) \, dx &\leq 0, \\
\int_{\R^{4}} |\nabla v|^2 (\eta_n - 2 \zeta_{n} - x \cdot \nabla \zeta_{n} - 2) \, dx &
= \int_{\R^{4}} |\nabla v|^2 (2 \zeta_{n} - 2) \, dx
\leq 0, \\
\int_{\R^{4}} |v|^2 (-4\zeta_{n} + \eta_n - x \cdot \nabla \zeta_{n}) \, dx &= 0, \\
\int_{\R^{4}} |u|^2 v (4 - \eta_n) \, dx &\leq  4 \int_{|x| \geq n} |u|^2 |v| \, dx, \\
\int_{\R^{4}} \left|\p_{t} v\right|^2 (x \cdot \nabla \zeta_{n} + 4 \zeta_{n} - \eta_n) \, dx &= 0.
\end{align*}
Furthermore, by radial symmetry and Lemma \ref{lem-cutoff} \textrm{(iv)},
\begin{align*}
2 \mathrm{Re} \int_{\R^{4}} (\nabla u \cdot \nabla\zeta_{n}) (x \cdot \nabla\overline{u}) \, dx 
= 2 \int_{0}^{\infty} r \zeta_{n}^{\prime} |u_{r}|^2 r^3 dr 
&\leq 0, \\
2 \int_{\R^{4}} (x \cdot \nabla v) (\nabla v \cdot \nabla\zeta_{n}) \, dx  
= 2 \int_{\R^{4}} (x \cdot \nabla \zeta_{n}) 
|\nabla v|^2 \, dx &\leq 0. 
\end{align*}
The remaining error terms satisfy
\[
\begin{split}
\left| - 3 \int_{\R^{4}} |u|^{2} \Delta \zeta_{n} \, dx \right| 
+ \left| -\frac{1}{2} \int_{\R^{4}} |v|^2 \Delta \eta_n \, dx \right| 
+ \left| -\frac{1}{2} \int_{\R^{4}} [x \cdot \nabla (\Delta \zeta_{n})] |u|^{2} \, dx \right| 
\leq \frac{C}{n^2} 
\left( \|u\|_{L^2}^2 + \|v\|_{L^2}^2 \right).
\end{split}
\]    
Thus, we see that \eqref{v-eq30} holds. 
\end{proof}
    
We are now in a position to prove Theorem \ref{thm-grow-blow}. 
\begin{proof}[Proof of Theorem \ref{thm-grow-blow}]
We shall show Theorem \ref{thm-grow-blow} 
by contradiction. 
Suppose to the contrary that there exists 
$(u_{0}, v_{0}, v_{1}) \in \mathcal{R}_{-}(m)$ 
such that 
$(u,v,\partial_t v)$ is uniformly bounded in $X$, where 
$(u, v, \partial_{t} v)$ is a solution to \eqref{kgs} with 
$(u, v, \partial_{t} v)|_{t = 0} = (u_{0}, v_{0}, v_{1})$. 
Then, there exists $C>0$ such that 
    \begin{equation} \label{eq2:virial}
    \sup_{t \in [0,\infty)} 
    \|(u(t), v(t), \partial_t v(t))\|_{X}  
    \leq C. 
    \end{equation}
Since $(u, v, \partial_{t} v) \in \mathcal{R}_{-}$, 
it follows from Lemma \ref{lem-uni-neg} that 
there exists $\delta > 0$ such that 
$\sup_{t \in [0,\infty)} \mathcal{K}(u(t), v(t)) 
\leq - \delta$. 
Therefore, by \eqref{v-eq30}, H\"{o}lder inequality, and \cite[Lemma 1]{OgTs91},
we can take sufficiently large 
$N \in \mathbb{N}$ such that 
\begin{equation*}
\begin{split}
\frac{d}{d t} I_{n}\left(u(t), v(t), \p_{t} v(t)\right) 
 & \leq 2\mathcal{K}(u(t), v(t))
+ 4 \int_{|x| \geq n} |u|^{2} |v|(t, x) \, dx 
+ \frac{C_{0}^2}{n^{2}} 
    (\|u(t)\|_{L^{2}}^{2} + \|v(t)\|_{L^{2}}^{2}) \\
 & \leq 2\mathcal{K}(u(t), v(t))
+ 4 \|u(t)\|_{L^3(|x|\geq n)}^2\|v(t)\|_{L^3(|x|\geq n)}
+ \frac{C_{0}^2}{n^{2}} 
    (\|u(t)\|_{L^{2}}^{2} + \|v(t)\|_{L^{2}}^{2}) \\
 & \leq 2\mathcal{K}(u(t), v(t))
+ \frac{C}{n^\frac{3}{2}}\|u(t)\|_{H^1}^2\|v(t)\|_{H^1}
+ \frac{C_{0}^2}{n^{2}} 
    (\|u(t)\|_{L^{2}}^{2} + \|v(t)\|_{L^{2}}^{2}) \\
 & \leq - \delta + \frac{\delta}{2} = 
- \frac{\delta}{2}  
\end{split}
\end{equation*}
for all $n \geq N$. 
Integrating the above on the interval $[0, t]$, we obtain 
    \[
    \begin{split}
   I_{n}\left(u(t), v(t), \p_{t} v(t)\right)  - I_{n}\left(u_0, v_0, v_1\right)  \leq - \frac{\delta}{2} t. 
    \end{split}
    \]
Therefore, we obtain
\begin{equation} \label{eq3:virial}
\begin{split}
\lim_{t\to\infty} |I_{n}\left(u(t), v(t), \p_{t} v(t)\right) | = \infty. 
\end{split}
\end{equation}
On the other hand,
by H\"{o}lder's inequality,  
\eqref{eq2:virial} and 
Lemma \ref{lem-cutoff} \textrm{(i)} and 
\textrm{(v)}, we have
\begin{equation*}
\begin{split}
\left|I_{n}\left(u(t), v(t), \p_{t} v(t)\right) \right| 
\leq C n^4\left(\|u(t)\|_{H^1}^2 + \|v(t)\|_{H^1}^2 + 
\left\|\p_{t} v(t)
\right\|_{L^2}^2 \right) 
= C n^4\|(u(t), v(t), \p_{t} v(t))\|_{X}^2 
\leq C, 
\end{split}
\end{equation*}
which contradicts \eqref{eq3:virial}. 
This completes the proof.  
\end{proof}

\appendix 

\section{Cauchy problem}   
\label{sec:cauchy}
In this appendix, we shall give the proof of 
Theorem \ref{thm1:local-wp}. 
\begin{definition}\label{def1:ad}
We say that a pair $(q, r)$ is 
\textit{$L^{2}$-admissible} if 
$(q, r)$ satisfies 
\begin{align*}
   \frac{2}{q}=d\left(\frac{1}{2}-\frac{1}
   {r}\right),\quad 2\leq q \leq \infty,\quad 
   2 \leq r \leq 
   \frac{2d}{d - 2}.  
\end{align*}    
\end{definition}

We now recall the Strichartz estimate.
\begin{theorem}[Strichartz estimate]
\label{thm:ST1}
The following properties hold:
\begin{enumerate}
\item[\textrm{(i)}] For every $\phi \in L^2(\R^d)$, 
the function $t \to e^{i t \Delta} \phi$ belongs to
\begin{equation*}
L^q (\R,L^r (\R^d)) \cap C(\R, L^2(\R^d))
\end{equation*}
for every admissible pair $(q,r)$. 
Furthermore, there exists a positive constant $C_{\text{St}}$ such that
\begin{equation} \label{eq1-lw}
\|e^{i t \Delta}\phi\|_{L_t^q L_x^r (\R \times \R^d)} \leq C_{\text{St}}\|\phi\|_{L^2} 
\quad \text{for every} \; \phi \in L^2(\R^d).
\end{equation}
\item[\textrm{(ii)}] Let $I$ be an interval of $\R$ (bounded or not), $J = \overline{I}$, and $t_0 \in J$. If $(\gamma, \rho)$ is an admissible pair and $f \in L^{\gamma'}(I,L^{\rho'}(\R^d))$, then for every admissible pair $(q,r)$, the function
\begin{equation*}
t \to \Phi_{f}(t) = \int_{t_0}^{t} e^{i (t-s) \Delta} f(s)\, ds \quad \text{for} \; t \in I
\end{equation*}
belongs to $L^q (I,L^r (\R^d)) \cap C(J, L^2(\R^d))$. Furthermore, there exists a 
positive constant $C$ independent of $I$ such that
\begin{equation}\label{eq2-lw}
\|\Phi_{f}\|_{L_t^q L_x^r (I \times \R^d)} \leq C\|f\|_{L_t^{\gamma'}L_x^{\rho'} (I \times \R^d)} \quad \text{for every} \; f \in L^{\gamma'}L^{\rho'} (I \times \R^d). 
\end{equation}
\end{enumerate}
\end{theorem}
We now give the proof of Theorem \ref{thm1:local-wp}. 
\begin{proof}[Proof of Theorem \ref{thm1:local-wp}]
The Cauchy problem of the system \eqref{kgs} 
with $(u(t), v(t), \partial_{t} v(t))|_{t= 0} = 
(u_0, v_0, v_1)$ is equivalent to the following integral system: 
\begin{equation*}
\begin{cases}
   u(t) = e^{it\Delta}u_0 - 2i \int_0^t e^{i(t-s)\Delta} u(s)v(s) \, ds,  & \\
   v(t) = \cos(t\langle\nabla\rangle)v_0 + \langle\nabla\rangle^{-1}\sin(t\langle\nabla\rangle)v_1 + \int_0^t \langle\nabla\rangle^{-1}\sin((t-s)\langle\nabla\rangle)|u(s)|^2 \, ds, &
\end{cases}
\end{equation*} 
where $\langle\nabla\rangle=(1-\Delta)^{1/2}$.

We set $M := \max\{
    \|u_0\|_{H^{1}}, \|v_{0}\|_{H^{1}}, 
    \|v_{1}\|_{L^{2}}
    \}$
and define a function space $E$ by 
    \begin{equation}\label{eq2:function}
    \begin{split}
    E = 
    \Big\{ 
    & (u, v) \in L^{\infty}([0, T), H^{1}(\R^{4})) 
    \times L^{\infty}([0, T), H^{1}(\R^{4})) 
    \colon \\
    & \|u\|_{L^{\infty}
    ([0, T), H^{1}(\R^4))} \leq 3 M, \; 
    \|v\|_{L^{\infty}
    ([0, T), H^{1}(\R^4))} \leq 3 M
    \Big\},  
    \end{split}
    \end{equation}
which is equipped with the norm 
    \[
    \|(u, v)\|_{E} = 
    \|u\|_{L_{t}^{\infty}
    ([0, T), H_{x}^{1}(\R^4))} 
    + \|v\|_{L_{t}^{\infty}
    ([0, T), H_{x}^{1}(\R^4))}. 
    \]
We set 
    \begin{align*}
   & \mathcal{T}[u, v] 
   := e^{it\Delta}u_0 - 2i \int_0^t e^{i(t-s)\Delta} u(s)v(s) \, ds,   \\
   &\mathcal{S}[u, v]  
   := \cos(t\langle\nabla\rangle)v_0 + \langle\nabla\rangle^{-1}\sin(t\langle\nabla\rangle)v_1 + \int_0^t \langle\nabla\rangle^{-1}\sin((t-s)\langle\nabla\rangle)|u(s)|^2 \, ds. 
\end{align*}
We shall show that the map $(\mathcal{T}, 
\mathcal{S}): E \to E$ is a contraction. 
First, we claim that $\|\mathcal{S}[u, v]\|_{H^{1}} 
\leq 3 M$ for all $t \in [0, T)$ if $T > 0$ is 
sufficiently small. 
Observe that 
    \begin{equation}
    \begin{split}
    \|\mathcal{S}[u, v]\|_{H^{1}} 
    & = \left\| 
    \cos(t\langle\nabla\rangle)v_0 + \langle\nabla\rangle^{-1}\sin(t\langle\nabla\rangle)v_1 + \int_0^t \langle\nabla\rangle^{-1}\sin((t-s)\langle\nabla\rangle)|u(s)|^2 \, ds
    \right\|_{H^{1}} \\
    & \leq \|v_{0}\|_{H^{1}} 
    + \|v_{1}\|_{L^{2}} 
    + \|u\|_{L_{t}^{2}L_{x}^{4}}^2.
    \end{split}    
    \end{equation}
We have 
by the Sobolev inequality that 
 \begin{equation} \label{eq1:cauchy-1}
    \begin{split}
    \|u\|_{L_{t}^{2}L_{x}^{4}}^2 
    = \int_{0}^{T} \|u\|_{L_{x}^4}^2 \, ds 
    \leq C 
    \int_{0}^{T} \|\nabla u\|_{L_{x}^2}^{2} \, ds 
    \leq C \|u\|_{L_{t}^{\infty} 
   H_{x}^1}^{2}\int_{0}^{T} \, ds
   \leq C T \|u\|_{L_{t}^{\infty} 
   H_{x}^1}^{2}. 
    \end{split}
    \end{equation}
This implies that 
$\|\mathcal{S}[u, v]\|_{L_{t}^{\infty}H_{x}^{1}} 
\leq 3 M$ if $T > 0$ is 
sufficiently small. 
    
Next, we shall show that $\|\mathcal{T}[u, v]\|_{L_{t}^{\infty}H_{x}^{1}} 
\leq 3 M$ if $T > 0$ is 
sufficiently small. 
Using the Strichartz estimate 
\eqref{eq2-lw} with 
$(q, r) = (\infty, 2)$ and $(\gamma, \rho) = 
(2, 4)$,  
we obtain 
    \begin{equation} \label{eq2:cauchy}
    \begin{split}
    \|\mathcal{T}[u, v]\|_{L_{t}^{\infty}H_{x}^{1}} 
    & = \left\| 
    e^{it\Delta}u_0 - 2i \int_0^t e^{i(t-s)\Delta} u(s)v(s) \, ds
    \right\|_{L_{t}^{\infty}H_{x}^{1}} \\
   &\leq \|u_0\|_{H^1} 
   + 2 \left\|\int_0^t e^{i(t-s)\Delta} 
   \langle \nabla \rangle (u(s)v(s)) \, ds
   \right\|_{L_{t}^{\infty}(I:L_{x}^2)}   \\
   &\leq \|u_0\|_{H_{x}^1} 
   + C \Big(\|v \langle \nabla \rangle u\|_{L_{t}^{2}(I:L_{x}^{\frac{4}{3}})} 
   + \|u \langle \nabla \rangle v\|_{L_{t}^{2}(I:L_{x}^{\frac{4}{3}})} \Big). 
    \end{split}    
    \end{equation}
Then, using the H\"{o}lder and Sobolev inequalities,  
we can estimate the second and third 
terms of \eqref{eq2:cauchy} as follows: 
\begin{equation}\label{eq4:cauchy}
\begin{split}
& \quad \|v \langle \nabla \rangle u
\|_{L_{t}^{2}L_{x}^{\frac{4}{3}}}
=\| \|v \langle \nabla \rangle u
\|_{L_{x}^{\frac{4}{3}}} \|_{L_{t}^2} 
\leq \| \|\langle \nabla \rangle u\|_{L_{x}^2} 
\|v\|_{L_{x}^4} \|_{L_{t}^2} 
\leq C \| \|\langle \nabla \rangle u\|_{L_{x}^2} 
\|v\|_{H_{x}^1} \|_{L_{t}^2}   \\
&\leq C \|v\|_{L_{t}^{\infty}H_{x}^1} 
\|\langle \nabla \rangle u\|_{L_{t}^{\infty}L_{x}^2} 
\| 1 \|_{L_{t}^2} 
\leq C T^{\frac{1}{2}} \|v\|_{L_{t}^{\infty}
H_{x}^1} 
\|u\|_{L_{t}^{\infty}H_{x}^1}.    
\end{split}
\end{equation}
Interchanging the roles of $u$ and $v$ 
in \eqref{eq4:cauchy} yields
\begin{equation}\label{eq5:cauchy}
\|u\langle \nabla \rangle v
\|_{L_{t}^{2}L_{x}^{\frac{4}{3}}}
\leq C T^{\frac{1}{2}} 
\|u\|_{L_{t}^{\infty}H_{x}^1} 
\|v\|_{L_{t}^{\infty}H_{x}^1}.  
\end{equation}
Thus, we see from \eqref{eq2:cauchy}--\eqref{eq5:cauchy} that $\|\mathcal{T}[u, v]\|_{L_{t}^{\infty}H_{x}^{1}} 
\leq 3 M$ if $T > 0$ is 
sufficiently small. 


Next, we shall show that $(\mathcal{T}, 
\mathcal{S})$ is a contraction. 
By H\"{o}lder's inequality, one has
\begin{equation}\label{eq11:cauchy}
\begin{split}
\|\mathcal{S}[u_1, v_1]-\mathcal{S}[u_2, v_2]\|_{L_{t}^{\infty} H_{x}^{1}}  
&= \left\|\int_0^t \langle\nabla\rangle^{-1}\sin((t-s)\langle\nabla\rangle)(|u_1(s)|^2-|u_2(s)|^2) \, ds \right\|_{L_{t}^{\infty} H_{x}^{1}}   \\
&\leq \int_0^T \||u_1|^2-|u_2|^2\|_{L_{x}^2} 
\, ds  \\
&\leq \int_0^T (\|u_1\|_{L_{x}^4} 
+\|u_2\|_{L_{x}^4}) \|u_1-u_2\|_{L_{x}^4} \, ds  \\
&\leq T (\|u_1\|_{L_{t}^{\infty}H_{x}^1} 
+\|u_2\|_{L_{t}^{\infty}H_{x}^1}) 
\|u_1-u_2\|_{L_{t}^{\infty}H_{x}^1} \\
& \leq C_{M} T \|u_1-u_2\|_{E}
\end{split}
\end{equation}
Next, it follows from the Strichartz estimate 
\eqref{eq2-lw} with 
$(q, r) = (\infty, 2)$ and $
(\gamma, \rho) = (2, 4)$ that 
\begin{equation} \label{eq7:cauchy}
    \begin{split}
    \|\mathcal{T}[u_1, v_1]-\mathcal{T}[u_2, v_2]\|_{L_{t}^{\infty}H_{x}^1}  
    &\leq 
    2 \left\|\int_0^t e^{i(t-s)\Delta} 
    \langle \nabla \rangle (u_1(s)v_1(s)-u_2(s)v_2(s)) \, ds\right\|_{L_{t}^{\infty}L_{x}^2}   \\
   &\leq 
   C\|\langle \nabla \rangle (u_1v_1-u_2v_2)\|_{L^{2}(I:L^{\frac{4}{3}})}   \\
   &\leq
   C (\|(\langle \nabla \rangle u_1)v_1 - 
   (\langle \nabla \rangle u_2)v_2 \|_{L_{t}^{2}(I:L_{x}^{\frac{4}{3}})} 
   + \| u_1(\langle \nabla \rangle v_1) - 
   u_2(\langle \nabla \rangle v_2)\|_{L_{t}^{2}(I:L_{x}^{\frac{4}{3}})}).   
    \end{split}    
\end{equation}
We consider the first term of 
\eqref{eq7:cauchy}. 
As in \eqref{eq4:cauchy}, 
using the H\"{o}lder and Sobolev inequalities, 
one has 
\begin{equation} \label{eq9:cauchy}
\begin{split}
\|(\langle \nabla \rangle u_1)v_1 - 
(\langle \nabla \rangle u_2)v_2\|_{L_{t}^{2}
(I:L_{x}^{\frac{4}{3}})}   
&\leq \|\langle \nabla \rangle u_1(v_1-v_2)\|_{L_{t}^{2}(I:L_{x}^{\frac{4}{3}})} 
+ \|(\langle \nabla \rangle u_1- 
\langle \nabla \rangle u_2)
v_2\|_{L_{t}^{2}(I:L^{\frac{4}{3}})}   \\
&\leq C \|u_1\|_{L_{t}^{\infty}
H_{x}^1} \|v_1-v_2\|_{L_{t}^{\infty}H_{x}^1} \|1\|_{L_{t}^2(I)} 
+ C \|u_1 - u_2\|_{L_{t}^{\infty}(I)H_{x}^1} 
\|v_2\|_{L_{t}^{\infty}H_{x}^1} 
\|1\|_{L_{t}^2}   \\
&\leq C T^{\frac{1}{2}} \|u_1\|_{L_{t}^{\infty}H_{x}^1} 
\|v_1-v_2\|_{L^{\infty}H^1}  
+ C T^{\frac{1}{2}} 
\|v_2\|_{L_{t}^{\infty}H_{x}^1} 
\|u_1 - u_2\|_{L_{t}^{\infty}H_{x}^1} \\
& \leq C T^{\frac{1}{2}} 
\|(u_{1} - u_{2}, v_{1} - v_{2})\|_{E}.  
\end{split}
\end{equation}
By swapping $u_i$ and $v_i$ ($i = 1, 2$) in \eqref{eq9:cauchy}, we have
\begin{equation} \label{eq10:cauchy}
\|u_1(\langle \nabla \rangle v_1)- 
u_2(\langle \nabla \rangle v_2)
\|_{L_{t}^{2}(I:L_{x}^{\frac{4}{3}})}  
\leq C T^{\frac{1}{2}} 
\|(u_{1} - u_{2}, v_{1} - v_{2})\|_{E}.
\end{equation}
It follows from \eqref{eq7:cauchy} and 
\eqref{eq10:cauchy} that 
    \begin{equation}\label{eq12:cauchy}
    \|\mathcal{T}[u_1, v_1]-\mathcal{T}[u_2, v_2]\|_{L_{t}^{\infty}H_{x}^1}
    \leq C T^{\frac{1}{2}} 
\|(u_{1} - u_{2}, v_{1} - v_{2})\|_{E}. 
    \end{equation}
From \eqref{eq11:cauchy} and \eqref{eq12:cauchy}, 
we see that 
    \[
    \|
    (\mathcal{T}[u_1, v_1], \mathcal{S}[u_1, v_1]) - (\mathcal{T}[u_2, v_2], 
    \mathcal{S}[u_2, v_2])\|_{E} 
    \leq \frac{1}{2} \|(u_{1} - u_{2}, v_{1} - v_{2})\|_{E}
    \]
for sufficiently small $T>0$, 
so that $(\mathcal{T}, 
\mathcal{S})$ is a contraction mapping. 
This completes the proof. 
\end{proof}

\section{Regularity of solution to \eqref{sp} and Pohozaev Identity}
\label{sec:technical}
In this appendix, we consider the 
following semilinear elliptic system:
\begin{equation} \label{asp}
	\begin{cases}
	- \Delta \varphi + \ell \varphi = 2 \varphi \psi 
	& \qquad \mbox{in $\R^{d}$}, \\
	- \Delta \psi + k \psi = |\varphi|^{2} 
	& \qquad \mbox{in $\R^{d}$},  
	\end{cases}
	\end{equation}
where $\ell \geq 0$ and $k \geq 0$. 
We put 
    \[
    H_{k}(\R^{4}) := 
    \begin{cases}
   H^{1}(\R^{4})
    & \quad \mbox{when $k > 0$}, \\
     \dot{H}^{1}(\R^{4})
     & \quad \mbox{when $k = 0$},  
    \end{cases}
    \qquad 
    H_{\ell, k} := H_{\ell}(\R^{4}) 
    \times H_{k}(\R^{4}). 
    \]
\begin{lemma}\label{lem1-regu}
Assume that $(\Phi, \Psi) 
\in H_{\ell, k}$ is a weak solution 
to \eqref{asp}. 
Then, $\Phi, \Psi \in 
C^{2, \gamma}(\R^{4})$ 
for all $0 < \gamma < 1$ and 
$|\Phi(x)|, |\Psi(x)| 
\to 0$ as $|x| \to \infty$. 
\end{lemma}
\begin{proof}
By the embedding $\dot{H}^{1}(\R^{4}) 
\hookrightarrow L^{4}(\R^{4})$ 
that $|\Phi|^2 \in L^{2}(\R^{4})$. 
This together with  the second equation of 
\eqref{asp} 
implies that $\Psi\in \dot{H}^{2, 2}
(\R^{4})$. 
Then, we have by the Gagliardo-Nirenberg inequality 
that $\Psi\in L^{q_{1}}(\R^{4})$ for $4 
\leq q_{1} < \infty$. 
By the first equation of \eqref{asp}
and $\Phi \Psi\in L^{r_{1}}(\R^{4})$ for 
$2 \leq r_{1} < 4$, 
we have  $\Phi \in 
\dot{W}^{2, r_{1}}(\R^{4})$. 
This together with the 
Gagliardo-Nirenberg inequality and 
$\Phi \in L^{4}(\R^{4})$
implies that 
$\Phi\in L^{q_{1}}(\R^{4})$ for 
$4 \leq q_{1} < \infty$. 
Thus, we obtain $\Phi \Psi\in 
L^{r_{2}}(\R^{4})$ for $2 
\leq r_{2} < \infty$. 
This together with the first equation of \eqref{asp} yields that 
$\Phi \in \dot{W}^{2, r_{2}}(\R^{4})$, so that 
$\Phi\in C^{1, \gamma}(\R^{4})$ by the 
Morrey inequality. 
It follows that $\Phi
\p_{x_{j}}\Phi 
\in L^{r_{2}}(\R^{4})$ for $1 \leq j \leq 4$. 
Then, from the second equation of \eqref{asp}, we obtain  
$\p_{x_{j}} \Psi \in \dot{W}^{2, r_{2}}(\R^{4})$, so that $\Psi \in \dot{W}^{3, r_{2}}(\R^{4})$. 
This together with the Morrey inequality 
implies that 
$\Psi\in C^{2, \gamma}(\R^{4})$ for all $0 < \gamma < 1$. 
Then, we have $\p_{x_{j}}(\Phi \Psi) 
\in L^{r_{2}}(\R^{4})$. Then,
 from the first equation of 
\eqref{asp} it follows that $\p_{x_{j}} \Phi 
\in \dot{W}^{2, r_{2}}(\R^{4})$. 
Hence, $\Phi \in 
\dot{W}^{3, r_{2}}(\R^{4})$. 
By the Sobolev 
embedding, we obtain
$\Phi \in C^{2, \gamma}(\R^{4})$ 
for all $0 < \gamma < 1$. 
In addition, since $\Phi \in 
C^{2, \gamma} \cap L^{4}(\R^{4})$
and $\Psi\in 
C^{2, \gamma} \cap L^{4}(\R^{4})$, 
we have $|\Phi(x)|, |\Psi(x)| 
\to 0$ as $|x| \to \infty$. 
This completes the proof. 
\end{proof}
\begin{lemma}\label{lem2-pohozaev}
Let $(\Phi, \Psi) \in H_{\ell, k}$ be a weak solution to \eqref{asp}. 
Then, the following identities hold: 
\begin{align}
    & \|\nabla \Phi\|_{L^{2}}^{2} 
    + 2 \ell \|\Phi\|_{L^{2}}^{2}
    - 4 \int_{\R^{4}} |\Phi|^{2} \Psi\, dx 
    - \int_{\R^{4}} |\Phi|^{2} (x \cdot \nabla \Psi)\, dx
    = 0, \label{eq2-poho} \\ 
    & \|\nabla \Psi\|_{L^{2}}^{2} 
    + 2 k \|\Psi\|_{L^{2}}^2
    + \int_{\R^{4}} |\Phi|^{2} (x \cdot \nabla \Psi) \, dx 
    = 0. \label{eq3-poho}
\end{align}
\end{lemma}

\begin{proof}
We follow the argument of Berestycki and P.-L. ~Lions~\cite[Proposition 1]{BeLi83}. 
By Lemma \ref{lem1-regu}, $(\Phi, \Psi)$ is smooth enough to justify the following calculations.
Let $R > 0$ and $B_{R} = \{x \in \R^{4} \colon |x| < R \}$. 
Multiplying the first equation of \eqref{asp} by $x \cdot \nabla \Phi = \sum_{i=1}^{4} x_{i} \p_{x_{i}} \Phi$ and integrating over $B_{R}$, we have
\[
   - \sum_{i, j = 1}^{4} \int_{B_{R}} \p_{x_{j} x_{j}} \Phi x_{i} \p_{x_{i}} \Phi \, dx
   + \ell \sum_{i = 1}^{4} \int_{B_{R}} \Phi x_{i} \p_{x_{i}} \Phi \, dx
   - 2 \sum_{i = 1}^{4} \int_{B_{R}} \Phi \Psi x_{i} \p_{x_{i}} \Phi \, dx = 0. 
\]
Integrating by parts yields
\[
\begin{split}
\sum_{i = 1}^{4} \int_{B_{R}} \Phi \Psi x_{i} \p_{x_{i}} \Phi \, dx 
&= -2 \int_{B_{R}} \Phi^2 \Psi\, dx 
- \frac{1}{2} \int_{B_{R}} \Phi^{2} (x \cdot \nabla \Psi) \, dx 
+ \frac{1}{2} \int_{\p B_{R}} \Phi^2 \Psi(x \cdot \nu) \, dS, \\
\sum_{i = 1}^{4} \int_{B_{R}} \Phi x_{i} \p_{x_{i}} \Phi \, dx 
&= - 2 \int_{B_{R}} \Phi^2 \, dx 
+ \frac{1}{2} \int_{\p B_{R}} \Phi^2 (x \cdot \nu) \, dS,
\end{split}
\]
where $\nu = (\nu_1, \dots, \nu_4) = x/|x|$ denotes the outer unit normal vector on $\p B_R$.
Similarly, integrating by parts again, we obtain
Integrating by parts again, we obtain
\[
\begin{split}
- \sum_{i, j = 1}^{4} \int_{B_{R}} \p_{x_{j} x_{j}} \Phi x_{i} \p_{x_{i}} \Phi \, dx 
&= \sum_{i, j = 1}^{4} \int_{B_{R}} \p_{x_{j}} \Phi \left( \delta_{i j} \p_{x_{i}} \Phi + x_{i} \p_{x_{i} x_{j}} \Phi \right) \, dx 
- \int_{\p B_{R}} \left( \sum_{j=1}^{4} \p_{x_{j}} \Phi \nu_{j} \right) \left( \sum_{i=1}^{4} x_{i} \p_{x_{i}} \Phi \right) \, dS \\
&= \int_{B_{R}} |\nabla \Phi|^{2} \, dx 
+ \sum_{i = 1}^{4} \int_{B_{R}} x_{i} \p_{x_{i}} \left( \frac{|\nabla \Phi|^2}{2} \right) \, dx 
- \int_{\p B_{R}} \left| \frac{\p \Phi}{\p \nu} \right|^2 R \, dS \\
&= \int_{B_{R}} |\nabla \Phi|^{2} \, dx 
- 2 \int_{B_{R}} |\nabla \Phi|^{2} \, dx 
+ \frac{1}{2} \int_{\p B_{R}} |\nabla \Phi|^2 R \, dS 
- \int_{\p B_{R}} \left| \frac{\p \Phi}{\p \nu} \right|^2 R \, dS \\
&= - \int_{B_{R}} |\nabla \Phi|^{2} \, dx 
+ \int_{\p B_{R}} \left( \frac{1}{2} |\nabla \Phi|^{2} - \left| \frac{\p \Phi}{\p \nu} \right|^2 \right) R \, dS. 
\end{split}
\]
Combining these identities gives
\[
\begin{split}
& \int_{B_{R}} |\nabla \Phi|^{2} \, dx 
+ 2 \omega \int_{B_{R}} \Phi^2 \, dx 
- 4 \int_{B_{R}} \Phi^2 \Psi\, dx 
- \int_{B_{R}} \Phi^{2} (x \cdot \nabla \Psi) \, dx \\
&= - \int_{\p B_{R}} \left( \left|\frac{\p \Phi}{\p \nu}\right|^{2} - \frac{1}{2} |\nabla \Phi|^{2} \right) R \, dS
+ \frac{\ell}{2} \int_{\p B_{R}} \Phi^2 R \, dS 
- \int_{\p B_{R}} \Phi^2 \Psi R \, dS,
\end{split}
\]
where we used $x \cdot \nu = R$ on $\p B_R$.
Since $(\Phi, \Psi) \in H_{\ell, k}$, 
it holds that
\begin{equation} \label{eq1-poho}
\int_{\R^{4}} \left\{ |\nabla \Phi|^{2} + 
\ell \Phi^{2} + \Phi^{2} |\Psi| \right\} \, dx < \infty. 
\end{equation}
Then, we claim that there exists a sequence $\{R_{n}\} \subset (0, \infty)$ with $\lim_{n \to \infty} R_{n} = \infty$ such that
\[
\lim_{n \to \infty} R_{n} \int_{\p B_{R_{n}}} \left\{ |\nabla \Phi|^{2} + \ell 
\Phi^{2} + \Phi^{2} |\Psi| \right\} \, dS = 0. 
\]
Suppose to the contrary that 
\[
\liminf_{R \to \infty} R \int_{\p B_{R}} \left\{ |\nabla \Phi|^{2} + \ell \Phi^{2} + \Phi^{2} |\Psi| \right\} \, dS = \alpha > 0. 
\]
Then, there exists $R_{1} > 0$ such that 
\[
\int_{\p B_{R}} \left\{ |\nabla \Phi|^{2} + \ell 
\Phi^{2} + \Phi^{2} |\Psi| \right\} \, dS \geq 
\frac{\alpha}{2 R} \quad \text{for all } R > R_{1}. 
\]
This yields
\[
\int_{\R^{4}} \left\{ |\nabla \Phi|^{2} + 
\ell \Phi^{2} + \Phi^{2} |\Psi| \right\} \, dx
= \int_{0}^{\infty} \int_{\p B_{R}} \left\{ |\nabla \Phi|^{2} + \ell \Phi^{2} 
+ \Phi^{2} |\Psi| \right\} \, dS \, dR
\geq \int_{R_{1}}^{\infty} \frac{\alpha}{2 R} \, dR = \infty, 
\]
which contradicts \eqref{eq1-poho}. Taking the limit along the sequence $\{R_n\}$, we conclude that \eqref{eq2-poho} holds. 

Identity \eqref{eq3-poho} can be obtained similarly. This completes the proof.
\end{proof}
\begin{lemma}\label{lem1-unique}
Assume that $(\Phi, \Psi) \in H_{\ell, k}$ 
is a weak solution to \eqref{asp}. 
Then, the following identity hold: 
    \begin{equation} \label{eq1-unique}
    \ell \|\Phi\|_{L^{2}}^{2}
    = \|\nabla \Psi
    \|_{L^{2}}^{2}.  
    \end{equation}
\end{lemma}
\begin{proof}
It follows from \eqref{eq2-poho} and \eqref{eq3-poho} that 
    \[
    \|\nabla \Phi
    \|_{L^{2}}^{2} 
    + 2 \ell 
    \|\Phi
    \|_{L^{2}}^{2}
    +\|\nabla  \Psi
    \|_{L^{2}}^{2} 
    + 2 k \|\Psi
    \|_{L^{2}}^2
    - 4 \int_{\R^{4}} 
    |\Phi|^{2} 
    \Psi \, dx  
    = 0.
    \]
Next, multiplying 
the first equation of \eqref{asp} by 
$\Phi$, 
the second 
equation of \eqref{asp} by 
$\Psi$, 
and integrating the resulting equation, respectively, one has 
    \begin{align}
    & \|\nabla \Phi
    \|_{L^{2}}^{2} 
    + \ell 
    \|\Phi
    \|_{L^{2}}^{2}
    - 2 \int_{\R^{4}} 
    |\Phi|^{2} 
    \Psi \, dx 
    = 0, \notag \\
    &  
    \|\nabla  \Psi
    \|_{L^{2}}^{2} 
    +  k \|\Psi
    \|_{L^{2}}^2
    - \int_{\R^{4}} 
    |\Phi|^{2} 
    \Psi \, dx 
    = 0.
    \label{eq38:large}
    \end{align}
From these identities, we obtain     
\[
   \ell
    \|\Phi
    \|_{L^{2}}^{2}
    + k \|\Psi
    \|_{L^{2}}^2
    - \int_{\R^{4}} 
    |\Phi|^{2} 
    \Psi \, dx  
    = 0.
    \]
Using this and \eqref{eq38:large} again, 
we obtain \eqref{eq1-unique}. 
\end{proof}

\begin{thank}
M.H. was supported by JSPS KAKENHI 
Grant Number JP25K17285 and JST CREST
Grant Number JPMJCR24Q6.
H.K. was supported by JSPS KAKENHI 
Grant Number JP25H01453 and JP25K07089. 
K. T. was supported by JSPS KAKENHI Grant Number JP25KJ0325.
\end{thank}

 \begin{availability}
Data sharing not applicable to this article as no datasets were generated or analyzed during the current study.
\end{availability}

\begin{interest}
The authors declare that they have no conflict of interest.    
\end{interest}

\begin{ai}
During the preparation of this manuscript,
the authors used ChatGPT and Google Gemini for language
editing and to support the development and
verification of mathematical arguments.
Specifically, ChatGPT helped identify a potential
issue in the proof of Theorem \ref{thm-conv2}
and suggested a strategy to address it,
which the authors subsequently incorporated
into the proof.
The authors have carefully reviewed and verified
all content, including the mathematical statements
and proofs, and take full responsibility
for the final manuscript.
\end{ai}

\bibliographystyle{plain}
\bibliography{References} 

\noindent 
Masaru Hamano
\\[6pt]
Center for Science Adventure and Collaborative Research Advancement, 
Graduate School of
Science, 
\\[6pt]
Kyoto University, 
\\[6pt]
Kitashirakawa Oiwakecho, Sakyo-ku, 
Kyoto 606-8502, JAPAN.
\\[6pt]
E-mail: hamano.masaru.3c@kyoto-u.ac.jp

\vspace{0.5cm}

\noindent
Hiroaki Kikuchi
\\[6pt]
Department of Mathematics,
\\[6pt]
Tsuda University,
\\[6pt]
2-1-1 Tsuda-machi, Kodaira-shi, Tokyo 187-8577, JAPAN
\\[6pt]
E-mail: hiroaki@tsuda.ac.jp

\vspace{0.5cm}

\noindent 
Kenta Tomioka
\\[6pt] 
Institute for Mathematics and Computer Science, 
\\[6pt]
Tsuda University, 
\\[6pt]
2-1-1 Tsuda-machi, Kodaira-shi, Tokyo 187-8577, JAPAN
\\[6pt]
kenta.tomioka@gm.tsuda.ac.jp
\end{document}